\documentclass[english,11pt,a4paper]{article}

\usepackage[english]{babel}
\usepackage{enumitem} % ermoeglicht manuelle Aenderung von itemize und enumerate
\usepackage{wrapfig}
\usepackage[final]{graphicx}
\usepackage{epsfig}
\usepackage{textcomp}
\usepackage{amsmath}
\usepackage{amssymb}
\usepackage{amsfonts}
\usepackage{mathtools} 
\usepackage{mathbbol} % fuer Befehl \Eins
\usepackage{psfrag}
\usepackage{latexsym}
\usepackage{eurosym}
\usepackage[pdftex,dvipsnames]{xcolor} % coloured text etc.

\usepackage{tikz}
\usepackage{stmaryrd} % fuer dreifache Normstriche mit \interleave
\usepackage[normalem]{ulem}
\usepackage{verbatim} % \begin{comment}

\usepackage{mathabx} %\widecheck
\usepackage{dsfont} % 1|
\usepackage{lmodern}	% necessary for font sharpness 
\usepackage[numbers,sort&compress]{natbib}
\numberwithin{equation}{section} % Gleichungsnummern im Format ``(Section.Nummer)'' statt ''(Nummer)''

\usepackage{hyphenat}

\usepackage{tikz}
\usepackage{pgfplots}
\usetikzlibrary{matrix, backgrounds}

\tikzset{ 
table/.style={
  matrix of nodes,
  row sep=-\pgflinewidth,
  column sep=-\pgflinewidth,
  nodes={rectangle,text width=3em,align=center},
  text depth=1.25ex,
  text height=2.5ex,
  nodes in empty cells
},
}

\usepackage{hyperref}

\newcommand{\epsi}{\varepsilon}
\newcommand{\IC}{\mathbb{C}} 
\newcommand{\IN}{\mathbb{N}}

\newcommand{\IR}{\mathbb{R}}

\newcommand{\tend}{t_{\mbox{\tiny end}}}
\newcommand{\tstar}{t_\star}

\newcommand{\diag}{\mbox{diag}}

\newcommand{\hs}[1]{\hspace*{#1 mm}}

\newcommand{\Ord}[1]{\mathcal{O} \! \left(#1\right)}

\newcommand{\F}{\mathcal{F}}
\newcommand{\G}{\mathcal{G}}

\newcommand{\J}{\mathcal{J}}
\renewcommand{\L}{\mathcal{L}}
\renewcommand{\P}{\mathcal{P}}
\renewcommand{\S}{\mathcal{S}}
\newcommand{\T}{\mathcal{T}}
\newcommand{\TT}{S}
\newcommand{\FF}{F} % {\mathbf{F}}

\newcommand{\dd}{\mathrm{d}}
\newcommand{\ee}{\mathrm{e}}
\newcommand{\ii}{\mathrm{i}}

\newcommand{\w}{\omega}
\newcommand{\Chat}{\widehat{C}}
\newcommand{\uhat}{\widehat{u}}
\newcommand{\vhat}{\widehat{v}}
\newcommand{\phat}{\widehat{p}}
\newcommand{\qhat}{\widehat{q}}
\newcommand{\Uhat}{\widehat{U}}
\newcommand{\uu}{u} % \newcommand{\uu}{\mathbf{u}}
\newcommand{\uutilde}{\widetilde{u}} % \newcommand{\uutilde}{\mathbf{\widetilde{u}}}

\newcommand{\vtilde}{\widetilde{v}} % \newcommand{\uutilde}{\mathbf{\widetilde{u}}}
\newcommand{\Pperp}{P^\perp}

\newcommand{\pt}{\partial_{t}}
\newcommand{\px}{\partial_{x}}

\newcommand{\intl}{\int\limits}

\newcommand{\weg}[1]{}
\newcommand{\mytodo}[1]{{\fontfamily{pag}\selectfont\color{red}#1 \color{black} }}
\newcommand{\mynote}[1]{{\fontfamily{pag}\selectfont\color{blue}#1 \color{black} }}
\newcommand{\myfootnote}[1]{\footnote{\color{blue} #1} \color{black}}
\renewcommand{\mytodo}[1]{} % Notizen verschwinden
\renewcommand{\mynote}[1]{} % Notizen verschwinden
\renewcommand{\myfootnote}[1]{} % Fussnoten verschwinden
\newcommand{\proof}{\noindent\textbf{Proof. }}
\newcommand{\prooftext}[1]{\noindent\textbf{Proof#1}}
\newcommand{\proofof}[1]{\noindent\textbf{Proof of #1. }}
\def\qed {\hfill \rule{0,2cm}{0,2cm} \bigskip \\}

\newcommand{\il}{{\interleave}}

\newtheorem{Theorem}{Theorem}[section]
\newtheorem{Proposition}[Theorem]{Proposition}
\newtheorem{Lemma}[Theorem]{Lemma}
\newtheorem{Corollary}[Theorem]{Corollary}
\newtheorem{Assumption}[Theorem]{Assumption}

\newtheorem{Remark}[Theorem]{Remark}

\newcounter{fig}

\begin{document}

\title{Improved error bounds for approximations of high-frequency wave propagation in nonlinear dispersive media\thanks{Funded by the Deutsche Forschungsgemeinschaft (DFG, German Research Foundation) -- Project--ID 258734477 -- SFB 1173.}}
\author{Julian Baumstark\thanks{Karlsruher Institut f\"ur Technologie, Fakult\"at f\"ur Mathematik, Institut f\"ur Angewandte und Numerische Mathematik, Englerstr.~2, D-76131 Karlsruhe,
\texttt{julian.baumstark@gmx.de}, \texttt{tobias.jahnke@kit.edu}} \and Tobias Jahnke\footnotemark[2]}
\date{\today}
\maketitle

\abstract{
\noindent
High-frequency wave propagation is often modelled by
nonlinear Friedrichs systems where both the differential equation 
and the initial data contain the inverse of a small parameter $\epsi$,
which causes oscillations with wavelengths proportional to $\epsi$ in time and space.
A prominent example is the Maxwell--Lorentz system, which is a well-established model for the propagation of light in nonlinear media.
In diffractive optics, such problems have to be solved on long time 
intervals with length proportional to $1/\epsi$. 
Approximating the solution of such a problem numerically with a standard method is hopeless, because traditional methods require an extremely fine resolution in time and space, which entails unacceptable computational costs.
A possible alternative is to replace the original problem by a new system of PDEs which 
is more suitable for numerical computations but still yields a sufficiently accurate approximation.
Such models are often based on the \emph{slowly varying envelope approximation} or generalizations thereof.
Results in the literature state that the error of the slowly varying envelope approximation is of $\Ord{\epsi}$. In this work, however, we prove that the error is even proportional
to $\epsi^2$, which is a substantial improvement, and which explains the error behavior observed in numerical experiments.
For a higher-order generalization of the slowly varying envelope approximation we improve the error bound 
from $\Ord{\epsi^2}$ to $\Ord{\epsi^3}$.
Both proofs are based on a careful analysis of 
the nonlinear interaction between oscillatory and non-oscillatory error terms,
and on \textit{a priori} bounds for certain ``parts'' of the approximations which are defined by suitable projections. 
As an important technical tool we use an advantageous transformation of the coefficient functions which appear in the approximations.
}


% --------------------------------------------------------------------------------
\section{Introduction}\label{Sec.introduction}
% --------------------------------------------------------------------------------

High-frequency wave propagation in nonlinear, dispersive media can be modeled by 
Friedrichs systems of the form
\begin{subequations}
\label{PDE.uu}
\begin{align}
\label{PDE.uu.a}
 \pt \uu +A(\partial)\uu +\frac{1}{\epsi} E\uu &= \epsi T(\uu,\uu,\uu),
 & & t\in(0,\tend/\epsi], \; x\in\IR^d,
 \\
\label{PDE.uu.b}
 \uu(0,x) &= p(x)\ee^{\ii (\kappa \cdot x)/\epsi} + c.c.,
\end{align}
\end{subequations}
with a trilinear nonlinearity $T: \IR^n \times \IR^n \times \IR^n \rightarrow \IR^n $ and a differential operator 
\begin{align}
\label{Def.A.partial}
A(\partial) = \sum_{\ell=1}^d A_\ell \partial_\ell
\end{align}
($d,n\in\IN$).
We assume that the matrices $A_1, \ldots, A_d \in \IR^{n \times n}$ in \eqref{Def.A.partial} are symmetric,
and that $ E\in\IR^{n\times n} $ in \eqref{PDE.uu.a} is skew-symmetric.
In the initial data a smooth and localized envelope function $p: \IR^d \rightarrow \IR^n$ is multiplied by a phase with a given wave vector $\kappa \in \IR^d\setminus\{0\}$. 
Here and below, ``$X$ + c.c.'' means $X+\overline{X}$, where $\overline{X}$ is the complex conjugate of 
$X$.
An important example in this class of problems is the Maxwell--Lorentz system,
which is a classical model for the propagation of light in a Kerr medium; 
cf.~\cite{donnat-rauch:97a,colin-lannes:09,lannes:11,lannes:98,donnat-rauch:97b,joly-metivier-rauch:96,colin-gallica-laurioux:05}.

% Both the PDE \eqref{PDE.uu.a} and the initial data in \eqref{PDE.uu.b} involve a small positive parameter
% $ \epsi \ll 1 $, which accounts for different scales in time and space. 

The PDE \eqref{PDE.uu.a}, the initial data in \eqref{PDE.uu.b}, and the time interval 
involve a small positive parameter $ \epsi \ll 1 $.
Although the nonlinearity in \eqref{PDE.uu.a} is multiplied by $\epsi$, the problem 
\eqref{PDE.uu} is strongly nonlinear, because the length of the time interval is proportional to $\epsi^{-1}$.
In fact, by rescaling $\tau=\epsi t$ and $ w(\tau,x)=u(t,x) $, we could convert \eqref{PDE.uu.a}
into the equivalent form
\begin{align*}
 \partial_\tau w + \frac{1}{\epsi}A(\partial)w +\frac{1}{\epsi^2} Ew &= T(w,w,w),
 & & \tau\in(0,\tend], \; x\in\IR^d,
\end{align*}
where the nonlinear term and the time interval do not depend on $\epsi$ anymore.
However, we will consider the original version \eqref{PDE.uu}, which is the representation considered, e.g., 
in \cite{colin-lannes:09,lannes:11,baumstark:22,baumstark-jahnke-2023,baumstark-jahnke-lubich-2024}.

The small parameter $ \epsi $ accounts for different scales in time and space.
The terms $ \ee^{\pm \ii (\kappa \cdot x)/\epsi} $ in the initial data cause spatial oscillations with wavelength of $\Ord{\epsi}$, whereas $p$ changes on a scale of $\Ord{1}$, roughly speaking.
As a consequence, the solution $u(t,x)$ is a wave packet with a high-frequency carrier wave modulated by a smooth envelope.
Concerning the evolution in time, the initial value problem \eqref{PDE.uu} is scaled in such a way that nonlinear and diffractive effects appear on long time intervals of length $\tend/\epsi$ for some $\tend>0$, whereas the envelope of the wave packet propagates with speed $\Ord{1}$.  
The solution itself, however, evolves on a third scale, because the linear part $A(\partial)\uu +\frac{1}{\epsi} E\uu$ of the PDE causes rapid oscillations in time with wavelength of $\Ord{\epsi}$.
Because of the highly oscillatory nature and the long time interval, an attempt to approximate the vector-valued solution $\uu: [0,\tend/\epsi] \times \IR^d \rightarrow \IR^n$ of \eqref{PDE.uu} numerically with a traditional method is bound to fail, because the time and space discretizations would require an extremely fine resolution and hence an impracticable runtime.

A feasible approach is to replace \eqref{PDE.uu} by a different model which can be solved numerically with significantly less efforts and at the same time provides a decent approximation to $\uu$.
Such models are often based on the \emph{slowly varying envelope approximation} (SVEA) or generalizations thereof, which are derived as follows.
For every $\beta \in \IR^d$ the matrix
\begin{align*}
% \label{Def.Ak}
 A(\beta)&=\sum_{\ell=1}^d \beta_\ell A_\ell \in \IR^{n \times n}
\end{align*}
is symmetric, and
\begin{align}
 \label{Def.L}
\L(\alpha,\beta) &= -\alpha I + A(\beta)-\ii E \in \IC^{n\times n}
\end{align}
is Hermitian for all $\alpha \in \IR$ and $\beta \in \IR^d$.
Let $\kappa \in \IR^d\setminus\{0\}$ be the (given) wave vector which appears in 
\eqref{PDE.uu.b}, and 
\begin{align}
\label{Def.w}
\text{let $\w=\w(\kappa)$ be an eigenvalue of $A(\kappa)-\ii E$. }
\end{align}
Then, $\L(\w,\kappa)$ has a non-trivial kernel, and the pair $(\w,\kappa)$ 
is said to fulfill the dispersion relation. We assume the following.

\pagebreak

\begin{Assumption}\label{Ass:polarization}
\text{ } \nopagebreak
\begin{enumerate}[label=(\roman*)]
\item \label{Ass:polarization.item.i}
The kernel of $\L(\w,\kappa)$ is one-dimensional. 

\item \label{Ass:polarization.item.ii}
The function $p$ in \eqref{PDE.uu.b} has the structure
\begin{align}
\label{p.p0.p1}
 p=p_0+\epsi p_1  \qquad 
 \text{with } p_0(x) \in \text{ker}\big(\L(\w,\kappa)\big) \; a.e.
\end{align}
and $ p_0, p_1 \in L^\infty(\IR^d,\IC^n) $.

\end{enumerate}
\end{Assumption}
Assumption~\ref{Ass:polarization.item.i} is only made in order to keep the notation simple;
cf.~Remark~\ref{Remark:eigenspace.one-dimensional} below.
Assumption~\ref{Ass:polarization.item.ii} is a polarization condition, which was also imposed in a similar way  
in \cite[Theorem 1]{colin-lannes:09}, \cite[Theorem 2.15]{lannes:11}, \cite{baumstark-jahnke-2023}, and other works.

\bigskip

\noindent
As in \cite{baumstark-jahnke-2023} we seek an approximation of the form
\begin{align}
\label{Ansatz}
 \uu(t,x) &\approx \uutilde^{(m)}(t,x) = 
 \sum_{j\in \J^{(m)}}
% \ee^{\ii j(\kappa \cdot x - \w t)/\epsi} u_j(t,x),
\qquad u_{-j} = \overline{u_j}
\end{align}
for $ \J^{(m)} = \{\pm 1, \pm 3, \ldots, \pm m\}$, where $m\in\IN$ is an odd integer.
If we substitute \eqref{Ansatz} into \eqref{PDE.uu}, then the trilinear nonlinearity generates \emph{higher harmonics}, i.e.\ terms with prefactor $\ee^{\ii j(\kappa \cdot x -\w t)/\epsi}$
for $|j|>m$. These terms appear only on the right-hand side of \eqref{PDE.uu.a}, because all terms on the 
left-hand side are linear.
Ignoring higher harmonics and then comparing terms on both sides yields the PDE system
\begin{subequations}
\label{PDE.mfe}
\begin{align}
\label{PDE.mfe.a}
\pt u_j + \frac{\ii}{\epsi}\L(j \w, j \kappa)u_j + A(\partial)u_j
= \epsi \sum_{j_1+j_2+j_3=j} T(u_{j_1},u_{j_2},u_{j_3}) &
\\
\notag
\text{for } j\in \J_+^{(m)}=\J^{(m)} \cap\IN, \quad t\in (0,\tend/\epsi], \quad x\in\IR^d &,
\end{align}
with initial conditions
\begin{align}
\label{PDE.mfe.b}
u_1(0,\cdot) = p, \qquad u_j(0,\cdot) = 0 \text{ for } j\in\J_+^{(m)}\setminus\{1\}.
\end{align}
\end{subequations}
The sum on the right-hand side of \eqref{PDE.mfe.a} is to be taken over the set 
\begin{align*}
 \Big\{J=(j_1,j_2,j_3) \in (\J^{(m)})^3 : \#J:= j_1+j_2+j_3=j \Big\},
\end{align*}
and $T$ is now the trilinear extension of the real nonlinearity from \eqref{PDE.uu.a}
to $T:\IC^n \times \IC^n \times \IC^n \to \IC^n $.
It is sufficient to consider positive $j\in \J_+^{(m)}=\J^{(m)} \cap\IN$ instead of $j\in\J^{(m)}$ in \eqref{PDE.mfe.a}, because the $u_j$ with negative index $j$ 
are obtained from the condition $u_{-j} = \overline{u_j}$.
For $m=1$ and $\J^{(1)}=\{-1,1\}$, we obtain the SVEA

\begin{align}
\label{Ansatz.SVEA}
 \uu(t,x) \approx \uutilde^{(1)}(t,x) &= \ee^{\ii (\kappa \cdot x - \w t)/\epsi} u_1(t,x) + c.c.,
\end{align}
with $u_1$ being the solution of
\begin{subequations}
\label{SVEA}
\begin{align}
\label{SVEA.a}
 \pt u_1 + \frac{\ii}{\epsi}\L(\w,\kappa)u_1 + A(\partial)u_1
&= \epsi \sum_{j_1+j_2+j_3=1} T(u_{j_1},u_{j_2},u_{j_3})
\\
\notag
&= \epsi \Big( T(u_1,u_1,u_{-1})+T(u_1,u_{-1},u_1)+T(u_{-1},u_1,u_1) \Big),
 \\
\label{SVEA.b}
u_1(0,\cdot) &= p
\end{align}
\end{subequations}
as a special case of \eqref{Ansatz} and \eqref{PDE.mfe}.
Note that the initial data in \eqref{PDE.mfe.b} are smooth, non-oscillatory functions, in contrast to \eqref{PDE.uu.b}. Hence, solutions to \eqref{PDE.mfe} can be discretized in space on a $\epsi$-independent grid, which is a significant advantage over \eqref{PDE.uu}. However, typical solutions of \eqref{PDE.mfe} do still oscillate in time due to the term $ \tfrac{\ii}{\epsi}\L(j\w,j\kappa)u_j $ in \eqref{PDE.mfe.a}.

For the error of the SVEA \eqref{Ansatz.SVEA}--\eqref{SVEA} the bound
\begin{align}
\label{error.bound.svea.cl}
\sup_{t\in[0,\tend/\epsi]} \| \uu(t,\cdot) - \uutilde^{(1)}(t,\cdot) \|_{L^\infty(\IR^d,\IC^n)} \leq C\epsi
\end{align}
was shown in \cite[Section 2.2]{colin-lannes:09}. 
Under additional assumptions, one can replace the PDE \eqref{SVEA.a} by a nonlinear Schr\"odinger equation without spoiling the error bound \eqref{error.bound.svea.cl}; 
cf.~\cite[Corollary 2]{colin-lannes:09} and also \cite{colin:02,lannes:11,donnat-joly-metivier-rauch:96,joly-metivier-rauch:98,kirrmann-schneider-mielke:92,schneider-uecker:17}. 
This nonlinear Schr\"odinger equation has the advantage that it does not involve $\epsi$ at all
when considered in a co-moving coordinate system, and that it only has to be solved on the $\epsi$-independent time interval $[0,\tend]$. Hence, standard numerical methods can be used to solve 
the nonlinear Schr\"odinger equation numerically, which then yields an $\Ord{\epsi}$-approximation to 
$ \uutilde^{(1)} $ and, via \eqref{error.bound.svea.cl}, to the solution of \eqref{PDE.uu}.

In this paper, we consider the situation where an error of $\Ord{\epsi}$ is yet too large.
For the approximation 
\begin{align*}
  \uu(t,x) \approx \uutilde^{(3)}(t,x) &= 
  \left( \ee^{\ii (\kappa \cdot x - \w t)/\epsi} u_1(t,x) + 
  \ee^{3\ii (\kappa \cdot x - \w t)/\epsi} u_3(t,x) \right) + c.c.
\end{align*}
we have already shown the error bound
\begin{align}
\label{error.bound.BJ.3}
\sup_{t\in[0,\tstar/\epsi]} \| \uu(t,\cdot) - \uutilde^{(3)}(t,\cdot) \|_{L^\infty(\IR^d,\IC^n)} \leq C\epsi^2
\end{align}
for some $ \tstar \in (0,\tend]$ in \cite{baumstark-jahnke-2023}.
The proof is rather long and technical because of the complicated self-interaction of the oscillatory solution via the nonlinearity. Moreover, the approximation $\uutilde^{(3)}$ is more complicated than \eqref{Ansatz.SVEA} because of the additional coefficient function $u_3$.
Numerical experiments show, however, that the estimates \eqref{error.bound.svea.cl} and \eqref{error.bound.BJ.3} \emph{are both not optimal}; see Sections~\ref{Subsec:NumEx.01} and \ref{Subsec:NumEx.02} below.
In this work, we will prove the improved error bounds
\begin{align}
\label{error.bound.BJ.1.improved}
\sup_{t\in[0,\tend/\epsi]} \| \uu(t,\cdot) - \uutilde^{(1)}(t,\cdot) \|_{L^\infty(\IR^d,\IC^n)} 
&\leq C\epsi^2,
\\
\label{error.bound.BJ.3.improved}
\sup_{t\in[0,\tstar/\epsi]} \| \uu(t,\cdot) - \uutilde^{(3)}(t,\cdot) \|_{L^\infty(\IR^d,\IC^n)} 
&\leq C\epsi^3.
\end{align}
The result \eqref{error.bound.BJ.1.improved} explains the error behavior which appears in numerical examples
where a reference solution can be computed.
Moreover, this inequality shows that the SVEA yields a \emph{significantly higher accuracy than the classical
nonlinear Schr\"odinger approximation}, which has an error of $\Ord{\epsi}$. 
This fact was apparently not known until now.
The second error bound \eqref{error.bound.BJ.3.improved} states that in applications where an error of
$\Ord{\epsi^2}$ is still too large, the refined approximation $ \uutilde^{(3)} $ offers the possibility of reducing the error down to $\Ord{\epsi^3}$ at the cost of higher computational work.

In \cite{rauch-book:12,donnat-rauch:97b,joly-metivier-rauch:00,joly-metivier-rauch:93} and other contributions, asymptotic expansions of solutions to problems similar to \eqref{PDE.uu} have been analyzed in the regime of geometric optics, i.e.\ for time intervals of length $\Ord{1}$.
This differs from the regime of diffractive geometric optics, where the PDE system has to be solved on time intervals of length $\Ord{1/\epsi}$, which is the situation we consider here.
Approximations in diffractive geometric optics have been constructed in 
\cite{donnat-joly-metivier-rauch:96} 
and \cite{joly-metivier-rauch:98} for semilinear and quasilinear systems with a more general nonlinearity, but with $\epsi E$ instead of $E/\epsi$ in \cite{donnat-joly-metivier-rauch:96} and with $E=0$
in \cite{joly-metivier-rauch:98}. 
Quasilinear systems with dispersion and dispersive problems with bilinear nonlinearity are approximated in 
\cite{lannes:98} and \cite{colin:02}, respectively, 
but without an explicit rate of convergence.
The book \cite{schneider-uecker:17} provides an extensive analysis of the approximation of PDEs by nonlinear Schr\"odinger equations and other modulation equations.

In \cite{baumstark-jahnke-lubich-2024} we have constructed modulated Fourier expansions for \eqref{PDE.uu}
with \emph{nonlinear} polarization of the initial data. This approach is likewise based on the ansatz \eqref{Ansatz} and \eqref{PDE.mfe},
but the nonlinear polarization considered in \cite{baumstark-jahnke-lubich-2024} means that
$p_1$ depends on $p_0$ in \eqref{p.p0.p1}, which excludes, e.g., the case $ p_1=0$.
In the present work, $p_1$ and $p_0$ are completely independent.

In this paper, we consider wave packets where the wavelength of the oscillations is much shorter
than the scale on which the envelope varies.
This assumption excludes short or chirped pulses. Since it is known that the nonlinear Schr\"odinger approximation does not provide a reasonable approximation for such pulses, many improved models have been proposed and analyzed, e.g., in \cite{barrailh-lannes:02,colin-lannes:09,lannes:11,colin-gallica-laurioux:05,alterman-rauch:00,alterman-rauch:03,chung-jones-schaefer-wayne:05}. 

In Section~\ref{Sec.analytical.setting} we specify the analytical framework, we review results on local well-posedness of \eqref{PDE.uu} and \eqref{PDE.mfe}, and we introduce a transformation of the coefficient functions $u_j$ 
which was already employed in \cite{baumstark-jahnke-2023}.
The proofs of the error bounds \eqref{error.bound.BJ.1.improved} and \eqref{error.bound.BJ.3.improved} rely on the fact 
that for a certain projection $\P_\epsi$ the Fourier transform $\uhat_1$ of the coefficient function $u_1$ can be decomposed into an essentially non-oscillatory part $\P_\epsi \uhat_1$ and an oscillatory but ``small'' part 
$(I-\P_\epsi) \uhat_1$. For the SVEA (i.e. for $m=1$) we compile the corresponding results in
Section~\ref{Sec.error.SVEA.preparations}. Then, in Section~\ref{Sec.error.bound.SVEA}, we prove the error bound \eqref{error.bound.BJ.1.improved} for the SVEA, and we corroborate this result by a numerical experiment. 
In Section~\ref{Sec.error.bound.approximation.m3}, we turn to the case $m=3$. We show the error bound \eqref{error.bound.BJ.3.improved} 
and we give reasons why we observe an even better rate of convergence
in a numerical example with a one-dimensional Klein--Gordon system.

\paragraph{Notation.}
Throughout the text, $v \cdot w = v^* w $ is the Euclidean scalar product of vectors $v,w \in \IC^n$,
and $ |v|_q $ is the $q-$norm of $v$.
The identity matrix and the identity operator are both denoted by $I$.
For space- and time-dependent functions $f=f(t,x)$ we will often denote the mapping $ x \mapsto f(t,x) $ by $f(t)$ instead of $f(t,\cdot)$.
Likewise, we will omit the second argument of the Fourier transform $\widehat{f}(t,k)$ of such a function.
From now on, we will use the short-hand notation $L^1$ and $L^\infty$ for 
$ L^1(\IR^d,\IC^n) $ and $ L^\infty(\IR^d,\IC^n), $ respectively.
The symbol $\ii=\sqrt{-1}$ is the imaginary unit, whereas $i$ appears as an index in a few formulas.

% --------------------------------------------------------------------------------
\section{Analytical setting}\label{Sec.analytical.setting}
% --------------------------------------------------------------------------------

% In this section we compile a number of analytical tools, concepts and results 
% for later use in our analysis.

\paragraph{Wiener algebra and evolution equations in Fourier space.}

As in \cite{baumstark:22,colin-lannes:09,baumstark-jahnke-2023,baumstark-jahnke-lubich-2024,lannes:11} 
we will analyze the accuracy in the Wiener algebra  
\begin{align}
\label{Def.W}
W &= \left\{f \in \big(\S'(\IR^d)\big)^n: \widehat{f}\in L^1
\right\}, &
\|f\|_W &= \|\widehat{f}\|_{L^1} = \intl_{\IR^d} |\widehat{f}(k)|_2 \; \dd k
\end{align}
of vector-valued functions.
Here and below, $\widehat{f} = \F f $ denotes the Fourier transform
\begin{align*}
(\F f)(k) := (2\pi)^{-d/2}\intl_{\IR^d} f(x)\ee^{-\ii k\cdot x} \dd x
\end{align*}
 of $f$. For $s\in\IN_0$, we define
\begin{align*}
W^s &= \{f \in W: \partial^\alpha f \in W \text{ for all }
\alpha \in \IN_0^d, |\alpha|_1 \leq s\}, \\
\|f\|_{W^s} &= \sum_{|\alpha|_1 \leq s} \|\partial^\alpha f\|_W.
\end{align*}
It is well-known that $W^s$ is a Banach algebra with continuous embedding 
$ W \hookrightarrow L^\infty$, 
cf.~\cite[Proposition 1]{colin-lannes:09} and \cite[Proposition 3.2]{lannes:11}.

In order to work in the Wiener algebra, we apply the Fourier transform to the PDE system \eqref{PDE.mfe.a}. 
This yields
\begin{subequations}
\label{PDE.mfe.Fourier}
\begin{align}
\label{PDE.mfe.Fourier.a}
\pt \uhat_j(t,k) + \frac{\ii}{\epsi}\L_j(\epsi k)\uhat_j(t,k)
=
\epsi \sum_{\#J=j} \T\big(\uhat_{j_1},\uhat_{j_2},\uhat_{j_3}\big)(t,k), &
\\
\notag
 j\in \J_+^{(m)}, \quad t\in(0,\tend/\epsi], \quad k \in \IR^d &
\end{align}
with initial conditions
\begin{align}
\label{PDE.mfe.Fourier.b}
\uhat_1(0,\cdot) = \phat, \qquad \qquad 
\uhat_j(0,\cdot) = 0 \quad \text{for } j\in\J_+^{(m)}\setminus\{1\}
\end{align}
\end{subequations}
and the notation
\begin{align}
\label{Def.L.j}
\L_j(\theta) &= \L(j\w, j\kappa + \theta) = \L_j(0) + A(\theta),
\qquad j\in \J_+^{(m)},
\\
\notag
% \label{Def.T.Fourier}
\T\left(\uhat_{j_1},\uhat_{j_2},\uhat_{j_3}\right)(k)
&=
\F\Big(T(u_{j_1},u_{j_2},u_{j_3})\Big)(k) 
\\
\notag
&= 
(2\pi)^{-d}
\intl_{\IR^d} \intl_{\IR^d} 
T(\uhat_{j_1}(k^{(1)}),\uhat_{j_2}(k^{(2)}),\uhat_{j_3}(k-k^{(1)}-k^{(2)}))
\; \dd k^{(2)} \; \dd k^{(1)},
\end{align}
cf.~\cite[Section 2.2]{baumstark-jahnke-2023}.
In \eqref{Def.L.j} we have used that by definition the mapping $\beta \mapsto A(\beta)$ is linear.
With the shorthand notation
\begin{subequations}
\label{Shorthand.notation.K}
\begin{align}
K=\big(k^{(1)},k^{(2)},k^{(3)}\big)\in \IR^d \times \IR^d \times \IR^d,
\qquad \qquad
\#K:=k^{(1)}+k^{(2)}+k^{(3)}\in \IR^d,
\end{align}
and
\begin{align}
& \intl_{\#K = k} T(\uhat_{j_1}(k^{(1)}),\uhat_{j_2}(k^{(2)}),\uhat_{j_3}(k^{(3)}))
\; \dd K
\\
\notag
&=
\intl_{\IR^d} \intl_{\IR^d} 
T(\uhat_{j_1}(k^{(1)}),\uhat_{j_2}(k^{(2)}),\uhat_{j_3}(k-k^{(1)}-k^{(2)}))
\; \dd k^{(2)} \; \dd k^{(1)},
\end{align}
\end{subequations}
the Fourier transform of the nonlinearity can be expressed as
\begin{align}
\label{Def.T.Fourier}
\T\left(\uhat_{j_1},\uhat_{j_2},\uhat_{j_3}\right)(k)
&= 
(2\pi)^{-d} \intl_{\#K = k} 
T(\uhat_{j_1}(k^{(1)}),\uhat_{j_2}(k^{(2)}),\uhat_{j_3}(k^{(3)}))
\; \dd K.
\end{align}
Later we will often use that
\begin{align}
\big\| \T(\widehat{f}_1,\widehat{f}_2,\widehat{f}_3) \big\|_{L^1}
\label{Bound.T.Fourier}
\leq 
C_\T \| \widehat{f}_1 \|_{L^1} \| \widehat{f}_2 \|_{L^1} \| \widehat{f}_3 \|_{L^1}
\end{align}
with a constant $C_\T$ which depends on $\T$ and on $n$.
Via trilinearity, we obtain that
\begin{align}
\label{Bound.T.difference.Fourier}
\big\| \T(\widehat{f}_1,\widehat{f}_2,\widehat{f}_3) - \T(\widehat{g}_1,\widehat{g}_2,\widehat{g}_3) \big\|_{L^1}
&\leq 
C_\T \| \widehat{f}_1 - \widehat{g}_1 \|_{L^1} \| \widehat{f}_2 \|_{L^1} \| \widehat{f}_3 \|_{L^1}
\\
\notag
&\quad
+
C_\T \| \widehat{g}_1 \|_{L^1} \| \widehat{f}_2 - \widehat{g}_2 \|_{L^1} \| \widehat{f}_3 \|_{L^1}
\\
\notag
&\quad
+
C_\T \| \widehat{g}_1 \|_{L^1} \| \widehat{g}_2 \|_{L^1} \| \widehat{f}_3 - \widehat{g}_3 \|_{L^1}.
\end{align}
We set $u_{-j} = \overline{u_j}$ throughout, which implies that 
$\uhat_{-j}(t,k) = \overline{\uhat_j(t,-k)}$.
The system \eqref{PDE.mfe.Fourier.a} can be extended to  
$ j\in \J^{(m)} $ (including negative indices) if we define 
\begin{align}
\label{Def.L.j.negative}
\L_{-j}(\theta)&= -\overline{\L_j(-\theta)}
\qquad \text{for } j\in \J_+^{(m)}.
\end{align}

\paragraph{Local well-posedness.}

The polarization condition (Assumption~\ref{Ass:polarization}\ref{Ass:polarization.item.ii}) 
is not needed to prove existence and uniqueness of solutions 
to the original problem \eqref{PDE.uu} and the PDE system \eqref{PDE.mfe}.
For the sake of consistency, however, we always allow for $\epsi$-dependent initial data of the form
\begin{align}
\label{Def:p}
p = p_0 + \epsi p_1 \qquad \text{with } p_0, p_1 \in W^\sigma
\end{align}
for some $\sigma \in \IN$. The value of $\sigma$ will be specified whenever we refer to \eqref{Def:p}.

\begin{Lemma}[Local well-posedness of \eqref{PDE.uu}]\label{Lemma.wellposedness.uu}
If $ p_0, p_1 \in W $, then there is a $\tend>0$ such that for every $\epsi\in(0,1]$ 
the original problem \eqref{PDE.uu} with $p = p_0 + \epsi p_1$ has a unique mild solution 
$ \uu \in C([0,\tend/\epsi),W) $
which is uniformly bounded, i.e.\ there is a constant $c>0$ such that
\begin{align*}
\sup_{t\in[0,\tend/\epsi]}\| \uu(t) \|_{W} \leq c
\qquad \text{for all } \epsi \in (0,1]. 
\end{align*}
\end{Lemma}
We omit the proof, because Lemma~\ref{Lemma.wellposedness.uu} can be shown with the usual fixed-point argument.
Other proofs for well-posedness of \eqref{PDE.uu} via approximation by the SVEA are given in 
\cite[Theorem 1]{colin-lannes:09} and \cite[Theorem 3.8]{lannes:11}.

\begin{Lemma}[Local well-posedness of \eqref{PDE.mfe}]\label{Lemma.wellposedness}

Let $ m \in\IN $ be an odd integer.

\begin{enumerate}[label=(\roman*)]
\item
\label{Lemma.wellposedness.0}
If \eqref{Def:p} holds with $\sigma=0$ and $ C_{u,0}>\|p_0\|_W + \|p_1\|_W$, then there is a $\tend>0$ such that for every $\epsi\in(0,1]$ the system \eqref{PDE.mfe} has a unique mild solution 
\begin{align*}
\{u_j\}_{j\in\J_+^{(m)}}, \qquad
u_j \in C([0,\tend/\epsi),W)
\end{align*}
which is uniformly bounded, i.e.
\begin{align*}
% \label{Lemma.wellposedness.bound.0}
\sup_{t\in[0,\tend/\epsi]}\| u_j(t) \|_{W} \leq C_{u,0}
\qquad \text{for all } j\in\J_+^{(m)} \text{ and all } \epsi \in (0,1]. 
\end{align*}

\item
\label{Lemma.wellposedness.1}
If \eqref{Def:p} holds with $\sigma=1$, then the mild solution on $[0,\tend/\epsi]$ is a classical solution with
\begin{align*}
% \notag
& u_j \in C^1([0,\tend/\epsi],W) \cap C([0,\tend/\epsi],W^1),
&& j\in\J_+^{(m)},
\\
% \label{Lemma.wellposedness.bound.1}
& \sup_{t\in[0,\tend/\epsi]}\| u_j(t) \|_{W^1} \leq C_{u,1}.
\end{align*}

\item
\label{Lemma.wellposedness.23}
If \eqref{Def:p} holds with $\sigma \in\{2,3\}$, then 
\begin{align}
\notag
u_j \in C^{\sigma-\ell}([0,\tend/\epsi],W^\ell) & \quad \text{for every } \ell = 0, \ldots, \sigma,
&& j\in\J_+^{(m)},
\\
\label{Lemma.wellposedness.bound.23}
\sup_{t\in[0,\tend/\epsi]} \| u_j(t) \|_{W^\sigma} \leq C_{u,\sigma}.
\end{align}
\end{enumerate}
The constants $C_{u,\sigma}$, $\sigma\in\{0,1,2,3\}$, depend on the nonlinearity $T$ and 
on $\tend$, $\|p_0\|_{W^\sigma}$, $\|p_1\|_{W^\sigma}$, but not on $\epsi\in(0,1]$.
\end{Lemma}
For $m=3$ a slightly different version of this result was shown in \cite[Lemma 2.3]{baumstark-jahnke-2023}.
The extension to arbitrary 
odd $m$ is straightforward. Wellposedness of the SVEA ($m=1$) was proven in
\cite[Theorem 1]{colin-lannes:09} and \cite[Theorem 3.8]{lannes:11}.

Although $\tend$ does in general not have the same value in Lemma~\ref{Lemma.wellposedness.uu} and
Lemma~\ref{Lemma.wellposedness}, 
we will henceforth assume that solutions to \eqref{PDE.uu} and \eqref{PDE.mfe} exist on the \emph{same} interval 
$[0,\tend/\epsi]$, as suggested by our notation. 
This is not a restriction as one can always consider the smaller one of the two possibly different intervals.

\paragraph{Eigendecompositions.}

The highly oscillatory behavior of the coefficient functions $ \uhat_j $ originates from the linear part $ \tfrac{\ii}{\epsi}\L_j(\epsi k)\uhat_j(t,k) $ in \eqref{PDE.mfe.Fourier.a}.
It is thus not surprising that the eigendecomposition of 
$\L_j(\theta) = \L(j\w, j\kappa + \theta)$ plays a crucial role in our analysis.
As in \cite[Assumption 2.2]{baumstark-jahnke-2023} we assume the following.

\begin{Assumption}\label{Ass:L.properties}
\text{ }
\begin{enumerate}[label=(\roman*)]
\item \label{Ass:L.properties.item.i}
The matrix $\L(0,\beta) = A(\beta)-\ii E$ has a smooth eigendecomposition:
if $\w_\ell(\beta)$ is an eigenvalue of $\L(0,\beta)$ for some $\ell\in \{1, \ldots, n\}$, then $ \w_\ell \in C^\infty(\IR^d\setminus\{0\}, \IR) $,
and there is a corresponding eigenvector $ \psi_\ell(\beta) $ with 
$ \psi_\ell \in C^\infty(\IR^d\setminus\{0\}, \IC^{n}). $
With no loss of generality, we assume that $ |\psi_\ell(\beta)|_2=1 $ for all $\beta$ and all $\ell=1, \ldots, n$.
The enumeration is chosen in such a way that $\w=\w_1(\kappa)$ in \eqref{Def.w}.

\item \label{Ass:L.properties.item.ii}
Every eigenvalue $\w_\ell(\beta)$ of $\L(0,\beta)$ is globally Lipschitz continuous, i.e. there is a constant $C$ such that
\begin{align*}
| \w_\ell(\widetilde{\beta}) - \w_\ell(\beta) | \leq C | \widetilde{\beta} - \beta |_1 \qquad 
\text{for all } \; \widetilde{\beta}, \beta \in\IR^d \text{ and } \ell=1, \ldots, n.
\end{align*}

\item \label{Ass:L.properties.item.iii}
The eigenvalue $\w=\w_1(\kappa)$ is bounded away from the other eigenvalues:
There is a constant $C$ such that 
\begin{align*}
| \w - \w_\ell(\beta) | \geq C  \qquad 
\text{for all } \; \beta \in\IR^d \text{ and } \ell=2, \ldots, n.
\end{align*}
\end{enumerate}

\end{Assumption}
Assumption \ref{Ass:L.properties.item.i} corresponds to Assumption 2 in \cite{colin-lannes:09}, 
whereas Assumption \ref{Ass:L.properties.item.iii} is a part of Assumption 3 in \cite{colin-lannes:09}.

\begin{Remark}\label{Ass:L.properties.remark}
Explicit formulas for the eigenvalues in case of the Maxwell--Lorentz system and the Klein--Gordon system are given in \cite[Example 3 and 4]{colin-lannes:09}, and one can check that the 
assumptions~\ref{Ass:L.properties.item.i} and \ref{Ass:L.properties.item.ii} on the eigenvalues are true.
Assumption~\ref{Ass:L.properties.item.iii} is true if we choose $\w$ to be the largest or smallest eigenvalue in \eqref{Def.w}.
\end{Remark}

For $j\in\J_+^{(m)}$ and every $\theta\in\IR^d$ let
\begin{subequations}
\label{eigendecomposition}
\begin{align}
 \L_j(\theta) &= \Psi_j(\theta)\Lambda_j(\theta)\Psi_j^*(\theta) 
\end{align}
be the eigendecomposition of \eqref{Def.L.j}:
the real diagonal matrix 
\begin{align}
\Lambda_j(\theta) = \diag(\lambda_{j1}(\theta), \ldots, \lambda_{jn}(\theta)) \in\IR^{n \times n} 
\end{align}
contains the eigenvalues $\lambda_{j\ell}(\theta) \in \IR $
of $ \L_j(\theta) $, and 
\begin{align}
\Psi_j(\theta)=
\begin{pmatrix} \psi_{j1}(\theta) \mid \cdots \mid \psi_{jn}(\theta) \end{pmatrix} \in\IC^{n \times n}
\end{align}
\end{subequations}
is unitary with the corresponding normalized eigenvectors $\psi_{j\ell}(\theta) \in \IC^n $ in its columns.
By Assumption~\ref{Ass:polarization}\ref{Ass:polarization.item.i} $ \L_1(0)=\L(\w,\kappa) $ has a one-dimensional kernel, and we choose the enumeration of the eigenvalues and eigenvectors in such a way that
$ \lambda_{11}(0) = 0 $ and $ \ker \L_1(0)=\text{span}\{\psi_{11}(0)\}.$
Equation~\eqref{Def.L.j.negative} implies that $\Psi_{-j}(\theta)= \overline{\Psi_j(-\theta)}$ and 
$\Lambda_{-j}(\theta)= -\overline{\Lambda_j(-\theta)}= -\Lambda_j(-\theta)$.

The matrices $\L(0,j\kappa+\theta)$ and 
$ \L(j\w,j\kappa+\theta) = -j\w I + \L(0,j\kappa+\theta) $ 
have the same eigenvectors, and their eigenvalues $ \w_\ell(j\kappa+\theta) $
and $ \lambda_{j\ell}(\theta) = -j\w + \w_\ell(j\kappa+\theta) $ differ only by a shift.
Hence, it follows from Assumption~\ref{Ass:L.properties} that
$ \lambda_{j\ell} \in C^\infty(\IR^d\setminus\{-j\kappa\}, \IR) $ and 
$ \psi_{j\ell} \in C^\infty(\IR^d\setminus\{-j\kappa\}, \IC^{n}) $ with
\begin{align}
\label{Ass:lambda.Lipschitz}
| \lambda_{j\ell}(\widetilde{\theta}) - \lambda_{j\ell}(\theta) | 
&\leq C | \widetilde{\theta} - \theta |_1 
&&
\text{for all } \; \widetilde{\theta}, \theta \in\IR^d,
\\
\label{Ass:lambda.1.ell.bounded.below}
| \lambda_{1\ell}(\theta) | &\geq C  
&&
\text{for all } \; \theta \in\IR^d \text{ and } \ell=2, \ldots, n.
\end{align}

\paragraph{Transformation of the coefficient functions.}

The strategy in the proofs of \eqref{error.bound.BJ.1.improved} and \eqref{error.bound.BJ.3.improved} is, roughly speaking, to distinguish the oscillatory ``parts'' of the solution from the non-oscillatory ones, and to carefully analyze how these parts interact in the nonlinearity. 
For this purpose, the following transformation was introduced in \cite{baumstark-jahnke-2023}.

Let $\Uhat^{(m)} = \{\uhat_j\}_{j\in\J_+^{(m)}}$ 
be the solution of \eqref{PDE.mfe.Fourier} for $ \epsi\in(0,1]$.
For every $t\geq 0$ and $k\in\IR^d$ we define
\begin{align}
\label{Def.z}
z_{j}(t,k) &= \TT_{j,\epsi}(t,k)\uhat_j(t,k), \qquad
z_{-j}(t,k) = \overline{z_{j}(t,-k)}, \qquad j\in\J_+^{(m)}
\end{align}
with transformation matrix
\begin{subequations}
\label{Def.TT}
\begin{align}
 \TT_{j,\epsi}(t,k)&= \exp\big(\tfrac{\ii t}{\epsi}\Lambda_{j}(\epsi k)\big)
 \Psi_{j}^*(\epsi k)
 = \Psi_{j}^*(\epsi k) \exp\big(\tfrac{\ii t}{\epsi}\L_{j}(\epsi k)\big), 
 \qquad j\in\J_+^{(m)}, 
 \\
 \TT_{-j,\epsi}(t,k) &:= \overline{\TT_{j,\epsi}(t,-k)}.
\end{align}
\end{subequations}
It follows from \eqref{Def.z} and \eqref{PDE.mfe.Fourier.a} that
\begin{align}
\label{zdot.in.terms.of.uhat}
\pt z_{j}(t) 
&=
\epsi \sum_{\#J=j} \FF\left(t,\Uhat^{(m)},J\right),
\qquad \Uhat^{(m)} = \{\uhat_j\}_{j\in\J_+^{(m)}},
\end{align}
where $\FF$ is given by 
\begin{align}
\label{Def.F.uhat}
\FF(t,\Uhat^{(m)},J) = \TT_{j,\epsi}(t) \T\big(\uhat_{j_1},\uhat_{j_2},\uhat_{j_3}\big)(t), 
\qquad & J=(j_1,j_2,j_3) \in (\J^{(m)})^3, \quad j=\#J.
\end{align}
By means of the inverse transform $ \uhat_j(t,k) = \TT_{j,\epsi}^*(t,k)z_{j}(t,k) $
we could turn \eqref{zdot.in.terms.of.uhat} into a closed system of evolution equations for
$ \{z_j\}_{j\in\J^{(m)}} $, but with a rather complicated right-hand side.
The initial conditions are
\begin{align}
\label{PDE.z.initial.data}
z_{j}(0,k) &= \TT_{j,\epsi}(0,k) \uhat_j(0,k)
  =
  \begin{cases}
  \Psi_1^*(\epsi k) \phat(k) & \text{if } j=1, \\ 0 & \text{if } j\in\J_+^{(m)}\setminus\{1\}
   \end{cases}
\end{align}
according to \eqref{PDE.mfe.Fourier.b}, \eqref{Def.z} and \eqref{Def.TT}.

The transformation \eqref{Def.z} and \eqref{Def.TT} is motivated by the fact that in the \emph{linear} case the exact solution of \eqref{PDE.mfe.Fourier} is 
\begin{align*}
\uhat_j(t,k) = \TT_{j,\epsi}^*(t,k)z_j(0,k) \qquad
 \text{for } \T(\cdot,\cdot,\cdot)=0,
\end{align*}
because $ z_j(t)=z_j(0) $ is constant in time 
for $\T(\cdot,\cdot,\cdot)=0$ according to \eqref{zdot.in.terms.of.uhat} and \eqref{Def.F.uhat}.
But even in the nonlinear case $\T(\cdot,\cdot,\cdot)\not=0$ the right-hand side of \eqref{zdot.in.terms.of.uhat} is formally only $\Ord{\epsi}$ instead of $\Ord{1/\epsi}$ in \eqref{PDE.mfe.Fourier.a}, because the linear part
$ \tfrac{\ii}{\epsi}\L_j(\epsi k)\uhat_j(t,k) $ is cancelled by the transformation.
The transformed functions $z_j$ do still oscillate in time, but the oscillations appear on a much smaller scale, and in this sense, $z_j$ is smoother than $\uhat_j$.

\paragraph{Projectors.}

Recall that by Assumption~\ref{Ass:polarization}\ref{Ass:polarization.item.i},
the matrix $\L(\w,\kappa)=\L_1(0)$ has a one-dimensional kernel spanned by 
$\psi_{11}(0)$. This is the reason why the first eigenspace of the matrix
$ \L_1(\epsi k) = \L_1(0) + \epsi A(k)$ which appears in \eqref{PDE.mfe.Fourier.a}
will play a special role in our analysis.
We denote the orthogonal projection onto this eigenspace by
\begin{align}
 \label{Def.P2}
 \widehat{w} \mapsto \P_\epsi \widehat{w}, \qquad
\P_\epsi(k) = \psi_{11}(\epsi k)\psi_{11}^*(\epsi k) \in \IC^{n \times n}
\end{align}
and the projector onto the orthogonal complement by $\P_\epsi^\perp = I-\P_\epsi.$
Assumption \ref{Ass:polarization}\ref{Ass:polarization.item.ii} is equivalent to
$ \P_0^\perp\phat=\epsi \P_0^\perp \widehat{p}_1$, and for $ p_0, p_1 \in W^1 $ 
it was shown in the proof of Lemma 3 in \cite{colin-lannes:09} that
\begin{align}
\label{Bound.Pperp.p}
\| \P_\epsi^\perp \widehat{p} \|_{L^1} &\leq C \epsi (\| p_1 \|_W + \| \nabla p \|_{W})
\leq C \epsi (\| p_0 \|_{W^1} + \| p_1 \|_{W^1}).
\end{align}
For the transformed function \eqref{Def.z} we obtain from \eqref{Def.TT} that
\begin{align}
\label{Relation.between.projections}
\P_\epsi(k)\uhat_1(t,k) 
= 
\psi_{11}(\epsi k) \exp\big(-\tfrac{\ii t}{\epsi}\lambda_{11}(\epsi k)\big) z_{11}(t,k)
=
\TT_{1,\epsi}^*(t,k) Pz_1(t,k),
\end{align}
where $ z_{11}(t,k) $ is the first entry of $ z_{1}(t,k)\in\IC^n $ and 
\begin{align}
\label{Def.P1}
P: \IC^n \rightarrow \IC^n, \quad
(w_1, \ldots, w_n)^\top \mapsto (w_1, 0, \ldots, 0)^\top
\end{align}
is the orthogonal projection of a vector $w$ onto $\text{span}\{(1,0,\ldots,0)^\top\}$.
For $\Pperp=(I-P)$ the estimate \eqref{Bound.Pperp.p} yields 
\begin{align}
 \label{Bound.projection.initial.value}
 \| \Pperp z_1(0,\cdot) \|_{L^1} \leq C \epsi (\| p_0 \|_{W^1} + \| p_1 \|_{W^1}),
 \end{align}
because with \eqref{Relation.between.projections} we obtain
\begin{align*}
 \Pperp z_1(0,\cdot) = z_1(0,\cdot) - Pz_1(0,\cdot)
 = \TT_{1,\epsi}(0,k)\P_\epsi^\perp(k)\uhat_1(0,k).
\end{align*}

\paragraph{Useful identities and inequalities.}

Throughout, we will frequently use the following facts.
Since we have chosen the Euclidean vector norm $ | \cdot |_2 $ to define $ \| \cdot \|_{L^1} $
in \eqref{Def.W}, the norm $ \| \widehat{f} \|_{L^1} $ of $\widehat{f} \in L^1$ is invariant under multiplication of 
$ \widehat{f}(k) \in \IC^n $ with a unitary matrix $ S(k) \in \IC^{n \times n}$.
This means, in particular, that for the transformed functions 
$ z_{j}(t,k) = \TT_{j,\epsi}(t,k)\uhat_j(t,k) $ from \eqref{Def.z} the identities
\begin{align}
\label{Norm.uhat.z}
| z_j(t,k) |_2 &= | \uhat_j(t,k) |_2,
&
\| z_{j}(t) \|_{L^1} &= \| \uhat_{j}(t) \|_{L^1} = \| u_j(t) \|_W
\end{align}
and, via \eqref{Relation.between.projections}, the equations
\begin{align}
\notag
% \label{Norm.uhat.z.P}
| Pz_1(t,k) |_2 &= | \P_\epsi(k) \uhat_1(t,k) |_2, &
\| Pz_1(t) \|_{L^1} &= \| \P_\epsi \uhat_1(t) \|_{L^1}, 
\\
\label{Norm.uhat.z.Pperp}
| \Pperp z_1(t,k) |_2 &= | \P_\epsi^\perp(k) \uhat_1(t,k) |_2, &
\| \Pperp z_1(t) \|_{L^1} &= \| \P_\epsi^\perp \uhat_1(t) \|_{L^1}
\intertext{hold for all $t\geq 0$, $k\in\IR^d$, and $\epsi\in(0,1]$.
Moreover, we will use that for all $ w\in\IC^n $ and $f \in L^1$
the inequalities
}
\label{P.decrease.norm}
| Pw |_2 &\leq |w|_2, &
\| Pf \|_{L^1} &\leq \| f \|_{L^1}, 
\\
\label{Pepsi.decrease.norm}
| \P_\epsi(k) w |_2 &\leq | w |_2, &
\| \P_\epsi f \|_{L^1} &\leq \| f \|_{L^1}
\end{align}
hold, as well as the same inequalities with $P$ and $\P_\epsi$ replaced by $\Pperp$ and $\P_\epsi^\perp$, respectively.

% ---------------------------------------------------------------------------------------
\section{Why $ \P_\epsi \uhat_1 $ is smooth and $ \P_\epsi^\perp \uhat_1 $ is small
in the slowly varying envelope approximation}\label{Sec.error.SVEA.preparations}
% ---------------------------------------------------------------------------------------

In this and the next section we analyze the SVEA \eqref{Ansatz.SVEA}--\eqref{SVEA}, 
which corresponds to setting
\begin{align*}
j=m=1, \qquad \J^{(m)}=\J^{(1)}=\{-1,1\}, \qquad \J_+^{(1)}=\{1\} 
\end{align*}
in \eqref{Ansatz}--\eqref{PDE.mfe} and in \eqref{PDE.mfe.Fourier}, respectively.
Our main goal is to prove the error bound \eqref{error.bound.BJ.1.improved},
which will be achieved in Section~\ref{Sec.error.bound.SVEA}; 
cf.~Theorem~\ref{Theorem.error.bound.SVEA} below.
This proof is based on a number of auxiliary results, 
which we compile now.
We start by quoting two important inequalities from \cite{colin-lannes:09}.

\begin{Lemma}\label{Lemma.bound.uhat1dot}
Let $m = 1$, let $\sigma=1$ in \eqref{Def:p}, and let $u_1$ be the classical solution of \eqref{SVEA} 
which was established in Lemma~\ref{Lemma.wellposedness}\ref{Lemma.wellposedness.1}.
Under Assumptions \ref{Ass:L.properties} and 
\ref{Ass:polarization}\ref{Ass:polarization.item.i}, there is a constant $C$ such that
 \begin{align} \notag
\sup_{t\in[0,\tend/\epsi]} 
\| \pt \P_\epsi \uhat_1(t) \|_{L^1} \leq C.
\end{align}
The constant $C$ depends on $C_{u,1}$ from \eqref{Lemma.wellposedness.bound.23} and thus also on $\tend$, but not on $\epsi\in(0,1]$.
\end{Lemma}
\prooftext{:} See \cite[Lemma 2]{colin-lannes:09}.

\begin{Proposition}\label{Proposition01.SVEA}
Let $m=1$ and let $u_1$ be the classical solution of \eqref{SVEA} with initial data of the form \eqref{Def:p} with $\sigma=1$.
Under the Assumptions~\ref{Ass:polarization} and \ref{Ass:L.properties}, there is a constant $C$ such that
\begin{align}
\label{Bound.uhat.1}
\sup_{t\in[0,\tend/\epsi]}\| \P_\epsi^\perp \uhat_1(t) \|_{L^1} 
&\leq C \epsi
\end{align}
for all $\epsi\in(0,1]$. 
\end{Proposition}
\prooftext{:} See~\cite[Lemma 3]{colin-lannes:09}. 
In \cite{baumstark:22} a similar result was shown without 
Assumption~\ref{Ass:L.properties}\ref{Ass:L.properties.item.iii}, but on a possibly smaller interval
$ [0,\tstar/\epsi] $ for some $ \tstar \leq \tend $.

\bigskip

\noindent
These results can be interpreted as follows.
The term $ \tfrac{\ii}{\epsi}\L_1(\epsi k)\uhat_1(t,k) $ in \eqref{PDE.mfe.Fourier.a} suggests that
formally $ \pt \uhat_1 = \Ord{1/\epsi} $.
Lemma~\ref{Lemma.bound.uhat1dot} shows, however, that the time derivative of the \emph{projected} part 
$ \P_\epsi \uhat_1 $ is bounded uniformly in $\epsi$. 
Hence, we can consider $ \P_\epsi \uhat_1 $ as ``the non-oscillatory part of $\uhat_1$'', although strictly speaking this interpretation is not correct, because oscillations in $ \P_\epsi \uhat_1 $ can still be detected on a very small scale;
cf.~Remark~\ref{Remark:interpretation.Puhat1.smooth} at the end of this subsection.

For the time derivatives of the other part $ \P_\epsi^\perp \uhat_1 = \uhat_1 - \P_\epsi \uhat_1 $ a corresponding result does not hold, which means that $ \pt \P_\epsi^\perp \uhat_1 = \Ord{1/\epsi} $ in general.
Proposition~\ref{Proposition01.SVEA} shows, however, 
that $ \| \P_\epsi^\perp \uhat_1(t) \|_{L^1} = \Ord{\epsi} $ 
even on the long time interval $ [0,\tend/\epsi] $.
Hence, we can think of $ \P_\epsi^\perp \uhat_1 $ as ``small but oscillatory'' in the sense that its time derivative is much larger than $ \P_\epsi^\perp \uhat_1 $ itself.
Exploiting the different properties of $ \P_\epsi \uhat_1 $ and $ \P_\epsi^\perp \uhat_1 $ will be crucial 
in the proof of Theorem~\ref{Theorem.error.bound.SVEA} in Section~\ref{Sec.error.bound.SVEA}.
Before that, we have to extend Lemma~\ref{Lemma.bound.uhat1dot} and Proposition~\ref{Proposition01.SVEA}
to a stronger norm.

Let $ D_\mu $ denote the Fourier multiplier 
$(D_\mu \widehat{w})(k) = \ii k_\mu \widehat{w}(k)$ for $\mu\in\{1, \ldots, d\}$.
We want to show that under stronger regularity assumptions 
Proposition~\ref{Proposition01.SVEA} remains true when $ \P_\epsi^\perp \uhat_1(t) $ 
is replaced by $ D_\mu \P_\epsi^\perp \uhat_1(t) $; 
cf.~Proposition~\ref{Proposition02.SVEA} below.
This corresponds to an extension of the inequality \eqref{Bound.uhat.1} from 
\begin{align*}
\| \P_\epsi^\perp \uhat_1(t) \|_{L^1} &= \| \F^{-1}(\P_\epsi^\perp \uhat_1(t)) \|_W 
\intertext{to the stronger norm }
\| \P_\epsi^\perp \uhat_1(t) \|_{L^1} + \sum_{\mu=1}^d \| D_\mu \P_\epsi^\perp \uhat_1(t) \|_{L^1}
&= \| \F^{-1}(\P_\epsi^\perp \uhat_1(t)) \|_{W^1}.
\end{align*}
As a first step, we prove the following counterpart of Lemma~\ref{Lemma.bound.uhat1dot}.

\begin{Lemma}\label{Lemma.bound.Dmu.uhat1dot}
Let $m = 1$, let $\sigma=2$ in \eqref{Def:p}, and let $u_1$ be the classical solution of \eqref{SVEA} 
which was established in Lemma~\ref{Lemma.wellposedness}\ref{Lemma.wellposedness.23}.
Under Assumptions \ref{Ass:L.properties} and 
\ref{Ass:polarization}\ref{Ass:polarization.item.i}, there is a constant $C$ such that
\begin{align*}
\sup_{t\in[0,\tend/\epsi]} 
\| \pt D_\mu \P_\epsi \uhat_1(t) \|_{L^1} \leq C.
\end{align*}
The constant $C$ depends on $C_{u,2}$ from \eqref{Lemma.wellposedness.bound.23} and thus also on $\tend$, but not on $\epsi\in(0,1]$.
\end{Lemma}

\proof
The proof is similar to the proof of Lemma~\ref{Lemma.bound.uhat1dot}. 
We choose $\mu\in\{1, \ldots, d\}$ and apply $D_\mu \P_\epsi(k) $ 
to both sides of \eqref{PDE.mfe.Fourier.a} with $j=m=1$.
This yields
\begin{align}
 \label{Lemma.bound.Dmu.uhat1dot.proof.1}
\pt D_\mu \P_\epsi(k) \uhat_1(t,k) 
= - \frac{\ii}{\epsi}D_\mu\P_\epsi(k)\L_1(\epsi k)\uhat_1(t,k)
+
\epsi D_\mu \P_\epsi(k) \sum_{\#J=1} \T\big(\uhat_{j_1},\uhat_{j_2},\uhat_{j_3}\big)(t,k) 
\end{align}
for all $ t\in(0,\tend/\epsi] $ and $ k \in \IR^d $. 
The first term on the right-hand side is
\begin{align}
 -\frac{\ii}{\epsi}D_\mu \P_\epsi(k) \L_1(\epsi k)\uhat_1(t,k)
 &=
  \label{Lemma.bound.Dmu.uhat1dot.proof.2}
 -\frac{\ii}{\epsi}\lambda_{11}(\epsi k)\P_\epsi(k) D_\mu \uhat_1(t,k)
 \end{align}
because of \eqref{Def.P2} and \eqref{eigendecomposition}. 
The Lipschitz continuity \eqref{Ass:lambda.Lipschitz} of the eigenvalues 
and the fact that $\lambda_{11}(0)=0$ yield
\begin{align*}
| \lambda_{11}(\epsi k) | = | \lambda_{11}(\epsi k) - \lambda_{11}(0) | \leq C\epsi |k|_1,
\end{align*}
and together with \eqref{Lemma.bound.Dmu.uhat1dot.proof.2} and \eqref{Pepsi.decrease.norm}, this gives
\begin{align}
\notag
\Big\| \frac{\ii}{\epsi}D_\mu \P_\epsi \L_1(\epsi \, \cdot)\uhat_1(t) \Big\|_{L^1}
& \leq 
C \intl_{\IR^d}
|k|_1 |\P_\epsi(k) D_\mu \uhat_1(t,k)|_2 \; \dd k
\\
\label{Lemma.bound.Dmu.uhat1dot.proof.3}
& \leq
C \intl_{\IR^d}
|k|_1^2 |\uhat_1(t,k)|_2 \; \dd k
=
C \cdot C_{u,2}
\end{align}
with $ C_{u,2} $ from
Lemma~\ref{Lemma.wellposedness}\ref{Lemma.wellposedness.23}.
For the nonlinear term on the right-hand side of \eqref{Lemma.bound.Dmu.uhat1dot.proof.1}, we have 
\begin{align}
\label{Product.rule.DT}
D_\mu \T\big(\uhat_{j_1},\uhat_{j_2},\uhat_{j_3}\big)
&=
\T\big(D_\mu\uhat_{j_1},\uhat_{j_2},\uhat_{j_3}\big)
+ \T\big(\uhat_{j_1},D_\mu\uhat_{j_2},\uhat_{j_3}\big) + \T\big(\uhat_{j_1},\uhat_{j_2},D_\mu\uhat_{j_3}\big),
\end{align}
which corresponds to the product rule.
Since there are three multi-indices $J\in(\J^{(1)})^3$ with $\#J=1$,
namely $(1,1,-1),(1,-1,1),(-1,1,1)$, we obtain with \eqref{Bound.T.Fourier}
\begin{align}
\label{Lemma.bound.Dmu.uhat1dot.proof.4}
\epsi \big\| D_\mu \P_\epsi \sum_{\#J=1} \T\big(\uhat_{j_1},\uhat_{j_2},\uhat_{j_3}\big)(t)  \big\|_{L^1}
&\leq
3 \epsi \big\| D_\mu \T\big(\uhat_{j_1},\uhat_{j_2},\uhat_{j_3}\big)(t)  \big\|_{L^1}
\leq
9 \epsi C_\T C_{u,1}^3.
\end{align}
The assertion follows by combining 
\eqref{Lemma.bound.Dmu.uhat1dot.proof.1}, \eqref{Lemma.bound.Dmu.uhat1dot.proof.3}, 
\eqref{Lemma.bound.Dmu.uhat1dot.proof.4}, and using that $\epsi\leq 1$ by assumption.
\qed

\noindent
With Lemma~\ref{Lemma.bound.Dmu.uhat1dot}, we can now show the following extension of Proposition~\ref{Proposition01.SVEA}.

\begin{Proposition}\label{Proposition02.SVEA} 
Let $m=1$ and let $u_1$ be the classical solution of \eqref{SVEA} with initial data of the form \eqref{Def:p} with $\sigma=2$.
Under the assumptions of Proposition~\ref{Proposition01.SVEA} there is a constant $C$ such that
\begin{align}
\label{Bound.Dmu.uhat.1}
\sup_{t\in[0,\tend/\epsi]}\| D_\mu \P_\epsi^\perp \uhat_1(t) \|_{L^1} 
&\leq C \epsi
\end{align}
for all $\epsi\in(0,1]$ and all $\mu\in\{1, \ldots, d\}$.
\end{Proposition}

\proof
Choose a fixed $\mu\in\{1, \ldots, d\}$ and set
\begin{align}
\label{Proposition02.SVEA.Def.vhat} 
\vhat(t,k) = (D_\mu \P_\epsi^\perp \uhat_1)(t,k) = \ii k_\mu \P_\epsi^\perp(k) \uhat_1(t,k).
\end{align}
We apply $ D_\mu \P_\epsi^\perp $ to \eqref{PDE.mfe.Fourier} with $m=j=1$
and use that $ \P_\epsi^\perp $ commutes with $ \L_1(\epsi k) $.
This yields
\begin{align*}
\pt \vhat(t) + \frac{\ii}{\epsi}\L_1(\epsi \, \cdot)\vhat(t)
&=
\epsi \sum_{\#J=1} D_\mu \P_\epsi^\perp \T\big(\uhat_{j_1},\uhat_{j_2},\uhat_{j_3}\big)(t),
\\
\vhat(0) &= D_\mu \P_\epsi^\perp \phat
\end{align*}
with $\L_1(\epsi \, \cdot)$ denoting $ k \mapsto \L_1(\epsi k)$.
Now we adapt the proof of Proposition~\ref{Proposition01.SVEA}.
With Duhamel's formula and the short-hand notation
\begin{align*}
\T(\uhat_J) = \T(\uhat_{j_1},\uhat_{j_2},\uhat_{j_3}) \qquad \text{for } J=(j_1,j_2,j_3),
\end{align*}
we obtain
\begin{align*}
% \label{Proposition02.SVEA.proof.1} 
\vhat(t) &= \vhat^{[1]}(t) + \vhat^{[2]}(t) + \vhat^{[3]}(t) 
\end{align*}
with the three terms
\begin{align*}
\vhat^{[1]}(t)
&= 
\exp\left(-\frac{\ii t}{\epsi}\L_1(\epsi \, \cdot)\right) D_\mu \P_\epsi^\perp \phat,
\\
\vhat^{[2]}(t)
&=
\epsi
\sum_{\#J=1} \intl_0^t 
\exp\left(\frac{\ii (s-t)}{\epsi}\L_1(\epsi \, \cdot)\right)
D_\mu \P_\epsi^\perp \T\big(\P_\epsi \uhat_J(s)\big) \; \dd s,
\\
\vhat^{[3]}(t)
&=
\epsi
\sum_{\#J=1} \intl_0^t 
\exp\left(\frac{\ii (s-t)}{\epsi}\L_1(\epsi \, \cdot)\right)
D_\mu \P_\epsi^\perp 
\Big[\T\big(\uhat_J(s)\big) - \T\big(\P_\epsi \uhat_J(s)\big)\Big]
\; \dd s.
\end{align*}
We will show that 
\begin{align}
\label{Proposition02.SVEA.proof.goal} 
\| \vhat^{[\eta]}(t) \|_{L^1} 
\leq 
c_1\epsi + c_2 \epsi \intl_0^t \big\| \vhat(s) \big\|_{L^1} \; \dd s
\qquad \text{for } t\in[0,\tend/\epsi] \text{ and } \eta=1,2,3 
\end{align}
with constants $c_1\geq 0$ and $c_2\geq 0$ which do not depend on $ \epsi\in(0,1] $.
If \eqref{Proposition02.SVEA.proof.goal} is true, 
then applying Gronwall's lemma and using that $ \epsi t \leq \tend $ proves that
$ \sup_{t\in[0,\tend/\epsi]}\| \vhat(t) \|_{L^1} \leq C \epsi $ which,
via \eqref{Proposition02.SVEA.Def.vhat}, is equivalent to \eqref{Bound.Dmu.uhat.1}.

For the first term $\vhat^{[1]}(t)$, the inequality \eqref{Bound.Pperp.p} implies
\begin{align*}
% \label{Proposition02.SVEA.proof.2} 
\big\| \vhat^{[1]}(t) \big\|_{L^1}
&= 
\| \P_\epsi^\perp D_\mu \phat \|_{L^1}
\leq C \epsi (\| p_0 \|_{W^2} + \| p_1 \|_{W^2}),
\end{align*}
which verifies \eqref{Proposition02.SVEA.proof.goal} for $\eta=1$ (with $c_2=0$).

For the third term $\vhat^{[3]}(t)$ we infer with 
\eqref{Pepsi.decrease.norm}, \eqref{Product.rule.DT}, \eqref{Bound.T.difference.Fourier},
and Proposition~\ref{Proposition01.SVEA} that
\begin{align*}
\notag
\|  \vhat^{[3]}(t) \|_{L^1}
&\leq
3 \epsi
\intl_0^t 
\big\| D_\mu \T\big(\uhat_J(s)\big) - D_\mu \T\big(\P_\epsi \uhat_J(s)\big) \big\|_{L^1}
\; \dd s
\\
\notag
&\leq
C \epsi
\intl_0^t 
\Big( \big\| \uhat_1(s) - \P_\epsi \uhat_1(s) \big\|_{L^1} 
+ \big\| D_\mu \uhat_1(s) - D_\mu \P_\epsi \uhat_1(s) \big\|_{L^1} \Big)
\; \dd s
\\
\notag
&=
C \epsi
\intl_0^t 
\Big(\big\| \P_\epsi^\perp \uhat_1(s) \big\|_{L^1} + \big\| \vhat(s) \big\|_{L^1}\Big)
\; \dd s
\\
% \label{Proposition02.SVEA.proof.3} 
&\leq
c_1 \epsi + 
c_2 \epsi
\intl_0^t \big\| \vhat(s) \big\|_{L^1} \; \dd s
\end{align*}
with constants $c_1, c_2$ which depend on $ C_\T $, $ C_{u,1} $, 
and in case of $c_1$ also on the constant from \eqref{Bound.uhat.1}.

Now we consider the second term $\vhat^{[2]}(t)$.
Since there are three multi-indices $J\in(\J^{(1)})^3$ with $\#J=1$, we obtain
\begin{align}
\notag
\|  \vhat^{[2]}(t) \|_{L^1}
&\leq
3 \epsi
\intl_{\IR^d}
\Big| 
\intl_0^t 
\exp\left(\frac{\ii (s-t)}{\epsi}\L_1(\epsi k)\right)
\big(D_\mu \P_\epsi^\perp \T(\P_\epsi \uhat_J)\big)(s,k) \; \dd s
\Big|_2
\; \dd k
\\
\label{Proposition02.SVEA.proof.4} 
&=
3 \epsi
\intl_{\IR^d}
\Big| 
\intl_0^t 
\exp\left(\frac{\ii (s-t)}{\epsi}\L_1(\epsi k)\right)
\P_\epsi^\perp(k)
\qhat(s,k) \; \dd s
\Big|_2
\; \dd k
\end{align}
with the abbreviation
\begin{align}
\label{Proposition02.SVEA.proof.Def.qhat}
\qhat(s,k)=\big(D_\mu \T(\P_\epsi \uhat_J)\big)(s,k).
\end{align}
The goal is now to integrate by parts to gain one additional factor $\epsi$, which is then used to compensate the long time interval.
However, this requires some care, because the matrix $ \L_1(0) = \L(\w,\kappa) $ is singular;
see \eqref{Def.w} or Assumption~\ref{Ass:polarization}\ref{Ass:polarization.item.i}.
What saves us here is the projector $ \P_\epsi^\perp(k) $ in \eqref{Proposition02.SVEA.proof.4}.
For every $k\in\IR^d$, the restriction of $ \L_1(\epsi k) $ to the subspace $ \P_\epsi^\perp(k)\IR^n $
is given by
\begin{align*}
\L_1^\perp(\epsi k): \P_\epsi^\perp(k)\IR^n \to \P_\epsi^\perp(k)\IR^n, \qquad
\L_1^\perp(\epsi k) = \L_1(\epsi k)\P_\epsi^\perp(k) =
\sum_{\ell=2}^n \lambda_{1\ell}(\epsi k) \psi_{1\ell}(\epsi k)\psi_{1\ell}^*(\epsi k).
\end{align*}
By \eqref{Ass:lambda.1.ell.bounded.below}, this mapping is regular with uniformly bounded inverse
\begin{align*}
\big(\L_1^\perp(\epsi k)\big)^{-1}: \P_\epsi^\perp(k)\IR^n \to \P_\epsi^\perp(k)\IR^n, \qquad
\big(\L_1^\perp(\epsi k)\big)^{-1} =
\sum_{\ell=2}^n \frac{1}{\lambda_{1\ell}(\epsi k)} \psi_{1\ell}(\epsi k)\psi_{1\ell}^*(\epsi k).
\end{align*}
The presence of $ \P_\epsi^\perp(k) $ in \eqref{Proposition02.SVEA.proof.4}
allows us to replace $ \L_1(\epsi k) $ by $ \L_1^\perp(\epsi k) $ and to
integrate by parts in the inner integral of \eqref{Proposition02.SVEA.proof.4}.
This yields
\begin{align*}
& \Big| 
\intl_0^t 
\exp\left(\frac{\ii (s-t)}{\epsi}\L_1(\epsi k)\right)
\P_\epsi^\perp(k) \qhat(s,k)
\; \dd s
\Big|_2
\\
&\leq
\Big|
\frac{\epsi}{\ii}
\big(\L_1^\perp(\epsi k)\big)^{-1}
\qhat(t,k)
-
\frac{\epsi}{\ii}
\exp\left(-\frac{\ii t}{\epsi}\L_1^\perp(\epsi k)\right) 
\big(\L_1^\perp(\epsi k)\big)^{-1}
\qhat(0,k)
\Big|_2
\\
&\quad
+
\Big| 
\frac{\epsi}{\ii}
\intl_0^t 
\exp\left(\frac{\ii (s-t)}{\epsi}\L_1^\perp(\epsi k)\right) 
\big(\L_1^\perp(\epsi k)\big)^{-1}
\pt\qhat(s,k)
\; \dd s
\Big|_2
\\
&\leq
C\epsi \Big(|\qhat(t,k)|_2 + |\qhat(0,k)|_2\Big)
+ C \epsi 
\intl_0^t |\pt\qhat(s,k)|_2 \; \dd s,
\end{align*}
and substituting this into \eqref{Proposition02.SVEA.proof.4} leads to 
\begin{align}
\label{Proposition02.SVEA.proof.5}
\| \vhat^{[2]}(t) \|_{L^1}
&\leq 
C\epsi^2 \Big(\|\qhat(t)\|_{L^1} + \|\qhat(0)\|_{L^1}\Big)
+ C \epsi \tend
\sup_{s\in[0,\tend/\epsi]} \|\pt\qhat(s)\|_{L^1}.
\end{align}
With \eqref{Proposition02.SVEA.proof.Def.qhat}, \eqref{Product.rule.DT}, and \eqref{Bound.T.Fourier}
we obtain that
\begin{align}
\label{Proposition02.SVEA.proof.6}
\|\qhat(t)\|_{L^1} &= \|\big(D_\mu \T(\P_\epsi \uhat_J)\big)(t)\|_{L^1}
\leq 3C_\T C_{u,1}^3,
&
\|\qhat(0)\|_{L^1} &\leq 3C_\T C_{u,1}^3,
\end{align}
and that
\begin{align*}
\|\pt\qhat(s)\|_{L^1} 
&=
\|\pt D_\mu \T\big(\P_\epsi \uhat_J(s)\big)\|_{L^1}
\\
&\leq
6C_\T \|\pt \P_\epsi \uhat_1(s)\|_{L^1} \|D_\mu \P_\epsi \uhat_1(s)\|_{L^1} \|\P_\epsi \uhat_1(s)\|_{L^1} 
+
3C_\T \|\pt D_\mu \P_\epsi \uhat_1(s)\|_{L^1} \|\P_\epsi \uhat_1(s)\|_{L^1}^2.
\end{align*}
Since $ \|\pt \P_\epsi \uhat_1(s)\|_{L^1} $ and $ \|\pt D_\mu \P_\epsi \uhat_1(s)\|_{L^1} $ are uniformly
bounded by Lemmas~\ref{Lemma.bound.uhat1dot} and \ref{Lemma.bound.Dmu.uhat1dot}, respectively,
this shows that $\|\pt\qhat(s)\|_{L^1}$
is uniformly bounded in $ s\in[0,\tend/\epsi]$ and $\epsi\in(0,1]$.
Combining this with \eqref{Proposition02.SVEA.proof.6} and \eqref{Proposition02.SVEA.proof.5} yields
\eqref{Proposition02.SVEA.proof.goal} for $\eta=2$.
This completes the proof.
\qed

Before closing this section we prove that even the second time derivative of $\P_\epsi \uhat_1(t)$ is uniformly bounded. 
This somewhat simple observation will be crucial for showing the error bound for the SVEA;
cf.~\eqref{Step.6.bound.Tdotdot} in step~6 of the proof of Theorem~\ref{Theorem.error.bound.SVEA} below.

\begin{Lemma}\label{Lemma.bound.uhat1dotdot}
Let $m = 1$, let $\sigma=2$ in \eqref{Def:p}, and let $u_1$ be the classical solution of \eqref{SVEA}.
Under Assumptions \ref{Ass:L.properties} and \ref{Ass:polarization}\ref{Ass:polarization.item.i}, 
the second time derivative of $\P_\epsi \uhat_1(t)$ is uniformly bounded, 
i.e.\ there is a constant $C$ such that
 \begin{align} \notag
\sup_{t\in[0,\tend/\epsi]} 
\| \pt^2 \P_\epsi \uhat_1(t) \|_{L^1} \leq C.
\end{align}
The constant $C$ depends on the constant $C_{u,2}$ from \eqref{Lemma.wellposedness.bound.23} and thus also on $\tend$, but not on $\epsi$.
\end{Lemma}

\proof
Applying $ \P_\epsi(k) \pt $ on both sides of \eqref{PDE.mfe.Fourier.a} with $j=m=1$ gives
\begin{align*}
% \label{Lemma.bound.uhat1dot.proof.1}
\P_\epsi(k)\pt^2 \uhat_1(t,k) 
= - \frac{\ii}{\epsi}\P_\epsi(k)\L_1(\epsi k)\pt\uhat_1(t,k)
+
\epsi \P_\epsi(k) \sum_{\#J=1} \pt \T\big(\uhat_{j_1},\uhat_{j_2},\uhat_{j_3}\big)(t,k) 
\end{align*}
for $ t\in(0,\tend/\epsi] $ and $ k \in \IR^d $. 
By adapting the arguments from the proof of Lemma~\ref{Lemma.bound.Dmu.uhat1dot} we arrive at the bound
\begin{align} \notag
\Big| \frac{\ii}{\epsi} \P_\epsi \L_1(\epsi k)\pt \uhat_1(t,k) \Big|_2
&\leq 
C |k|_1 \, \big|\P_\epsi \pt \uhat_1(t,k) \big|_2
=
C \sum_{\mu=1}^d \big| \pt D_\mu \P_\epsi \uhat_1(t,k) \big|_2
\end{align}
for the first term. For the nonlinear term, the product rule yields
\begin{align*}
% \label{Product.rule.ptT}
\pt \T\big(\uhat_{j_1},\uhat_{j_2},\uhat_{j_3}\big)
&=
\T\big(\pt \uhat_{j_1},\uhat_{j_2},\uhat_{j_3}\big)
+ \T\big(\uhat_{j_1},\pt \uhat_{j_2},\uhat_{j_3}\big) + \T\big(\uhat_{j_1},\uhat_{j_2},\pt \uhat_{j_3}\big),
\end{align*} 
and with \eqref{Bound.T.Fourier} we obtain
\begin{align} \notag
\| \pt^2 \P_\epsi \uhat_1(t) \|_{L^1} 
&\leq  
C \sum_{\mu=1}^d \big\| \pt D_\mu \P_\epsi \uhat_1(t) \big\|_{L^1}
+ 9 \epsi C_\T \| \pt \uhat_1 (t) \|_{L^1} \| \uhat_1 (t) \|_{L^1}^2. 
\end{align}
Now the assertion follows from Lemma~\ref{Lemma.bound.Dmu.uhat1dot} and the fact that $\| \pt \uhat_1 (t) \|_{L^1} \leq C \epsi^{-1}$.
\qed

\begin{Remark}\label{Remark:interpretation.Puhat1.smooth}
By taking more derivatives of \eqref{PDE.mfe.Fourier.a} and proceeding as in the proof of 
Lemma~\ref{Lemma.bound.uhat1dotdot}, it can be shown that 
$ \| \pt^\ell \P_\epsi \uhat_1(t) \|_{L^1} = \Ord{\epsi^{2-\ell}} $ for 
$ \ell\geq 3$.
Hence, higher-order time derivatives are \emph{not} uniformly bounded, which means that 
our interpretation of $ \P_\epsi \uhat_1 $ as the non-oscillatory part of $\uhat_1$ 
is only true to a certain extent.
\end{Remark}

% --------------------------------------------------------------------------------
\section{Convergence analysis for the slowly varying envelope approximation}
\label{Sec.error.bound.SVEA}
% --------------------------------------------------------------------------------

With the results from the previous section we are now in a position to 
prove the error bound \eqref{error.bound.BJ.1.improved},
where $\uutilde^{(1)}$ is the SVEA \eqref{Ansatz.SVEA}--\eqref{SVEA}.
We assume that $p$ has the form \eqref{Def:p} with $\sigma=2$.
Then, by Lemmas~\ref{Lemma.wellposedness.uu} and \ref{Lemma.wellposedness}, 
there is a constant $ C_\uu$ such that
\begin{align} 
\label{Def.Cboldu}
\sup_{\epsi\in(0,1]} \sup_{t\in[0,\tend/\epsi]} \| \uu(t) \|_W \leq C_\uu
\qquad \text{and} \qquad 
 \sup_{\epsi\in(0,1]} \sup_{t\in[0,\tend/\epsi]} \| \uutilde^{(1)}(t) \|_W  \leq C_\uu.
\end{align} 
\noindent
The error bound requires the following assumption on the eigenvalues of $\L_j(0) = \L(j\w,j\kappa)$.

\begin{Assumption}[Non-resonance condition]
\label{Ass:Delta.lambda.j3}
The matrix $ \L(3\w, 3\kappa)$ is regular and has no common eigenvalues with
$\L_1(0)=\L(\w,\kappa)$, 
i.e. $ \lambda_{3i}(0) \neq \lambda_{1\ell}(0) $ for all $i, \ell \in \{1, \ldots, n\}$.
\end{Assumption}

\begin{Remark} 
As mentioned earlier, explicit formulas for the eigenvalues 
in case of the Klein--Gordon system and the Maxwell--Lorentz system
can be found in \cite[Example 3 and 4]{colin-lannes:09}.
For these applications, one can check that Assumption~\ref{Ass:Delta.lambda.j3} holds 
if the chosen eigenvalue $\w=\w(\kappa)$ is not constant with respect to $\kappa$.
\end{Remark}

% --------------------------------------------------------------------------------
\subsection{Improved error bound for the SVEA}
% --------------------------------------------------------------------------------

The following theorem is our first main result. It states that the SVEA converges with 
\emph{second} order.
We recall that the SVEA \eqref{Ansatz.SVEA}--\eqref{SVEA} is identical to \eqref{Ansatz}--\eqref{PDE.mfe} 
with $m = 1$ and $\J^{(1)}=\{-1,1\} $.

\begin{Theorem}[Error bound for the SVEA]\label{Theorem.error.bound.SVEA}
Suppose that \eqref{Def:p} holds with $\sigma=2$, and let $\uu$ be the solution of \eqref{PDE.uu}.
Let $u_1$ be the classical solution of \eqref{SVEA}
established in part~\ref{Lemma.wellposedness.23} of Lemma~\ref{Lemma.wellposedness}, 
and let $\uutilde^{(1)}$ be the approximation defined in \eqref{Ansatz.SVEA}.
Under Assumptions~\ref{Ass:polarization}, \ref{Ass:L.properties}, and \ref{Ass:Delta.lambda.j3}
there is a constant such that 
\begin{align}
\label{Theorem.error.bound.SVEA.assertion01}
 \sup_{t\in[0,\tend/\epsi]} \| \uu(t) - \uutilde^{(1)}(t) \|_W
 &\leq C\epsi^2,
 \\
\label{Theorem.error.bound.SVEA.assertion02}
 \sup_{t\in[0,\tend/\epsi]} \| \uu(t) - \uutilde^{(1)}(t) \|_{L^\infty}
 &\leq C\epsi^2.
\end{align}
\end{Theorem}

\proof
The error bound \eqref{Theorem.error.bound.SVEA.assertion02} follows directly from 
\eqref{Theorem.error.bound.SVEA.assertion01} via the embedding $W \hookrightarrow L^\infty$.
The proof of \eqref{Theorem.error.bound.SVEA.assertion01}, however, is rather long.
The strategy, notation and presentation is very similar to the proof of Theorem 4.2 
in \cite{baumstark-jahnke-2023}, but there are some crucial differences which we point out below.

\paragraph{Step 1.} 
In the first step, we derive an evolution equation for the difference
$ \delta=\uu-\uutilde^{(1)} $ between the exact solution and its approximation.
Let 
\begin{align*}
% \label{Def.R}
R =
\epsi T(\uutilde^{(1)},\uutilde^{(1)},\uutilde^{(1)}) - 
\Big(\pt \uutilde^{(1)}(t,x) +A(\partial)\uutilde^{(1)} +\frac{1}{\epsi} E\uutilde^{(1)}\Big)
\end{align*}
be the residual of the approximation $\uutilde^{(1)}$.
Hence, $ \delta=\uu-\uutilde^{(1)} $ solves the problem
\begin{subequations}
\begin{align}
 \label{PDE.delta.a}
\pt \delta &=  
- A(\partial)\delta
- \frac{1}{\epsi} E\delta
+ \epsi\left[T(\uu,\uu,\uu)-T(\uutilde^{(1)},\uutilde^{(1)},\uutilde^{(1)})\right]
+ R,
\\
 \label{PDE.delta.b}
\delta(0) &= 0.
\end{align}
\end{subequations}
Next, we investigate the structure of the residual.
By \eqref{Ansatz.SVEA}, the approximation $\uutilde^{(1)}$ can be expressed as 
\begin{align*}
\uutilde^{(1)}(t,x) &= \ee^{\ii (\kappa \cdot x - \w t)/\epsi} u_1(t,x) 
+ \ee^{-\ii (\kappa \cdot x - \w t)/\epsi} u_{-1}(t,x)
= \sum_{j\in\J^{(1)}} \ee^{\ii j(\kappa \cdot x - \w t)/\epsi} u_j(t,x).
\end{align*}
Substituting this into the left-hand side of \eqref{PDE.uu} 
and using \eqref{SVEA.a} yields
\begin{align}
\nonumber
& \pt \uutilde^{(1)}(t,x) +A(\partial)\uutilde^{(1)}(t,x) +\frac{1}{\epsi} E\uutilde^{(1)}(t,x)
\\
\nonumber
&\quad =
 \sum_{j\in\J^{(1)}} 
 \ee^{\ii j (\kappa \cdot x - \w t)/\epsi}
\Big(
 \pt u_j(t,x) + \tfrac{\ii}{\epsi}\L(j\w,j\kappa)u_j(t,x) + A(\partial)u_j(t,x) 
\Big) 
 \\
 &\quad =
 \epsi
 \sum_{j\in\J^{(1)}} 
 \sum_{\#J=j} 
 \ee^{\ii j (\kappa \cdot x - \w t)/\epsi}
 T(u_{j_1},u_{j_2},u_{j_3})(t,x),
 \label{R.sum.1}
 \end{align}
whereas on the right-hand side of \eqref{PDE.uu} we obtain
\begin{align}
\nonumber
\epsi T(\uutilde^{(1)},\uutilde^{(1)},\uutilde^{(1)})(t,x)
&= 
\epsi 
\sum_{J\in(\J^{(1)})^3} 
 \ee^{\ii \#J (\kappa \cdot x - \w t)/\epsi}
T(u_{j_1},u_{j_2},u_{j_3})(t,x)
\\
&= 
\epsi
\sum_{\substack{j \text{ odd} \\ |j|\leq 3}}
\sum_{\#J=j} 
 \ee^{\ii j (\kappa \cdot x - \w t)/\epsi}
T(u_{j_1},u_{j_2},u_{j_3})(t,x).
\label{R.sum.2}
\end{align}
The only difference between \eqref{R.sum.1} and \eqref{R.sum.2} is that the terms with $j=\pm3$
are missing in \eqref{R.sum.1}. These terms are exactly the higher harmonics which were omitted in the derivation of \eqref{PDE.mfe} and hence of the SVEA.
The equations \eqref{R.sum.1} and \eqref{R.sum.2} yield the representation
\begin{align*}
R(t,x) &= \epsi \sum_{j\in\{\pm 3\}}
\sum_{\#J=j} 
 \ee^{\ii j (\kappa \cdot x - \w t)/\epsi}
 T(u_{j_1},u_{j_2},u_{j_3})(t,x)
\end{align*}
of the residual.

Since \eqref{Theorem.error.bound.SVEA.assertion01} is equivalent to
\begin{align}
\label{Theorem.error.bound.SVEA.assertion01.deltahat}
 \sup_{t\in[0,\tend/\epsi]} \| \widehat{\delta}(t) \|_{L^1}
 &\leq C\epsi^2,
\end{align}
we need an evolution equation for $\widehat{\delta} = \F \delta$.
In Fourier space, \eqref{PDE.delta.a} reads
\begin{align} 
 \label{PDE.delta.Fourier}
\pt \widehat{\delta}(t,k) =  
-\big(\ii A(k)+\tfrac{1}{\epsi} E\big)
\widehat{\delta}(t,k) 
+ \epsi \G\big(\F\uu,\F\uutilde^{(1)} \big)(t,k)
+ \widehat{R}(t,k) 
\end{align}
with 
\begin{align}
\nonumber
\G\big(\F\uu,\F\uutilde^{(1)}\big)&=
\T(\F\uu,\F\uu,\F\uu)-\T\left(\F\uutilde^{(1)},\F\uutilde^{(1)},\F\uutilde^{(1)}\right),
\\
\nonumber
\widehat{R}(t,k) 
&= 
\epsi
\sum_{j\in\{\pm 3\}} \sum_{\#J=j} 
\F\left(
 T(u_{j_1},u_{j_2},u_{j_3})
 \ee^{\ii j \kappa \cdot x /\epsi} \right)(t,k) \ee^{-\ii j \w t/\epsi} 
 \\
 \label{Def.Rhat}
 &= 
 \epsi 
 \sum_{j\in\{\pm 3\}}  \sum_{\#J=j} 
 \T(\uhat_{j_1},\uhat_{j_2},\uhat_{j_3})(t,k-\tfrac{j\kappa}{\epsi})
\ee^{-\ii j \w t/\epsi},  
\end{align}
and with $\T$ defined in \eqref{Def.T.Fourier}.

\paragraph{Step 2.} 
In this step, we identify the most challenging part of the proof of \eqref{Theorem.error.bound.SVEA.assertion01}.
For this purpose, we apply Duhamel's formula to \eqref{PDE.delta.Fourier}
and use that $\widehat{\delta}(0,k) = 0$ to obtain
\begin{align}
\nonumber
\widehat{\delta}(t,k) 
&=  
\epsi \intl_0^t
\exp\big((s-t)\big(\ii A(k)+\tfrac{1}{\epsi} E\big)\big)
\G\big(\F\uu(s),\F\uutilde^{(1)}(s) \big)(k) \; \dd s
\\[-3mm]
\label{Theorem.error.bound.SVEA.01}
\\[-3mm]
\nonumber
& \qquad +
\intl_0^t
\exp\big((s-t)\big(\ii A(k)+\tfrac{1}{\epsi} E\big)\big)
\widehat{R}(s,k) \; \dd s.
\end{align}
Our goal is to prove \eqref{Theorem.error.bound.SVEA.assertion01.deltahat} via Gronwall's lemma, which requires suitable bounds for the two terms on the right-hand side of \eqref{Theorem.error.bound.SVEA.01}.
For every $k\in\IR^d$ the matrix $\ii A(k)+E/\epsi$ is skew-Hermitian, and hence
$\exp\left(t\left(\ii A(k)+E/\epsi  \right) \right)$ is unitary for every $t \in \IR$.
The first term on the right-hand side of \eqref{Theorem.error.bound.SVEA.01} can thus be bounded in $L^1$ by
\begin{align}
\nonumber
&\epsi \intl_0^t
 \intl_{\IR^d}
 \big| \exp\big((s-t)\big(\ii A(k)+\tfrac{1}{\epsi} E\big)\big) 
 \G\big(\F\uu(s),\F\uutilde^{(1)}(s) \big)(k) \big|_2 \; \dd k\; \dd s
\\
\nonumber
& \quad =
 \epsi \intl_0^t
 \intl_{\IR^d}
\big| \G\big(\F\uu(s),\F\uutilde^{(1)}(s) \big)(k) \big|_2 \; \dd k\; \dd s
 \\
& \quad  \leq 
3 C_\T C_\uu^2 \; 
\epsi \intl_0^t
\| \widehat{\delta}(s) \|_{L^1} \; \dd s.
\label{Theorem.error.bound.SVEA.02}
\end{align}
The last step follows from \eqref{Bound.T.difference.Fourier} and \eqref{Def.Cboldu}.
Now suppose that for the second term of \eqref{Theorem.error.bound.SVEA.01} the inequality
\begin{align}
 \label{Theorem.error.bound.SVEA.crucial}
 \sup_{t\in[0,\tend/\epsi]}
 \Big\|
 \intl_0^t
\exp\big((s-t)\big(\ii A(\cdot)+\tfrac{1}{\epsi} E\big)\big)
\widehat{R}(s) \; \dd s
\Big\|_{L^1} \leq C \epsi^2
\end{align}
holds. Then
it follows from \eqref{Theorem.error.bound.SVEA.01}, \eqref{Theorem.error.bound.SVEA.02}, and \eqref{Theorem.error.bound.SVEA.crucial} that
\begin{align*}
\| \widehat{\delta}(t) \|_{L^1} 
&\leq 
C C_\uu^2 \; \epsi \intl_0^t
\| \widehat{\delta}(s) \|_{L^1} \; \dd s
+ C \epsi^2,
\end{align*}
and applying Gronwall's lemma yields the desired inequality \eqref{Theorem.error.bound.SVEA.assertion01.deltahat}
with a constant which depends on $C_\uu$ and $\tend$.

The central task is thus to prove \eqref{Theorem.error.bound.SVEA.crucial}.
Equation \eqref{Def.Rhat} shows that $ \big\| \widehat{R}(s) \big\|_{L^1} = \Ord{\epsi} $,
but straightforward estimates yield only
\begin{align*}
 \sup_{t\in[0,\tend/\epsi]}
 \Big\|
 \intl_0^t
\exp\big((s-t)\big(\ii A(\cdot)+\tfrac{1}{\epsi} E\big)\big)
\widehat{R}(s) \; \dd s
\Big\|_{L^1} 
&\leq
 \sup_{t\in[0,\tend/\epsi]}
 \intl_0^t \big\| \widehat{R}(s) \big\|_{L^1} \; \dd s
 \\
&\leq 
\frac{\tend}{\epsi} \sup_{t\in[0,\tend/\epsi]} \big\| \widehat{R}(s) \big\|_{L^1}
\leq C.
\end{align*}
Compared to this simple bound, we have to gain a factor of $\epsi^2$. 
This is where the real work starts.

\paragraph{Step 3.}
In this step, we express the integral term from \eqref{Theorem.error.bound.SVEA.crucial} in an appropriate way. 
We use \eqref{Def.L}, \eqref{Def.L.j}, and \eqref{Def.Rhat} to obtain
\begin{align*} 
& \intl_0^t
\exp\big((s-t)\big(\ii A(k)+\tfrac{1}{\epsi} E\big)\big)
\widehat{R}(s,k) \; \dd s
\\ 
&\quad = \epsi 
\intl_0^t
\exp\big((s-t)\big(\ii A(k)+\tfrac{1}{\epsi} E\big)\big) 
 \sum_{j \in\{\pm 3\}}  \sum_{\#J=j} 
 \T(\uhat_{j_1},\uhat_{j_2},\uhat_{j_3})(s,k-\tfrac{j\kappa}{\epsi})
\ee^{-\ii j \w s/\epsi} \; \dd s
\\ 
&\quad =
 \epsi \ee^{-\ii j \w t/\epsi}
 \sum_{j \in\{\pm 3\}}  \sum_{\#J=j} 
\intl_0^t
\exp\left(\tfrac{\ii}{\epsi}(s-t)\L(j\w,\epsi k)\right) 
\T(\uhat_{j_1},\uhat_{j_2},\uhat_{j_3})(s,k-\tfrac{j\kappa}{\epsi})
 \; \dd s
\\ 
&\quad =
 \epsi \ee^{-\ii j \w t/\epsi}
 \sum_{j \in\{\pm 3\}}  \sum_{\#J=j} 
\intl_0^t
\exp\left(\tfrac{\ii}{\epsi}(s-t)\L_j(\epsi k')\right) 
 \T(\uhat_{j_1},\uhat_{j_2},\uhat_{j_3})(s,k')
 \; \dd s
\end{align*}
with the shifted variable $k' = k-\tfrac{j\kappa}{\epsi}$.
In order to keep the notation simple, we write again $k$ instead of $k'$ in the following.
Since later we integrate over $k$, the difference between $k$ and $k'$ does not really matter.

With \eqref{Def.F.uhat} and \eqref{Def.TT} we can represent the integrand as
\begin{align} \notag
\exp\left(\tfrac{\ii}{\epsi}(s-t)\L_j(\epsi k)\right) \T(\uhat_{j_1},\uhat_{j_2},\uhat_{j_3})(s,k)
&= \exp\left(\tfrac{\ii}{\epsi}(s-t)\L_j(\epsi k)\right) \TT_{j,\epsi}^*(s,k) \FF(s,\uhat_1,J)(k)
\\ \notag
&=
\exp\left(-\tfrac{\ii t}{\epsi}\L_j(\epsi k)\right) \Psi_{j}(\epsi k) \FF(s,\uhat_1,J)(k)
\\ \notag
&=
\TT_{j,\epsi}^*(t,k) \FF(s,\uhat_1,J)(k).
\end{align}
Since $\TT_{j,\epsi}^*(t)$ is unitary and does not depend on $s$, it follows that 
the term which appears on the left-hand side of \eqref{Theorem.error.bound.SVEA.crucial} can be bounded by
\begin{align} 
\label{Theorem.error.bound.SVEA.03}
 \Big\|
 \intl_0^t
\exp\big((s-t)\big(\ii A(\cdot)+\tfrac{1}{\epsi} E\big)\big)
\widehat{R}(s) \; \dd s
\Big\|_{L^1} 
&\leq
 \epsi 
 \sum_{j \in\{\pm 3\}}  \sum_{\#J=j} 
\Big\| \intl_0^t \FF(s,\uhat_1,J)
\; \dd s \Big\|_{L^1}.
\end{align}

\paragraph{Step 4.}
The goal in this and the following steps is to prove that
\begin{align} 
\label{Theorem.error.bound.SVEA.main.difficulty}
 \sup_{t\in[0,\tend/\epsi]}
\sum_{j \in\{\pm 3\}}  \sum_{\#J=j} 
\Big\| \intl_0^t \FF(s,\uhat_1,J)
\; \dd s \Big\|_{L^1} \leq C \epsi.
\end{align}
If \eqref{Theorem.error.bound.SVEA.main.difficulty} holds, then the crucial inequality \eqref{Theorem.error.bound.SVEA.crucial} follows 
via\footnote{Note that the right-hand side of \eqref{Theorem.error.bound.SVEA.03} contains a factor $\epsi$, which was omitted on the left-hand side of \eqref{Theorem.error.bound.SVEA.main.difficulty}.
}
\eqref{Theorem.error.bound.SVEA.03},
which then completes the proof of \eqref{Theorem.error.bound.SVEA.assertion01}.
The sum in \eqref{Theorem.error.bound.SVEA.main.difficulty} is taken over multi-indices $J\in(\J^{(1)})^3=\{1,-1\}^3$ with $\#J=j\in\{3,-3\}$. There are only two possibilities, namely $ J=(1,1,1) $, $ j=3 $ and
$ J=-(1,1,1) $, $ j=-3. $ 
Since both cases can be treated \emph{mutatis mutandis}, we will only consider the first one,
i.e.\ $ J=(1,1,1) $, $ j=3 $, and thus
\begin{align*} 
\FF(s,\uhat_1,J) = \TT_{3,\epsi}(s) \T\big(\uhat_1,\uhat_1,\uhat_1\big)(s).
\end{align*}
We have to show that
\begin{align} 
\label{Main.difficulty.0}
 \sup_{t\in[0,\tend/\epsi]}
\Big\| \intl_0^t \TT_{3,\epsi}(s) \T\big(\uhat_1,\uhat_1,\uhat_1\big)(s)
\; \dd s \Big\|_{L^1} \leq C \epsi.
\end{align}
In order to use Proposition~\ref{Proposition01.SVEA}, we decompose the nonlinearity into eight parts
\begin{align*} 
\T\big(\uhat_1,\uhat_1,\uhat_1\big)
&=
\T\big(\P_\epsi \uhat_1,\P_\epsi\uhat_1,\P_\epsi\uhat_1\big)
+ \T\big(\P_\epsi \uhat_1,\P_\epsi\uhat_1,\P_\epsi^\perp\uhat_1\big)
+ \T\big(\P_\epsi \uhat_1,\P_\epsi^\perp\uhat_1,\P_\epsi\uhat_1\big)
\\ 
&\quad
+ \T\big(\P_\epsi^\perp \uhat_1,\P_\epsi\uhat_1,\P_\epsi\uhat_1\big)
+ \T\big(\P_\epsi \uhat_1,\P_\epsi^\perp\uhat_1,\P_\epsi^\perp\uhat_1\big)
+ \T\big(\P_\epsi^\perp \uhat_1,\P_\epsi\uhat_1,\P_\epsi^\perp\uhat_1\big)
\\ 
&\quad
+ \T\big(\P_\epsi^\perp \uhat_1,\P_\epsi^\perp\uhat_1,\P_\epsi\uhat_1\big)
+ \T\big(\P_\epsi^\perp \uhat_1,\P_\epsi^\perp\uhat_1,\P_\epsi^\perp\uhat_1\big).
\end{align*}
The last four terms are those where $\P_\epsi^\perp\uhat_1$ appears in at least two of the three arguments of 
$ \T(\cdot,\cdot,\cdot) $. 
These terms are $ \Ord{\epsi^2} $ because of Proposition~\ref{Proposition01.SVEA}, 
and their contribution to the left-hand side of \eqref{Theorem.error.bound.SVEA.main.difficulty} 
can be estimated in a straightforward way, for example
\begin{align*}
& \Big\| \intl_0^t 
\TT_{3,\epsi}(s) \T\big(\P_\epsi^\perp \uhat_1,\P_\epsi^\perp \uhat_1,\P_\epsi\uhat_{1}\big)(s)
\; \dd s \Big\|_{L^1}
\\
&\leq
C_\T \, \intl_0^t \Big(
\| \P_\epsi^\perp \uhat_1(s) \|_{L^1} 
\, \| \P_\epsi^\perp \uhat_1(s) \|_{L^1}
\, \| \P_\epsi\uhat_1(s) \|_{L^1}
\Big)
\; \dd s
\leq 
C t \epsi^2 
\leq 
C \tend \epsi.
\end{align*}
For the first four parts of $\T\big(\uhat_1,\uhat_1,\uhat_1\big)$ the analysis is much more involved. We have to prove that
\begin{align}
\label{Main.difficulty.2}
\Big\| \intl_0^t \TT_{3,\epsi}(s) \T\big(\P_\epsi^\perp \uhat_1,\P_\epsi\uhat_1,\P_\epsi\uhat_1 \big)(s) \; \dd s \Big\|_{L^1}
\leq C \epsi,
\\
\label{Main.difficulty.1}
\Big\| \intl_0^t \TT_{3,\epsi}(s) \T\big(\P_\epsi \uhat_1,\P_\epsi\uhat_1,\P_\epsi\uhat_1 \big)(s) \; \dd s \Big\|_{L^1}
\leq C \epsi,
\end{align}
because bounds for the two terms involving 
$\T\big(\P_\epsi\uhat_1,\P_\epsi^\perp \uhat_1,\P_\epsi\uhat_1 \big)$ and
$\T\big(\P_\epsi\uhat_1,\P_\epsi\uhat_1,\P_\epsi^\perp \uhat_1 \big)$ can be shown in the same way as 
\eqref{Main.difficulty.2}.

\paragraph{Step 5.}
In this step we prove \eqref{Main.difficulty.2}. 
To accomplish this, we have to identify the oscillatory ``parts'' of the integrand.
We use that \eqref{Relation.between.projections}, \eqref{Def.z}, and \eqref{Def.TT} 
yield the representation
\begin{align}
\P_\epsi^\perp \uhat_1(t,k) 
&= \TT_{1,\epsi}^*(t,k) \Pperp z_1(t,k)
= 
\label{representation.Pperpuhat.1}
\sum_{\ell=2}^n
\exp\big(-\tfrac{\ii t}{\epsi}\lambda_{1\ell}(\epsi k)\big)z_{1\ell}(t,k)\psi_{1\ell}(\epsi k),
\end{align}
where again $ \lambda_{1\ell}(\epsi k) $ is the $\ell$-th eigenvalue of $ \L_1(\epsi k) $ and
$ \psi_{1\ell}(\epsi k) $ is the corresponding eigenvector, as defined in \eqref{eigendecomposition}.
Combining \eqref{representation.Pperpuhat.1} with \eqref{Def.TT} and \eqref{Def.T.Fourier} results in
\begin{align*}
\notag
& 
\intl_0^t
\TT_{3,\epsi}(s,k) \T\big(\P_\epsi^\perp \uhat_1,\P_\epsi\uhat_1,\P_\epsi\uhat_1 \big)(s,k)
\; \dd s 
\\
\notag
&=
\frac{1}{(2\pi)^d} 
\intl_0^t
\exp\big(\tfrac{\ii s}{\epsi}\Lambda_3(\epsi k)\big)  \Psi_3^*(\epsi k)
\intl_{\#K = k} T\left(\P_\epsi^\perp \uhat_1(s,k^{(1)}),\P_\epsi\uhat_1(s,k^{(2)}),\P_\epsi\uhat_1(s,k^{(3)})\right)
\; \dd K
\; \dd s 
\\
% \label{Step5.Eq.1}
&= 
\sum_{\ell=2}^n \intl_{\#K = k}
\intl_0^t
\exp\big(\tfrac{\ii s}{\epsi}\big[\Lambda_3(\epsi k)-\lambda_{1\ell}(\epsi k^{(1)})I\big]\big)
f_{\epsi,\ell}(s,K)
\; \dd s 
\; \dd K 
\end{align*}
with the shorthand notation from \eqref{Shorthand.notation.K}, and with 
\begin{align*}
f_{\epsi,\ell}(s,K)&=
\frac{1}{(2\pi)^d} 
\Psi_3^*(\epsi k)
 T\left(z_{1\ell}(s,k^{(1)})\psi_{1\ell}(\epsi k^{(1)}),
 \P_\epsi\uhat_1(s,k^{(2)}),\P_\epsi\uhat_1(s,k^{(3)})\right),
 \quad \#K = k.
\end{align*}
Taking the norm yields
\begin{align}
\notag
& \Big\| \intl_0^t 
\TT_{3,\epsi}(s) \T\big(\P_\epsi^\perp \uhat_1,\P_\epsi \uhat_1,\P_\epsi\uhat_{1}\big)(s)
\; \dd s \Big\|_{L^1}
\\
\notag
&=
\intl_{\IR^d}
\Big|
\intl_0^t 
\TT_{3,\epsi}(s,k) \T\big(\P_\epsi^\perp \uhat_1,\P_\epsi \uhat_1,\P_\epsi\uhat_{1}\big)(s,k)
\; \dd s
\Big|_2 \; \dd k
\\
\label{Step5.Eq.2}
&\leq
\sum_{\ell=2}^n 
\intl_{\IR^d} \intl_{\#K = k}
\Big|
\intl_0^t
\exp\big(\tfrac{\ii s}{\epsi}\big[\Lambda_3(\epsi k)-\lambda_{1\ell}(\epsi k^{(1)})I\big]\big)
f_{\epsi,\ell}(s,K)
\; \dd s 
\Big|_2 \; \dd K \; \dd k.
\end{align}
Now we focus on the inner integral. 
The exponential function
$ s \mapsto \exp\big(\tfrac{\ii s}{\epsi}
\big[\Lambda_3(\epsi k)-\lambda_{1\ell}(\epsi k^{(1)})I\big]\big) $ 
in \eqref{Step5.Eq.2}
oscillates if all diagonal entries of the diagonal matrix 
$\Lambda_3(\epsi k)-\lambda_{1\ell}(\epsi k^{(1)})I$
are bounded away from zero,
but we cannot expect this to be true for all $ k, k^{(1)}\in\IR^d$.
For this reason, we define
\begin{align*}
\Delta_\ell(\theta,\theta^{(1)}) &= \Lambda_3(\theta)-\lambda_{1\ell}(\theta^{(1)})I
\qquad \text{for } \theta, \theta^{(1)} \in \IR^d,
\\
g_{\epsi,\ell}(s,K) &= 
\exp\left(\frac{\ii s}{\epsi}\big[
\Delta_\ell(\epsi k,\epsi k^{(1)}) - \Delta_\ell(0,0)\big]
\right)
f_{\epsi,\ell}(s,K)
\end{align*}
and reformulate the inner integral in \eqref{Step5.Eq.2} as
\begin{align}
\label{Step5.Eq.3}
\Big| \intl_0^t
\exp\big(\tfrac{\ii s}{\epsi}\big[\Lambda_3(\epsi k)-\lambda_{1\ell}(\epsi k^{(1)})I\big]\big)
f_{\epsi,\ell}(s,K)
\; \dd s \Big|_2
=
\Big| \intl_0^t
\exp\left(\frac{\ii s}{\epsi}\Delta_\ell(0,0)\right)
g_{\epsi,\ell}(s,K)
\; \dd s \Big|_2.
\end{align}
By Assumption~\ref{Ass:Delta.lambda.j3} the diagonal matrix 
\begin{align*}
\Delta_\ell(0,0) = \Lambda_3(0)-\lambda_{1\ell}(0)I = 
\diag\Big(\lambda_{31}(0)-\lambda_{1\ell}(0), \ldots , \lambda_{3n}(0)-\lambda_{1\ell}(0)\Big)
\end{align*}
is regular for all $\ell$. 
Hence, we can now integrate by parts to obtain
\begin{align}
\notag
&
\Big| \intl_0^t
\exp\left(\frac{\ii s}{\epsi}\Delta_\ell(0,0)\right)
g_{\epsi,\ell}(s,K)
\; \dd s \Big|_2
\\
\notag
&=
\Big| 
\frac{\epsi}{\ii} \Delta_\ell(0,0)^{-1}
\left(
\exp\left(\frac{\ii t}{\epsi}\Delta_\ell(0,0)\right)
g_{\epsi,\ell}(t,K)
-
g_{\epsi,\ell}(0,K)
\right)
\Big|_2
\\
\notag
&\quad
+
\Big| \frac{\epsi}{\ii} \Delta_\ell(0,0)^{-1}
\intl_0^t
\exp\left(\frac{\ii s}{\epsi}\Delta_\ell(0,0)\right)
\pt g_{\epsi,\ell}(s,K)
\; \dd s
\Big|_2
\\
\label{Step5.Eq.4}
&\leq
C \epsi \Big(| g_{\epsi,\ell}(t,K) |_2 + | g_{\epsi,\ell}(0,K) |_2 \Big)
+ 
C \epsi 
\intl_0^t
\Big| \pt g_{\epsi,\ell}(s,K) \Big|_2
\; \dd s.
\end{align}
By definition of $g_{\epsi,\ell}$, we have
\begin{align*}
\pt g_{\epsi,\ell}(s,K) 
&=
\frac{\ii}{\epsi}\big[\Delta_\ell(\epsi k,\epsi k^{(1)}) - \Delta_\ell(0,0)\big]
g_{\epsi,\ell}(s,K) 
\\
& \quad
+
\exp\left(\frac{\ii s}{\epsi}\big[
\Delta_\ell(\epsi k,\epsi k^{(1)}) - \Delta_\ell(0,0)\big]
\right)
\pt f_{\epsi,\ell}(s,K),
\end{align*}
and since $ | g_{\epsi,\ell}(s,K) |_2 = | f_{\epsi,\ell}(s,K) |_2 $
this yields
\begin{align*}
\big| \pt g_{\epsi,\ell}(s,K) \big|_2
&\leq
\frac{C}{\epsi}\big|\Delta_\ell(\epsi k,\epsi k^{(1)}) - \Delta_\ell(0,0)\big|_2
|f_{\epsi,\ell}(s,K)|_2
+
| \pt f_{\epsi,\ell}(s,K) |_2.
\end{align*}
With \eqref{Step5.Eq.4}, \eqref{Step5.Eq.3} and \eqref{Step5.Eq.2} we infer that
\begin{align*}
\Big\| \intl_0^t 
\TT_{3,\epsi}(s) \T\big(\P_\epsi^\perp \uhat_1,\P_\epsi \uhat_1,\P_\epsi\uhat_{1}\big)(s)
\; \dd s \Big\|_{L^1}
&\leq
C\epsi \Big(X_1(t,\epsi) + X_2(t,\epsi) + X_3(t,\epsi)\Big)
\end{align*}
with
\begin{align*}
X_1(t,\epsi) &= \sum_{\ell=2}^n 
\intl_{\IR^d} \intl_{\#K = k}
\Big(| g_{\epsi,\ell}(t,K) |_2 + | g_{\epsi,\ell}(0,K) |_2 \Big)
\; \dd K \; \dd k,
\\
X_2(t,\epsi) &=
\sum_{\ell=2}^n 
\intl_{\IR^d} \intl_{\#K = k}
\intl_0^t
\frac{1}{\epsi}\big|\Delta_\ell(\epsi k,\epsi k^{(1)}) - \Delta_\ell(0,0)\big|_2
|f_{\epsi,\ell}(s,K)|_2
\; \dd s
\; \dd K \; \dd k,
\\
X_3(t,\epsi) &=
\sum_{\ell=2}^n 
\intl_{\IR^d} \intl_{\#K = k}
\intl_0^t
| \pt f_{\epsi,\ell}(s,K) |_2
\; \dd s
\; \dd K \; \dd k.
\end{align*}
In order to complete the proof of \eqref{Main.difficulty.2} we have to show that 
$X_1(t,\epsi)$, $X_2(t,\epsi)$, and $X_3(t,\epsi)$ are uniformly bounded in 
$\epsi \in(0,1]$ and $t\in [0,\tend/\epsi]$. 
For $X_2(t,\epsi)$ and $X_3(t,\epsi)$, this is not obvious because of the integration over the possibly long time interval $[0,t]$ with $ t \leq \tend/\epsi $.
We use that
\begin{align*}
\sum_{\ell=2}^n | f_{\epsi,\ell}(s,K) |_2 
&\leq 
C \sum_{\ell=2}^n \Big| 
 T\Big(z_{1\ell}(s,k^{(1)})\psi_{1\ell}(\epsi k^{(1)}),
 \P_\epsi\uhat_1(s,k^{(2)}),\P_\epsi\uhat_1(s,k^{(3)})\Big)
\Big|_2
\\
&\leq
C 
\big| \P_\epsi^\perp \uhat_{1}(s,k^{(1)}) \big|_2
\big| \P_\epsi\uhat_1(s,k^{(2)}) \big|_2
\big| \P_\epsi\uhat_1(s,k^{(3)}) \big|_2
\end{align*}
holds, because of the normalization $ |\psi_{1\ell}(\epsi k^{(1)})|_2 = 1$
and the fact that 
\begin{align*}
\sum_{\ell=2}^n | z_{1\ell}(s,k^{(1)}) | = | \Pperp z_{1}(s,k^{(1)}) |_1
\leq C | \Pperp z_{1}(s,k^{(1)}) |_2 = C \big| \P_\epsi^\perp \uhat_{1}(s,k^{(1)}) \big|_2
\end{align*}
by \eqref{Norm.uhat.z.Pperp}.
With $ | g_{\epsi,\ell}(s,K) |_2 = | f_{\epsi,\ell}(s,K) |_2 $ this implies that
\begin{align*}
\notag
& \sum_{\ell=2}^n 
\intl_{\IR^d} \intl_{\#K = k}
| g_{\epsi,\ell}(s,K) |_2 
\; \dd K \; \dd k
\\
\notag
&\leq
C 
\intl_{\IR^d} \intl_{\#K = k}
\Big(
\big| \P_\epsi^\perp \uhat_{1}(s,k^{(1)}) \big|_2
\big| \P_\epsi\uhat_1(s,k^{(2)}) \big|_2
\big| \P_\epsi\uhat_1(s,k^{(3)}) \big|_2
\Big)
\; \dd K \; \dd k
\\
\notag
&=
C 
\intl_{\IR^d} 
\big| \P_\epsi^\perp \uhat_{1}(s,k^{(1)}) \big|_2
\; \dd k^{(1)}
\intl_{\IR^d} 
\big| \P_\epsi\uhat_1(s,k^{(2)}) \big|_2
\; \dd k^{(2)}
\intl_{\IR^d} 
\big| \P_\epsi\uhat_1(s,k^{(3)}) \big|_2
\; \dd k^{(3)}
\\
\notag
&=
C \| \P_\epsi^\perp \uhat_{1}(s) \|_{L^1}
\| \P_\epsi\uhat_1(s) \|_{L^1}^2
\\
% \label{Step5.Eq.5}
&\leq C \epsi
\end{align*}
for all $s\in [0,\tend/\epsi]$ due to Proposition~\ref{Proposition01.SVEA}.
This shows in particular that $ X_1(t,\epsi) $ is uniformly bounded\footnote{In fact, we have 
even shown that $ X_1(t,\epsi) \leq C \epsi $ for all $\epsi \in(0,1]$ and $t\in [0,\tend/\epsi]$.
} % footnote
in $\epsi \in(0,1]$ and $t\in [0,\tend/\epsi]$. 

For $X_2(t,\epsi)$ we use that the Lipschitz continuity \eqref{Ass:lambda.Lipschitz} of the eigenvalues 
yields
\begin{align*}
\big|\Delta_\ell(\epsi k,\epsi k^{(1)}) - \Delta_\ell(0,0)\big|_2
&\leq
\big| \Lambda_3(\epsi k)-\Lambda_3(0) \big|_2
+ \big| \lambda_{1\ell}(\epsi k^{(1)})-\lambda_{1\ell}(0) \big|
\\
&\leq
C \epsi (|k|_1 + |k^{(1)}|_1),
\end{align*}
and the $\epsi$ in the second line compensates the factor $1/\epsi$ in $X_2(t,\epsi)$.
For $ K=(k^{(1)},k^{(2)},k^{(3)}) $ with 
$ k = \#K=k^{(1)}+k^{(2)}+k^{(3)} $, we have that
$ |k|_1 \leq |k^{(1)}|_1 + |k^{(2)}|_1 + |k^{(3)}|_1 $.
Hence, it follows that
\begin{align*}
X_2(t,\epsi) &=
\sum_{\ell=2}^n 
\intl_{\IR^d} \intl_{\#K = k}
\intl_0^t
\frac{1}{\epsi}\big|\Delta_\ell(\epsi k,\epsi k^{(1)}) - \Delta_\ell(0,0)\big|_2
|f_{\epsi,\ell}(s,K)|_2
\; \dd s
\; \dd K \; \dd k
\\
&\leq
C
\sum_{\ell=2}^n 
\intl_{\IR^d} \intl_{\#K = k}
\intl_0^t
(|k^{(1)}|_1 + |k^{(2)}|_1 + |k^{(3)}|_1)
|f_{\epsi,\ell}(s,K)|_2
\; \dd s
\; \dd K \; \dd k
\end{align*}
and proceeding as before yields
\begin{align*}
X_2(t,\epsi) 
&\leq
C \sum_{\mu=1}^d \intl_0^t
\Big(
\| D_\mu \P_\epsi^\perp \uhat_{1}(s) \|_{L^1}
\| \P_\epsi\uhat_1(s) \|_{L^1}^2
\\
&
\hs{20} + 2
\| \P_\epsi^\perp \uhat_{1}(s) \|_{L^1}
\| D_\mu \P_\epsi\uhat_1(s) \|_{L^1}
\| \P_\epsi\uhat_1(s) \|_{L^1}
\Big)
\; \dd s
\\
&\leq
C \sum_{\mu=1}^d \intl_0^t
\| D_\mu \P_\epsi^\perp \uhat_{1}(s) \|_{L^1}
\; \dd s
+
C \intl_0^t
\| \P_\epsi^\perp \uhat_{1}(s) \|_{L^1}
\; \dd s.
\end{align*}
Since both integrands are $\Ord{\epsi}$ according to 
Propositions~\ref{Proposition01.SVEA} and \ref{Proposition02.SVEA}, respectively,
the right-hand side is uniformly bounded for $t\in [0,\tend/\epsi]$.

In a similar way, one can show that
\begin{align*}
X_3(t,\epsi) &=
\sum_{\ell=2}^n 
\intl_{\IR^d} \intl_{\#K = k}
\intl_0^t
| \pt f_{\epsi,\ell}(s,K) |_2
\; \dd s
\; \dd K \; \dd k
\\
&\leq
C \sum_{\ell=2}^n \intl_0^t
\Big(
\| \pt \Pperp z_{1}(s) \|_{L^1}
\| \P_\epsi\uhat_1(s) \|_{L^1}^2
\\
& \hs{20} + 2
\| \Pperp z_{1}(s) \|_{L^1}
\| \pt \P_\epsi\uhat_1(s) \|_{L^1}
\| \P_\epsi\uhat_1(s) \|_{L^1}
\Big)
\; \dd s.
\end{align*}
Since $ \| \Pperp z_{1}(s) \|_{L^1} = \| \P_\epsi^\perp \uhat_{1}(s) \|_{L^1} \leq C\epsi $ by Proposition~\ref{Proposition01.SVEA},
since $\| \pt \P_\epsi\uhat_1(s) \|_{L^1}$ is uniformly bounded by Lemma~\ref{Lemma.bound.uhat1dot},
and since
\begin{align*}
\| \pt \Pperp z_{1}(s) \|_{L^1} 
\leq
\| \pt z_{1}(s) \|_{L^1} 
& \leq
\epsi \sum_{\#J=1} \| \T\big(\uhat_{j_1},\uhat_{j_2},\uhat_{j_3}\big)(t) \|_{L^1} 
\leq C\epsi
\end{align*}
by \eqref{zdot.in.terms.of.uhat} and \eqref{Def.F.uhat}, we conclude that
$ X_3(t,\epsi) $ is uniformly bounded, too.
We have thus shown the inequality \eqref{Main.difficulty.2}.

\paragraph{Step 6.}

In this step, we prove \eqref{Main.difficulty.1}.
For the proof of \eqref{Main.difficulty.2} in the previous step, it was crucial that
$ \P_\epsi^\perp \uhat_1 $ appears in one of the arguments of $ \T $, 
because this allowed us to use Propositions~\ref{Proposition01.SVEA} and \ref{Proposition02.SVEA}.
In \eqref{Main.difficulty.1}, however, this is not possible, because all three arguments of $ \T $ are 
$ \P_\epsi \uhat_1 $ instead of $ \P_\epsi^\perp \uhat_1 $. Hence, we have to proceed in a different way.
The crucial observation is that $\P_\epsi \uhat_1$ and thus also
$ \T\big(\P_\epsi \uhat_1,\P_\epsi\uhat_1,\P_\epsi\uhat_1 \big) $
are non-oscillatory in the sense that the first two time derivatives of $\P_\epsi \uhat_1$
are uniformly bounded according to 
Lemmas~\ref{Lemma.bound.uhat1dot}, \ref{Lemma.bound.Dmu.uhat1dot}, and \ref{Lemma.bound.uhat1dotdot}.
The only oscillatory function on the left-hand side of \eqref{Main.difficulty.1}
is $ \TT_{3,\epsi}(s) $. The strategy is now to integrate by parts \emph{twice}, 
which generates a factor $\epsi$ each time. 
One of these factors is then used to compensate the long time interval.

We set $ \Delta_{3}(\epsi k) = \Lambda_{3}(\epsi k) - \Lambda_{3}(0) $ and 
\begin{align*} 
f_\epsi(t,k) =
\exp\big(\tfrac{\ii t}{\epsi}\Delta_{3}(\epsi k)\big) 
\Psi_{3}^*(\epsi k)\T\big(\P_\epsi \uhat_1,\P_\epsi \uhat_1,\P_\epsi \uhat_1\big)(t,k).
\end{align*}
With \eqref{Def.TT} we obtain the representation
\begin{align*}
\Big\|
\intl_0^t \TT_{3,\epsi}(s) \T\big(\P_\epsi \uhat_1,\P_\epsi\uhat_1,\P_\epsi\uhat_1 \big)(s) \; \dd s \Big\|_{L^1}
&=
\Big\|
\intl_0^t \exp\big(\tfrac{\ii s}{\epsi}\Lambda_{3}(0)\big) f_\epsi(s)
\; \dd s \Big\|_{L^1}
\end{align*}
of the left-hand side of \eqref{Main.difficulty.1}.
By Assumption~\ref{Ass:Delta.lambda.j3} the matrix 
$ \L_3(0) = \L(3\w, 3\kappa) $ and thus also $\Lambda_{3}(0)$ is invertible.
Hence, we can integrate by parts twice and obtain
\begin{align} 
\notag
& \Big\| \intl_0^t 
\exp\big(\tfrac{\ii s}{\epsi}\Lambda_{3}(0)\big) f_\epsi(s)
\; \dd s\Big\|_{L^1}
\leq 
C\epsi
\bigg[\| f_\epsi(0) \|_{L^1} + \| f_\epsi(t) \|_{L^1}
+ 
\Big\| \intl_0^t 
\exp\big(\tfrac{\ii s}{\epsi}\Lambda_{3}(0)\big)  \pt f_\epsi(s) \; \dd s\Big\|_{L^1}
\bigg]
\\
\notag
&\leq 
C\epsi \bigg[
\| f_\epsi(0) \|_{L^1} + \| f_\epsi(t) \|_{L^1}
+ 
\epsi \| \pt f_\epsi(0) \|_{L^1} + \epsi \| \pt f_\epsi(t) \|_{L^1}
\\[-3mm]
\label{Step.6.big.bracket}
\\[-3mm]
\notag
& \hs{10}
+
\epsi
\Big\| \intl_0^t 
\exp\big(\tfrac{\ii s}{\epsi}\Lambda_{3}(0)\big)  \pt^2 f_\epsi(s) \; \dd s\Big\|_{L^1}
\bigg].
\end{align}
Now we have to show that all terms inside the big bracket $[ \; \ldots \; ]$ are uniformly bounded 
in $t \in [0,\tend/\epsi]$ and $\epsi\in(0,1]$.

As a preparatory step, we note that applying \eqref{Bound.T.Fourier}, \eqref{Pepsi.decrease.norm}, and the product rule yields
\begin{align} 
\label{Step.6.bound.T}
  \big\| \T\big(\P_\epsi \uhat_1,\P_\epsi \uhat_1,\P_\epsi \uhat_1\big)(t) \big\|_{L^1}
&\leq C_\T \| \uhat_1(t) \|_{L^1}^3,
\\
\label{Step.6.bound.T.Dmu}
  \big\| D_\mu \T\big(\P_\epsi \uhat_1,\P_\epsi \uhat_1,\P_\epsi \uhat_1\big)(t) \big\|_{L^1}
&\leq 
3 C_\T \big\| D_\mu \uhat_1(t) \big\|_{L^1} \| \uhat_1(t) \|_{L^1}^2,
\\
\label{Step.6.bound.T.DmuDnu}
  \big\| D_\mu D_\nu \T\big(\P_\epsi \uhat_1,\P_\epsi \uhat_1,\P_\epsi \uhat_1\big)(t) \big\|_{L^1}
&\leq 
C_\T \Big(6 \big\| D_\mu \uhat_1(t) \big\|_{L^1}  
\big\| D_\nu \uhat_1(t) \big\|_{L^1}  \| \uhat_1(t) \|_{L^1}
\\
\notag
& \hs{15}
+ 3 \big\| D_\mu D_\nu \uhat_1(t) \big\|_{L^1} \big\| \uhat_1(t) \big\|_{L^1}^2
\Big)
\end{align}
for all $\mu,\nu\in\{1, \ldots, d\}$.
The right-hand side of \eqref{Step.6.bound.T}, \eqref{Step.6.bound.T.Dmu}, and
\eqref{Step.6.bound.T.DmuDnu} is uniformly bounded
by Lemma~\ref{Lemma.wellposedness}\ref{Lemma.wellposedness.23}.
In a similar way, we obtain the inequalities
\begin{align} 
\label{Step.6.bound.Tdot}
\big\| \pt \T\big(\P_\epsi \uhat_1,\P_\epsi \uhat_1,\P_\epsi \uhat_1\big)(t) \big\|_{L^1}
 &\leq 
 3 C_\T \big\| \pt \P_\epsi \uhat_1(t) \big\|_{L^1}  \| \uhat_1(t) \|_{L^1}^2,
 \\
%%%
\label{Step.6.bound.Tdot.Dmu}
  \big\| D_\mu \pt \T\big(\P_\epsi \uhat_1,\P_\epsi \uhat_1,\P_\epsi \uhat_1\big)(t) \big\|_{L^1}
&\leq 
C_\T \Big(6 \big\| D_\mu \uhat_1(t) \big\|_{L^1}  
\big\| \pt \P_\epsi \uhat_1(t) \big\|_{L^1}  \| \uhat_1(t) \|_{L^1}
\\
\notag
& \hs{15}
+ 3 \big\| D_\mu \pt \P_\epsi \uhat_1(t) \big\|_{L^1} \big\| \uhat_1(t) \big\|_{L^1}^2
\Big),
\\
%%%%%
\label{Step.6.bound.Tdotdot}
\big\| \pt^2 \T\big(\P_\epsi \uhat_1,\P_\epsi \uhat_1,\P_\epsi \uhat_1\big)(t) \big\|_{L^1}
 &\leq 
 3C_\T \Big(
 \big\| \pt \P_\epsi \uhat_1(t) \big\|_{L^1}^2  \| \uhat_1(t) \|_{L^1}
 \\
 \notag
 &\hs{15}
 +
 \big\| \pt^2 \P_\epsi \uhat_1(t) \big\|_{L^1}  \| \uhat_1(t) \|_{L^1}^2
 \Big)
\end{align}
and applying Lemmas~\ref{Lemma.bound.uhat1dot}, \ref{Lemma.bound.Dmu.uhat1dot}, and \ref{Lemma.bound.uhat1dotdot}
yields uniform boundedness of the right-hand sides of
\eqref{Step.6.bound.Tdot}, \eqref{Step.6.bound.Tdot.Dmu}, and \eqref{Step.6.bound.Tdotdot}.

Since the matrix $ \exp\big(\tfrac{\ii t}{\epsi}\Delta_{3}(\epsi k)\big) \Psi_{3}^*(\epsi k) $ is unitary,
\eqref{Step.6.bound.T} implies that $ \| f_\epsi(t) \|_{L^1} $ is uniformly bounded.
Taking the time derivative of $ f_\epsi(t) $ gives
\begin{align}
\label{Step.6.formula.fdot}
\pt f_\epsi(t,k) &= f_\epsi^{[1,1]}(t,k) + f_\epsi^{[1,2]}(t,k),
\\
\notag
f_\epsi^{[1,1]}(t,k) &=
\tfrac{\ii}{\epsi}\Delta_{3}(\epsi k) f_\epsi(t,k),
\\
\notag
f_\epsi^{[1,2]}(t,k) &=
\exp\big(\tfrac{\ii t}{\epsi}\Delta_{3}(\epsi k)\big) \Psi_{3}^*(\epsi k) \pt \T\big(\P_\epsi \uhat_1,\P_\epsi \uhat_1,\P_\epsi \uhat_1\big)(t,k).
\end{align}
The fact that $ \Lambda_3$ is globally Lipschitz continuous by \eqref{Ass:lambda.Lipschitz} yields
\begin{align*} 
 |\tfrac{\ii}{\epsi}\Delta_{3}(\epsi k)|_2 = \tfrac{1}{\epsi}|\Lambda_{3}(\epsi k) - \Lambda_{3}(0)|_2 \leq C|k|_1
\end{align*}
with a constant $C$ which does not depend on $\epsi$ and $k$.
Using again that $ \exp\big(\tfrac{\ii t}{\epsi}\Delta_{3}(\epsi k)\big) \Psi_{3}^*(\epsi k) $ is a unitary
matrix gives
\begin{align} 
\notag
|f_\epsi^{[1,1]}(t)|_2 
&\leq
C|k|_1 |f_\epsi(t,k)|_2
\leq
C|k|_1 \big|\T\big(\P_\epsi \uhat_1,\P_\epsi \uhat_1,\P_\epsi \uhat_1\big)(t,k)\big|_2
\\
\notag
% \label{Step.6.bound.g1}
&= C \sum_{\mu=1}^d \big|D_\mu \T\big(\P_\epsi \uhat_1,\P_\epsi \uhat_1,\P_\epsi \uhat_1\big)(t,k)\big|_2,
\\
\label{Step.6.bound.g2}
|f_\epsi^{[1,2]}(t)|_2 
&\leq
| \pt \T\big(\P_\epsi \uhat_1,\P_\epsi \uhat_1,\P_\epsi \uhat_1\big)(t,k) |_2
\end{align}
and by combining this with 
\eqref{Step.6.formula.fdot}, \eqref{Step.6.bound.T.Dmu}, and \eqref{Step.6.bound.Tdot}  
we infer that $ \epsi \| \pt f_\epsi(t) \|_{L^1} \leq C \epsi. $
% \begin{align*} 
% \epsi \| \pt f_\epsi(t) \|_{L^1} \leq C \epsi.
% \end{align*}
This $\Ord{\epsi}$ estimate is even better than the uniform boundedness which we require at this point.
Finally, we show uniform boundedness of the integral term in \eqref{Step.6.big.bracket}.
Since $ t\in[0,\tend/\epsi] $ we can use that
\begin{align*}
% \label{Step.6.big.bracket.integral.term}
\epsi
\Big\| \intl_0^t 
\exp\big(\tfrac{\ii s}{\epsi}\Lambda_{3}(0)\big)  \pt^2 f_\epsi(s) \; \dd s\Big\|_{L^1}
&\leq 
\tend \sup_{s\in[0,\tend/\epsi]} \| \pt^2 f_\epsi(s) \|_{L^1}
\end{align*}
with
\begin{align*}
\pt^2 f_\epsi(s) &= \pt^2 \Big(\exp\big(\tfrac{\ii t}{\epsi}\Delta_{3}(\epsi k)\big) 
\Psi_{3}^*(\epsi k)\T\big(\P_\epsi \uhat_1,\P_\epsi \uhat_1,\P_\epsi \uhat_1\big)(s,k) \Big)
\\
&= f_\epsi^{[2,1]}(s,k) + 2f_\epsi^{[2,2]}(s,k) + f_\epsi^{[2,3]}(s,k),
\\
f_\epsi^{[2,1]}(s,k) &=
\left(\tfrac{\ii}{\epsi}\Delta_{3}(\epsi k)\right)^2 f_\epsi(s,k),
\\
f_\epsi^{[2,2]}(s,k) &=
\tfrac{\ii}{\epsi}\Delta_{3}(\epsi k) f_\epsi^{[1,2]}(s,k),
\\
f_\epsi^{[2,3]}(s,k) &=
\exp\big(\tfrac{\ii t}{\epsi}\Delta_{3}(\epsi k)\big) \Psi_{3}^*(\epsi k) 
\pt^2 \T\big(\P_\epsi \uhat_1,\P_\epsi \uhat_1,\P_\epsi \uhat_1\big)(s,k).
\end{align*}
Proceeding as before yields 
\begin{align*}
|f_\epsi^{[2,1]}(s,k)|_2 
&\leq 
C |k|_1^2 |f_\epsi(s,k)|_2 
\\
&=
C |k|_1^2 \big|\T\big(\P_\epsi \uhat_1,\P_\epsi \uhat_1,\P_\epsi \uhat_1\big)(s,k)\big|_2 
\\
&= 
C \sum_{\mu=1}^d \sum_{\nu=1}^d \big|D_\mu D_\nu 
\T\big(\P_\epsi \uhat_1,\P_\epsi \uhat_1,\P_\epsi \uhat_1\big)(s,k)\big|_2,
\end{align*}
and hence uniform boundedness of $ \|f_\epsi^{[2,1]}(s)\|_{L^1} $ follows from \eqref{Step.6.bound.T.DmuDnu}.
In a similar way, we obtain with \eqref{Step.6.bound.g2}
\begin{align*}
|f_\epsi^{[2,2]}(s,k)|_2 &\leq C |k|_1 |f_\epsi^{[1,2]}(s,k)|_2 
\\
&\leq 
C |k|_1 | \pt \T\big(\P_\epsi \uhat_1,\P_\epsi \uhat_1,\P_\epsi \uhat_1\big)(s,k) |_2
\\
&= C \sum_{\mu=1}^d |D_\mu \pt \T\big(\P_\epsi \uhat_1,\P_\epsi \uhat_1,\P_\epsi \uhat_1\big)(s,k)|_2, 
\end{align*}
such that \eqref{Step.6.bound.Tdot.Dmu} yields uniform boundedness of $ \|f_\epsi^{[2,2]}(s)\|_{L^1} $.
Uniform boundedness of $ \|f_\epsi^{[2,3]}(s)\|_{L^1} $ follows from \eqref{Step.6.bound.Tdotdot}.
We have thus shown that all terms in the big bracket $[ \; \ldots \; ]$ in 
\eqref{Step.6.big.bracket} are uniformly bounded, which completes the proof of 
\eqref{Main.difficulty.1}.
\bigskip

According to step~4, the inequalities \eqref{Main.difficulty.2} and \eqref{Main.difficulty.1}
imply the bound \eqref{Theorem.error.bound.SVEA.main.difficulty}, which is equivalent to 
\eqref{Theorem.error.bound.SVEA.crucial}.
We have shown in step~2 that this concludes the proof of \eqref{Theorem.error.bound.SVEA.assertion01}
and hence of Theorem~\ref{Theorem.error.bound.SVEA}.
\qed

The proof shows that in general the error of the SVEA cannot be expected to be smaller than $\Ord{\epsi^2}$.
We have seen in step 2 that the accuracy is determined by the right-hand side of \eqref{Theorem.error.bound.SVEA.crucial},
and in order to improve this inequality, we have to replace \eqref{Main.difficulty.2} and 
\eqref{Main.difficulty.1} by something better\footnote{In addition, a number of terms which were estimated in a straightforward way in our proof would require a more sophisticated analysis.}.
But this is impossible, which can be seen in the proof of \eqref{Main.difficulty.1} in step~6.
Since $ \sup_{t\in[0,\tend/\epsi]} \| \P_\epsi \uhat_1(t) \|_{L^1} = \Ord{1} $
it follows that
\begin{align}
\label{Theorem.error.bound.SVEA.bottleneck.term}
\| f_\epsi(t) \|_{L^1} =
\| \T\big(\P_\epsi \uhat_1,\P_\epsi \uhat_1,\P_\epsi \uhat_1\big)(t) \|_{L^1} = \Ord{1},
\end{align}
and as a consequence, the right-hand side of \eqref{Step.6.big.bracket} cannot be smaller than $ \Ord{\epsi} $. 
We would like to point out that \eqref{Theorem.error.bound.SVEA.bottleneck.term} is not the only bottleneck in the proof, and that there are many terms for which a better bound is not feasible.
The only way to achieve a higher accuracy is thus to change the approximation, i.e. to use 
\eqref{Ansatz} and \eqref{PDE.mfe} with $m>1$. 
This is the topic of Section~\ref{Sec.error.bound.approximation.m3}.

\begin{Remark}
\label{Remark:eigenspace.one-dimensional}
We have assumed throughout that the kernel of $\L(\w,\kappa)$ is one-dimensional; see~Assumption~\ref{Ass:polarization}\ref{Ass:polarization.item.i}. 
In case of the Maxwell--Lorentz system, however, most of the eigenvalues of $\L(0,\kappa) = A(\kappa)-\ii E$ occur with multiplicity 2, as pointed out in \cite[Example~3.2.4]{baumstark:22},
and if the eigenvalue $\w$ chosen in \eqref{Def.w} has multiplicity 2, then 
$\L(\w,\kappa)=-\w I+\L(0,\kappa)$ has a two-dimensional kernel.
For this reason, we would like to emphasize that the only purpose of 
Assumption~\ref{Ass:polarization}\ref{Ass:polarization.item.i} is to keep the notation simple,
and that all results and proofs in this work could be adapted to cases where the dimension of the kernel is two or larger, as has been done in \cite{baumstark:22}.
If the kernel of $\L(\w,\kappa)$ has dimension 2, then 
$ \ell=2, \ldots, n$ in Assumption~\ref{Ass:L.properties}\ref{Ass:L.properties.item.iii}
has to be replaced by $ \ell=3, \ldots, n$, and the definitions of the projectors
\eqref{Def.P2} and \eqref{Def.P1} have to be modified in an obvious way.
Likewise,
$ \sum_{\ell=2}^n \ldots $ has to be replaced by $ \sum_{\ell=3}^n \ldots $
in the proofs of Proposition~\ref{Proposition02.SVEA} and Theorem~\ref{Theorem.error.bound.SVEA}.
\end{Remark}

% --------------------------------------------------------------------------------
\subsection{Numerical experiment} \label{Subsec:NumEx.01}
% --------------------------------------------------------------------------------

We illustrate Theorem~\ref{Theorem.error.bound.SVEA} by a numerical example. As a model problem, 
we use a Klein--Gordon system in one space dimension; cf.~Example 2 in \cite{colin-lannes:09} and Example 1.5 in \cite{lannes:11}. This system is a special case of \eqref{PDE.uu.a} with
\begin{align*}
d=1, \qquad n=2, \qquad 
A(\px) = \begin{pmatrix} 0 & 1 \\ 1 & 0 \end{pmatrix} \partial_x, \qquad
A(\kappa) = \begin{pmatrix} 0 & \kappa \\ \kappa & 0 \end{pmatrix}, \qquad
E=\begin{pmatrix} 0 & -\gamma \\ \gamma & 0 \end{pmatrix}.
\end{align*}
We set
\begin{align*}
\tend=1, \qquad \kappa=1.2, \qquad \gamma=0.7, \qquad
T(f_1,f_2,f_3) = (f_1 \cdot f_2)Ef_3.
\end{align*}
The eigenvalues of $A(\kappa)-\ii E \in \IC^{2\times 2}$ are 
$ \pm\sqrt{\kappa^2 + \gamma^2} \approx \pm 1.3892 $, and we select $\w=\w(\kappa)$ to be the one with the positive sign.
For the initial data in \eqref{PDE.uu.b} we choose $ p(x) = \ee^{-(x-0.5)^2}\nu $
with $ \nu \in \ker(\L(\w,\kappa))$, 
such that the polarization condition (Assumption~\ref{Ass:polarization}\ref{Ass:polarization.item.ii}) holds with 
$p=p_0$ and $p_1=0$.
The initial data and the values for $\kappa$ and $\gamma$ were chosen more or less arbitrarily. The numerical results reported below remained qualitatively the same for other parametrizations we have tested.

Since numerical approximations of \eqref{PDE.uu} and \eqref{SVEA} can only be computed on a bounded domain, 
we switch to co-moving coordinates 
\begin{align*}
\xi = x - c_g t, \qquad v(t,\xi) = \uu(t,x), \qquad v_1(t,\xi) = u_1(t,x)
\end{align*}
with group velocity $ c_g = \nabla \w(\kappa) = \kappa/\w(\kappa) $.
For $d=1$ this turns \eqref{PDE.uu} into
\begin{subequations}
\label{PDE.vv}
\begin{align}
\label{PDE.vv.a}
 \pt v +A(\partial_\xi)v - c_g \partial_\xi v + \frac{1}{\epsi} Ev &= \epsi T(v,v,v),
 & & t\in(0,\tend/\epsi], \; \xi\in\IR,
 \\
\label{PDE.vv.b}
 v(0,\xi) &= p(\xi)\ee^{\ii (\kappa \xi)/\epsi} + c.c.,
\end{align}
\end{subequations}
and \eqref{Ansatz.SVEA}--\eqref{SVEA} into
\begin{align}
\label{SVEA.v.a}
 v(t,\xi) \approx \vtilde^{(1)}(t,\xi) &= \ee^{\ii (\kappa \xi + (\kappa c_g - \w) t)/\epsi} v_1(t,\xi) + c.c.,
 \\
 \notag
 \pt v_1 + \frac{\ii}{\epsi}\L(\w,\kappa)v_1 + A(\partial_\xi)v_1 - c_g \partial_\xi v_1
&= \epsi \sum_{j_1+j_2+j_3=1} T(v_{j_1},v_{j_2},v_{j_3}),
 \\
 \notag
v_1(0,\cdot) &= p.
\end{align}
Then, we replace $ \xi \in \IR$ by $\xi \in [-64,64]$ with periodic boundary conditions
and approximate $v_1$ with a Strang splitting method with very small step-size ($\tend/10^5$) 
and mesh-width ($128/2^{14}=2^{-7}$).
Inserting this numerical approximation of $v_1$ into \eqref{SVEA.v.a} yields a numerical approximation to 
$ \vtilde^{(1)} $, which is then compared with a numerical approximation to the solution of 
\eqref{PDE.vv}. 
As we have explained in the introduction, such an approximation can unfortunately not be obtained by applying a standard method to \eqref{PDE.vv} in a straightforward way, because the highly oscillatory solution behavior imposes a very fine discretization in time and space, which causes huge computational costs even in one space dimension. As a remedy, we have used
\eqref{Ansatz} and \eqref{PDE.mfe} with $m=5$ to compute a reference solution in co-moving coordinates.

\begin{figure}[htb] 
\begin{center}
	% This file was created by matlab2tikz.
%
%The latest updates can be retrieved from
%  http://www.mathworks.com/matlabcentral/fileexchange/22022-matlab2tikz-matlab2tikz
%where you can also make suggestions and rate matlab2tikz.
%
\begin{tikzpicture}

\begin{axis}[%
width=10cm,
height=5cm,
at={(0cm,0cm)},
scale only axis,
xmode=log,
xmin=0.01,
xmax=0.1,
xminorticks=true,
xlabel style={font=\color{white!15!black}},
xlabel={$\varepsilon$},
ymode=log,
ymin=0.0001,
ymax=1.5*0.01,
yminorticks=true,
ylabel style={font=\color{white!15!black}},
ylabel={},
axis background/.style={fill=white},
legend style={legend cell align=left, align=left, draw=white!15!black},
legend pos={south east}
]
\addplot [color=blue, line width=1pt, mark=o, mark options={solid, blue}]
  table[row sep=crcr]{%
0.01	0.000126534870897532\\
0.0158489319246111	0.000317390163294018\\
0.0251188643150958	0.000810128901363449\\
0.0398107170553497	0.00200414884252531\\
0.0630957344480193	0.00493080630866571\\
0.1	0.0123589134146372\\
};
\addlegendentry{error of SVEA}

\addplot [color=black, dashed, line width=1pt]
  table[row sep=crcr]{%
0.01	0.0001\\
0.0158489319246111	0.000251188643150958\\
0.0251188643150958	0.000630957344480193\\
0.0398107170553497	0.00158489319246111\\
0.0630957344480193	0.00398107170553497\\
0.1	0.01\\
};
\addlegendentry{$\varepsilon^2$}
\end{axis}

%\begin{axis}[%
%width=5.833in,
%height=3.375in,
%at={(0in,0in)},
%scale only axis,
%xmin=0,
%xmax=1,
%ymin=0,
%ymax=1,
%axis line style={draw=none},
%ticks=none,
%axis x line*=bottom,
%axis y line*=left
%]
%\end{axis}
\end{tikzpicture}%
\begin{minipage}{12cm}
\caption{Accuracy of the SVEA for different values of $\epsi$. See text for details.}
  \label{Fig.error_SVEA}  
\end{minipage}
\end{center}
\end{figure}
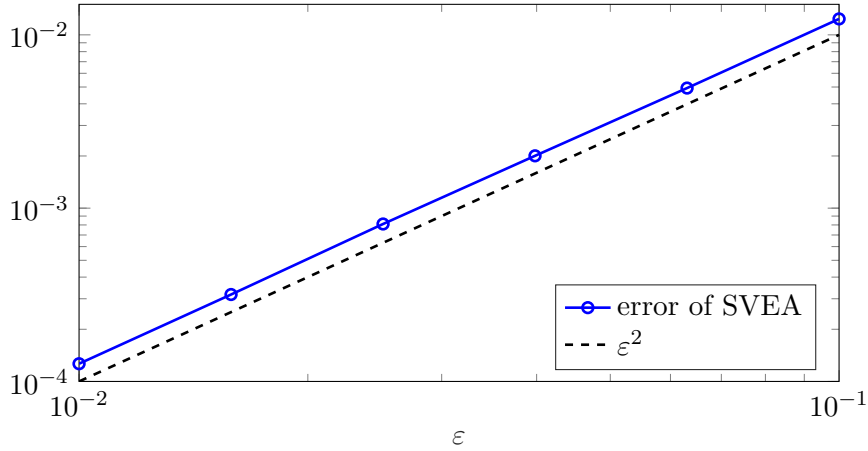

Figure \ref{Fig.error_SVEA} shows the numerical counterpart of 
\begin{align*}
\sup_{t\in[0,\tend/\epsi]} \| v(t,\cdot) - \vtilde^{(1)}(t,\cdot) \|_{L^\infty}
\end{align*}
for different values of $\epsi$ (blue line) in logarithmic axes. 
Comparing with $ \epsi \mapsto \epsi^2 $ (black dashed line) shows that the error is proportional to $\epsi^2$, as predicted by Theorem \ref{Theorem.error.bound.SVEA}. 

% --------------------------------------------------------------------------------
\section{Convergence analysis for \boldmath $m=3$}
\label{Sec.error.bound.approximation.m3}
% --------------------------------------------------------------------------------

In this section, we analyze the approximation \eqref{Ansatz} with 
\begin{align*}
m=3 \qquad \text{and} \qquad \J^{(3)} = \{\pm 1, \pm 3\}.
\end{align*}
As mentioned in the introduction, we have shown in \cite{baumstark-jahnke-2023} that $\uutilde^{(3)}$ approximates the exact solution $\uu$ of \eqref{PDE.uu} up to an error of $\Ord{\epsi^2}$; cf.~\eqref{error.bound.BJ.3}. In this section, we prove that actually the error is only $ \Ord{\epsi^3} $ if \eqref{Def:p} holds with $\sigma=3$.

By definition the approximation 
\begin{align*}
\uutilde^{(3)}(t,x) = 
 \sum_{j\in \J^{(3)}}
\ee^{\ii j(\kappa \cdot x - \w t)/\epsi} u_j(t,x)
=
\Big(
\ee^{\ii (\kappa \cdot x - \w t)/\epsi} u_1(t,x) + \ee^{3\ii (\kappa \cdot x - \w t)/\epsi} u_3(t,x)
\Big) + c.c.
\end{align*}
is based on two functions $u_1, u_3$ and their complex conjugates $ u_{-j} = \overline{u_j} $.
These functions $u_1, u_3$ are determined by the coupled system \eqref{PDE.mfe}, and thus
$u_1$ (which now depends on $u_3$) is \emph{not} the same as $u_1$ in the case $m=1$.
As a consequence, we cannot readily use the auxiliary results concerning $u_1$ which we have shown in Lemmas~\ref{Lemma.bound.uhat1dot}, \ref{Lemma.bound.Dmu.uhat1dot},
and \ref{Lemma.bound.uhat1dotdot}.
However, these results can be extended to the case $m=3$ with little effort. 
We summarize this in the following lemma.

\begin{Lemma}\label{Lemma.bounds.u1.m3}
Let $m = 3$, let $\sigma=1$ in \eqref{Def:p}, and let $\{u_1,u_3\}$ be the classical solution of \eqref{PDE.mfe}.  
\begin{enumerate}[label=(\roman*)]
\item 
Under Assumptions \ref{Ass:L.properties} and 
\ref{Ass:polarization}\ref{Ass:polarization.item.i}, 
there is a constant $C$ independent of $\epsi\in(0,1]$ such that
\begin{align} 
\label{Lemma.bound.uhat1dot.m3}
\sup_{t\in[0,\tend/\epsi]} 
\| \pt \P_\epsi \uhat_1(t) \|_{L^1} \leq C,
\\
\label{Lemma.bound.Dmu.uhat1dot.m3}
\sup_{t\in[0,\tend/\epsi]} 
\| \pt D_\mu \P_\epsi \uhat_1(t) \|_{L^1} \leq C.
\end{align}

\item
If in addition $\sigma=2$ in \eqref{Def:p}, then there is a constant $C$ such that
\begin{align} 
\label{Lemma.bound.uhat1dotdot.m3}
\sup_{t\in[0,\tend/\epsi]} 
\| \pt^2 \P_\epsi \uhat_1(t) \|_{L^1} \leq C.
\end{align}

\item
If in addition $\sigma=3$ in \eqref{Def:p}, then there is a constant $C$ such that
\begin{align} 
\label{Lemma.bound.Dmu.uhat1dotdot.m3}
\sup_{t\in[0,\tend/\epsi]} 
\| D_\mu \pt^2 \P_\epsi \uhat_1(t) \|_{L^1} \leq C.
\end{align}
\end{enumerate}
\end{Lemma}

\proof
The bound \eqref{Lemma.bound.uhat1dot.m3} was shown in \cite[Lemma 3.5]{baumstark-jahnke-2023}.
To show \eqref{Lemma.bound.Dmu.uhat1dot.m3} and \eqref{Lemma.bound.uhat1dotdot.m3}, the proofs
of Lemma~\ref{Lemma.bound.Dmu.uhat1dot} and \ref{Lemma.bound.uhat1dotdot}
carry over almost verbatim. The only difference is that for $m=3$ the sum
$ \sum_{\#J=1} \T(\uhat_{j_1},\uhat_{j_2},\uhat_{j_3}) $ contains more terms than for $m=1$,
for example $ \T(\uhat_3,\uhat_{-1},\uhat_{-1}) $, and thus the constants in the very last 
inequality of each proof change a bit.
The proof of \eqref{Lemma.bound.Dmu.uhat1dotdot.m3} is more complicated than the proof of
\eqref{Lemma.bound.uhat1dotdot.m3}, because new terms arise due to the presence of $D_\mu$, but these terms do not cause any essential new difficulty.
\qed

% --------------------------------------------------------------------------------
\subsection{Bounds on the coefficient functions}
% --------------------------------------------------------------------------------

As a first step, we prove that for $m=3$ it is still true that 
$ \| \P_\epsi^\perp \uhat_1(t) \|_{L^1} = \Ord{\epsi} $, and that in addition
$ \| \uhat_3(t) \|_{L^1} = \Ord{\epsi^2} $;
cf.~Corollary~\ref{Corollary01.m3} below.
For this purpose, we define the scaled norm $ \il \cdot \il_\epsi $ 
of a pair $ Y^{(3)} = \{y_1, y_3\} $ of functions $ y_j \in L^1$ by
\begin{align}
 \label{Def.scaled.norm.m3}
  \il Y^{(3)} \il_\epsi = 2\| Py_1 \|_{L^1} + \frac{2}{\epsi} \| \Pperp y_1 \|_{L^1}
  + \frac{2}{\epsi^2} \| y_3 \|_{L^1}.
\end{align}
In \cite[Equation (3.3)]{baumstark-jahnke-2023}
we have used a similar definition, but with factor $2/\epsi$ instead of $2/\epsi^2$ in the last term.
This difference is important.
The motivation for multiplying every term on the right-hand side with 2 is that then 
Equation~\eqref{Prop01.m3.proof.02} below holds true.

\begin{Proposition}\label{Proposition01.m3}
Suppose that the initial data in \eqref{PDE.mfe.b} have the form \eqref{Def:p} with $\sigma=2$.
Let $U^{(3)}=\{u_1,u_3\}$ be the classical solution of \eqref{PDE.mfe} with $m=3$ 
and let $\Uhat^{(3)}=\{\uhat_1,\uhat_3\}$.
Let $z_1$ and $z_3$ be the transformed functions defined in \eqref{Def.z},
and set $Z^{(3)}=\{z_1,z_3\}$.
For every sufficiently large $r>0$ there is a $\tstar\in(0,\tend]$ such that 
under the Assumptions~\ref{Ass:polarization} and \ref{Ass:L.properties}
\begin{align} 
\label{Proposition01.m3.inequality}
\sup_{t\in[0,\tstar/\epsi]} \il Z^{(3)}(t) \il_\epsi \leq r \qquad
\text{for all } \epsi\in(0,1].
\end{align}
The constant $\tstar$ depends on $\tend$, $r$, $C_{u,2}$, $C_\T$, on the inverse of the nonzero eigenvalues of $\Lambda_1(0)$, and on the Lipschitz constant in \eqref{Ass:lambda.Lipschitz}, but not on $\epsi$.
\end{Proposition}

\begin{Remark}
The proof yields an explicit formula for $\tstar$; cf.~\eqref{Prop01.m3.tstar.formula}.
Numerical computations indicate that this formula is way too pessimistic in most cases, 
but for our goals it is sufficient 
that for every $r$ there is a $\tstar$ such that \eqref{Proposition01.m3.inequality} holds, 
and that $\tstar$ does not depend on $\epsi$.
The number $\tstar$ obtained from \eqref{Prop01.m3.tstar.formula} is positive only if $ r > C_\bullet$, where $C_\bullet$ is a constant which appears in the proof. This is what we mean by ``sufficiently large $r$''.
\end{Remark}

\bigskip

\noindent
Before we prove Proposition~\ref{Proposition01.m3}, we note that the following corollary is an immediate consequence of \eqref{Norm.uhat.z.Pperp}, \eqref{Norm.uhat.z}, \eqref{Def.scaled.norm.m3}, and \eqref{Proposition01.m3.inequality}.

\begin{Corollary}\label{Corollary01.m3}
Under the assumptions of Proposition~\ref{Proposition01.m3} the bounds
\begin{alignat*}{2}
\sup_{t\in[0,\tstar/\epsi]}\| \P_\epsi^\perp \uhat_1(t) \|_{L^1} 
&=
\sup_{t\in[0,\tstar/\epsi]}\| \Pperp z_1(t) \|_{L^1} 
&&\leq C \epsi,
\\
\sup_{t\in[0,\tstar/\epsi]}\| \uhat_3(t) \|_{L^1} 
&=
\sup_{t\in[0,\tstar/\epsi]}\| z_3(t) \|_{L^1} 
&&\leq C \epsi^2
\end{alignat*}
hold with a constant independent of $\epsi\in(0,1]$. 
\end{Corollary}

Corollary~\ref{Corollary01.m3} reveals that Proposition~\ref{Proposition01.m3}
can be understood as an extension of Proposition~\ref{Proposition01.SVEA} from $m=1$ to the case $m=3$.
However, a substantial difference between the two cases is the fact that
the proof of Proposition~\ref{Proposition01.SVEA} (see \cite[Lemma 3]{colin-lannes:09}) is based on Gronwall's lemma,
whereas the proof of Proposition~\ref{Proposition01.m3} requires other techniques. The reason is, roughly speaking, that for $m=3$ there are two functions, $ \P_\epsi^\perp \uhat_1(t) $ and $ \uhat_3(t) $, 
which we have to estimate simultaneously.
This is also the reason why Proposition~\ref{Proposition01.m3} refers to a possibly smaller interval 
$[0,\tstar/\epsi]$ instead of $[0,\tend/\epsi]$.
\bigskip

\proofof{Proposition~\ref{Proposition01.m3}}
We integrate \eqref{zdot.in.terms.of.uhat} for $m=3$ from $0$ to $t\in[0,\tend/\epsi]$.
This yields
\begin{align}
\notag
\il Z^{(3)}(t) \il_\epsi 
 \leq
 \il Z^{(3)}(0)\il_\epsi
 &+ 2 \sum_{\# J=1} \left(
 \epsi \Big\| \intl_0^t P \FF(s,\Uhat^{(3)},J) \; \dd s \Big\|_{L^1}
 +
 \Big\| \intl_0^t \Pperp \FF(s,\Uhat^{(3)},J) \; \dd s \Big\|_{L^1}
 \right)
 \\
 \label{Prop01.m3.proof.01}
 & + \frac{2}{\epsi} \sum_{\# J=3} \Big\| \intl_0^t \FF(s,\Uhat^{(3)},J) \; \dd s \Big\|_{L^1}
\end{align}
with $\FF$ defined in \eqref{Def.F.uhat}. 
Since $z_3(0)=0$ by \eqref{PDE.z.initial.data}, it follows from \eqref{Bound.projection.initial.value} that
\begin{align*}
 \il Z^{(3)}(0)\il_\epsi
 =
2\| Pz_1(0) \|_{L^1} + \frac{2}{\epsi} \| \Pperp z_1(0) \|_{L^1}
\leq C (\| p_0 \|_{W^1} + \| p_1 \|_{W^1}).
\end{align*}
Now, we define
\begin{align}
\label{Def.a1.a3}
a_{\pm1}(t) &= \| Pz_1(t) \|_{L^1} + \frac{1}{\epsi} \| \Pperp z_1(t) \|_{L^1} 
& \text{and} &&
a_{\pm3}(t) &= \frac{1}{\epsi^2} \| z_3(t) \|_{L^1} 
\end{align}
and note that
\begin{align}
 \label{Prop01.m3.proof.02}
\sum_{j\in\J^{(3)}} a_j(s) = 2 a_1(s) + 2 a_3(s) = \il Z^{(3)}(s) \il_\epsi
\end{align}
by \eqref{Def.scaled.norm.m3}.
Our goal is to prove that there are constants $C_\star$ and $\Chat$ such that
for all $t\in[0,\tend/\epsi]$ the inequality 
\begin{align}
\label{Target.term.1}
\epsi \Big\| \intl_0^t P \FF(s,\Uhat^{(3)},J) \; \dd s \Big\|_{L^1}
+
\Big\| \intl_0^t \Pperp \FF(s,\Uhat^{(3)},J) \; \dd s \Big\|_{L^1}
\leq
C_\star +  \Chat \epsi \intl_0^t \prod_{i=1}^3 a_{j_i}(s) \; \dd s
\end{align}
holds for every $J=(j_1,j_2,j_3)\in (\J^{(3)})^3$ with $\#J=1$, and that
\begin{align}
\label{Target.term.3}
\frac{1}{\epsi}
\Big\| \intl_0^t \FF(s,\Uhat^{(3)},J) \; ds \Big\|_{L^1}
\leq
C_\star +  \Chat \epsi \intl_0^t \prod_{i=1}^3 a_{j_i}(s) \; \dd s
\end{align}
holds for every $J=(j_1,j_2,j_3)\in (\J^{(3)})^3$ with $\#J=3$.
Substituting \eqref{Target.term.1} and \eqref{Target.term.3} into \eqref{Prop01.m3.proof.01} yields
\begin{align}
\notag
 \il Z^{(3)}(t) \il_\epsi 
   &\leq 
  C_\bullet
  + 2\Chat \epsi \sum_{j\in\{1,3\}} \sum_{\# J=j} 
 \intl_0^t \prod_{i=1}^3 a_{j_i}(s) \; \dd s
 \\
 \notag
 &\leq 
  C_\bullet
  + \Chat \epsi \sum_{j\in\J^{(3)}} \sum_{\# J=j} 
 \intl_0^t \prod_{i=1}^3 a_{j_i}(s) \; \dd s
 \\
 \notag
  &=
  C_\bullet + \Chat \epsi 
 \intl_0^t \bigg(\sum_{j\in\J^{(3)}} a_j(s)\bigg)^3 \; \dd s
 \\
 \label{Prop01.m3.proof.03}
 &=
C_\bullet + \Chat \epsi \intl_0^t \il Z^{(3)}(s) \il_\epsi^3 \; \dd s
\end{align}
by \eqref{Prop01.m3.proof.02}. The constant $C_\bullet$ depends on $\|p_0\|_{W^1}$, $\|p_1\|_{W^1}$, $C_\star$ and the (finite) number of multi-indices $J$ with $\#J=1$ and $\#J=3$, respectively.
Now let $ \tstar \in (0,\tend] $ be a number to be determined below.
Then, \eqref{Prop01.m3.proof.03} implies that 
\begin{align} 
\label{Prop01.m3.proof.04}
 \il Z^{(3)}(t) \il_\epsi 
\leq
C_\bullet + \Chat \tstar \sup_{s\in[0,\tstar/\epsi]} \il Z^{(3)}(s) \il_\epsi^3
\qquad \text{for all } t \in[0,\tstar/\epsi].
\end{align}
If we choose $\tstar$ in such a way that the right-hand side of this inequality is not larger than $r$
for some $r>C_\bullet$,
then we can infer from \eqref{Prop01.m3.proof.04} that
$  \il Z^{(3)}(t) \il_\epsi \leq r $ for all $t\in[0,\tstar/\epsi]$. 
Hence, the desired inequality \eqref{Proposition01.m3.inequality} holds with
\begin{align}
\label{Prop01.m3.tstar.formula}
  \tstar = \frac{r-C_\bullet}{\Chat r^3}.
\end{align}

To prove the first inequality \eqref{Target.term.1} we can adapt the arguments 
from \cite[Section 3.2.2]{baumstark-jahnke-2023}, because the fact that $ a_3(t) $ was defined with a different prefactor in \cite[Eq. (3.7)]{baumstark-jahnke-2023} does not matter for this part.
To complete the proof of Proposition~\ref{Proposition01.m3}, we have to show \eqref{Target.term.3}.
In \cite[Section 3.2.1]{baumstark-jahnke-2023} we have proven such a bound, but without the factor
$1/\epsi$ on the left-hand side.
Let $\#J=3$ and recall that 
\begin{align*}
\FF(s,\Uhat^{(3)},J) 
&=
\TT_{j,\epsi}(s) \T\big(\uhat_{j_1},\uhat_{j_2},\uhat_{j_3} \big)(s),
\qquad j=\#J=3
  \end{align*}
according to \eqref{Def.F.uhat}. 
We first consider the (easy) case where $|J|_1 > \#J=3$ and thus $|J|_1 \geq 5$ because $|J|_1$ is an odd integer.
In this case, \eqref{Bound.T.Fourier} implies
 \begin{align*}
\frac{1}{\epsi}
\Big\| \intl_0^t \FF(s,\Uhat^{(3)},J) \; \dd s \Big\|_{L^1}
&=
\frac{1}{\epsi}
   \intl_0^t \Big\|\T\big(\uhat_{j_1},\uhat_{j_2},\uhat_{j_3} \big)(s) 
   \Big\|_{L^1} \; \dd s 
   \\
   &\leq
\frac{C_\T}{\epsi}
   \intl_0^t \prod_{i=1}^3 \|\uhat_{j_i}(s)\|_{L^1} \; \dd s 
   \intertext{ }
   &=
\frac{C_\T}{\epsi}
   \intl_0^t \prod_{i=1}^3 \|z_{j_i}(s)\|_{L^1}   \; \dd s   
\\
 &=
 C_\T \epsi^{|J|_1-4} 
  \intl_0^t 
  \prod_{i=1}^3 \Big( \epsi^{1-|j_i|} \| z_{j_i}(s) \|_{L^1} \Big)
  \; \dd s  
  \\
  &\leq
  C_\T \epsi \intl_0^t \prod_{i=1}^3 a_{j_i}(s) \; \dd s,  
   \end{align*} 
which is an estimate of the form \eqref{Target.term.3} with $C_\star=0$.  
   In the last step, we have used that $|J|_1-4 \geq 1$ and 
  $\epsi^{1-|j_i|} \| z_{j_i}(s) \|_{L^1} \leq a_{j_i}(s)$
by definition \eqref{Def.a1.a3}.

Now let $|J|_1 = \#J = 3$, which is only true for $ J=(1,1,1) $.
Since 
\begin{align*}
\FF\big(s,\Uhat^{(3)},(1,1,1)\big) = \TT_{3,\epsi}(s) \T\big(\uhat_1,\uhat_1,\uhat_1 \big)(s),
\end{align*}
we have to show that
\begin{align}
\label{Target.term.3.111}
\frac{1}{\epsi}
  \Big\| \intl_0^t \TT_{3,\epsi}(s) \T\big(\uhat_1,\uhat_1,\uhat_1 \big)(s) \; \dd s \Big\|_{L^1}
&\leq
C_\star +  \Chat \epsi \intl_0^t a_1^3(s) \; \dd s.
\end{align}
At this point, it seems that the inequality \eqref{Main.difficulty.0}, which we have shown in 
steps 4--6 of the proof of Theorem~\ref{Theorem.error.bound.SVEA}, readily implies \eqref{Target.term.3.111}
with $\Chat=0$.
This is not quite true, because \eqref{Main.difficulty.0} refers to the case $m=1$, not $m=3$, 
and we have pointed out at the beginning of this section that
$\uhat_1$ is not the same function in these two cases.
But the parts (i) and (ii) of Lemma~\ref{Lemma.bounds.u1.m3} ensure that for $m=3$ 
the function $\uhat_1$ has still all the properties which were used
to prove \eqref{Main.difficulty.0}, and this allows us to use that proof verbatim.
\qed

\noindent
Before we proceed, we have to extend Corollary~\ref{Corollary01.m3} to a stronger norm
as in Section~\ref{Sec.error.SVEA.preparations}.
The following result is the counterpart of Proposition~\ref{Proposition02.SVEA} in the case $m=3$.

\begin{Proposition}\label{Proposition02.m3}
Suppose that the assumptions of Proposition~\ref{Proposition01.m3} hold, and that in addition 
\eqref{Def:p} is true with $\sigma=3$.
Then, the bounds
\begin{alignat*}{2}
\sup_{t\in[0,\tstar/\epsi]}\| D_\mu \P_\epsi^\perp \uhat_1(t) \|_{L^1} 
&=
\sup_{t\in[0,\tstar/\epsi]}\| D_\mu \Pperp z_1(t) \|_{L^1} 
&&\leq C \epsi,
\\
\sup_{t\in[0,\tstar/\epsi]}\| D_\mu \uhat_3(t) \|_{L^1} 
&=
\sup_{t\in[0,\tstar/\epsi]}\| D_\mu z_3(t) \|_{L^1} 
&&\leq C \epsi^2
\end{alignat*}
hold with a constant independent of $\epsi\in(0,1]$. 
\end{Proposition}

\proof
Using the higher regularity and in particular \eqref{Lemma.bound.Dmu.uhat1dotdot.m3},
the bound
\begin{align*} 
\sup_{t\in[0,\tstar/\epsi]} \il D_\mu Z^{(3)}(t) \il_\epsi \leq C \qquad
\text{for all } \epsi\in(0,1], \mu\in\{1, \ldots, d\}.
\end{align*}
can be shown with standard techniques.
Then, the assertion follows from the definition \eqref{Def.scaled.norm.m3}.
\qed

% --------------------------------------------------
\subsection{Improved error bound for $\mathbf{m=3}$}
% --------------------------------------------------

For the error analysis of $\uutilde^{(3)}$ we need a second non-resonance condition similar to Assumption~\ref{Ass:Delta.lambda.j3}.

\begin{Assumption}[Non-resonance condition]
\label{Ass:Delta.lambda.j5}
The matrix $\L_{5}(0) = \L(5\w,5\kappa)$ is regular and has no common eigenvalues with
$\L_{3}(0)=\L(3\w,3\kappa)$, 
i.e. $ \lambda_{5i}(0) \neq \lambda_{3\ell}(0) $ for all $i, \ell =1, \ldots, n$.
\end{Assumption}

\noindent
We are now in a position to formulate and prove our second main result.

\begin{Theorem}[Error bound for $m=3$]\label{Theorem.error.bound.m3}
Let $ p $ have the form \eqref{Def:p} with $\sigma=3$ and let $\uu$ be the solution of \eqref{PDE.uu}.
Let $\uutilde^{(3)}$ be the approximation defined in \eqref{Ansatz} with $m=3$.
Under Assumptions~\ref{Ass:polarization}, \ref{Ass:L.properties}, \ref{Ass:Delta.lambda.j3}, 
and \ref{Ass:Delta.lambda.j5} there is a constant such that 
\begin{align}
\label{Theorem.error.bound.m3assertion01}
 \sup_{t\in[0,\tstar/\epsi]} \| \uu(t) - \uutilde^{(3)}(t) \|_W
 &\leq C\epsi^3,
 \\
\label{Theorem.error.bound.m3assertion02}
 \sup_{t\in[0,\tstar/\epsi]} \| \uu(t) - \uutilde^{(3)}(t) \|_{L^\infty}
 &\leq C\epsi^3.
\end{align}
\end{Theorem}

\proof 
We use the proofs of Theorem 4.2 in \cite{baumstark-jahnke-2023} and
of Theorem~\ref{Theorem.error.bound.SVEA} in the present paper 
as a blueprint and focus on what has to be changed.
In \cite[proof of Theorem 4.2]{baumstark-jahnke-2023} we have shown that the Fourier transform 
$\widehat{\delta}$ 
of $ \delta=\uu-\uutilde^{(3)} $ is the solution of
\begin{align*} 
 % \label{PDE.delta.Fourier.m3}
\pt \widehat{\delta}(t,k) =  
-\big(\ii A(k)+\tfrac{1}{\epsi} E\big)
\widehat{\delta}(t,k) 
+ \epsi \G\big(\F\uu,\F\uutilde^{(3)} \big)(t,k)
+ \widehat{R}(t,k) 
\end{align*}
with 
\begin{align*}
\G\big(\F\uu,\F\uutilde^{(3)}\big)&=
\T(\F\uu,\F\uu,\F\uu)-\T\left(\F\uutilde^{(3)},\F\uutilde^{(3)},\F\uutilde^{(3)}\right),
\\
\widehat{R}(t,k) 
&= 
 \epsi 
 \sum_{|j|\in\{5,7,9\}}  \sum_{\#J=j} 
 \T(\uhat_{j_1},\uhat_{j_2},\uhat_{j_3})(t,k-\tfrac{j\kappa}{\epsi})
\ee^{-\ii j \w t/\epsi},  
\end{align*}
and with $\T$ defined by \eqref{Def.T.Fourier}. 
Our main task is to prove that
\begin{align}
 \label{Theorem.error.bound.m3.crucial}
 \sup_{t\in[0,\tstar/\epsi]}
 \Big\|
 \intl_0^t
\exp\big((s-t)\big(\ii A(\cdot)+\tfrac{1}{\epsi} E\big)\big)
\widehat{R}(s) \; \dd s
\Big\|_{L^1} \leq C \epsi^3
\end{align}
uniformly in $\epsi\in(0,1]$. If \eqref{Theorem.error.bound.m3.crucial} holds, then the estimate \eqref{Theorem.error.bound.m3assertion01} can be shown by applying Duhamel's formula as in the proof of Theorem~\ref{Theorem.error.bound.SVEA},
and \eqref{Theorem.error.bound.m3assertion02} follows from the embedding 
$W \hookrightarrow L^\infty$.

In \cite[proof of Theorem 4.2]{baumstark-jahnke-2023}, we have already derived the inequality
\begin{align*}
 \Big\|
 \intl_0^t
\exp\big((s-t)\big(\ii A(\cdot)+\tfrac{1}{\epsi} E\big)\big)
\widehat{R}(s) \; \dd s
\Big\|_{L^1} 
&\leq
 \epsi 
 \sum_{|j| \in\{5,7,9\}}  \sum_{\#J=j} 
\Big\| \intl_0^t \FF(s,\Uhat^{(3)},J)
\; \dd s \Big\|_{L^1}
\\
&=
\epsi \sum_{|j| \in\{5,7,9\}}  \sum_{\#J=j} 
\Big\| \intl_0^t \TT_{j,\epsi}(s) \T\big(\uhat_{j_1},\uhat_{j_2},\uhat_{j_3} \big)(s)
\; \dd s \Big\|_{L^1}. 
\end{align*}
In order to prove \eqref{Theorem.error.bound.m3.crucial}, we thus have to show that
\begin{align}
\label{Theorem.error.bound.m3.basic}
 \sum_{|j| \in\{5,7,9\}}  \sum_{\#J=j} 
\Big\| \intl_0^t \TT_{j,\epsi}(s) \T\big(\uhat_{j_1},\uhat_{j_2},\uhat_{j_3} \big)(s)
\; \dd s \Big\|_{L^1} \leq C \epsi^2
\end{align}
with a constant $C$ which does not depend on $\epsi$ nor on $ t\in[0,\tstar/\epsi]$.

As before, we consider several cases. First, suppose that $|j|\in\{7,9\} $.
If $J=(j_1,j_2,j_3)\in(\J^{(3)})^3$ with $\#J=j$, then at least two of the three entries 
must have a modulus of 3, 
such that with Corollary~\ref{Corollary01.m3} we even obtain the bound
\begin{align}
\nonumber
\Big\| \intl_0^t \TT_{j,\epsi}(s) \T\big(\uhat_{j_1},\uhat_{j_2},\uhat_{j_3} \big)(s)
\; \dd s \Big\|_{L^1} 
&\leq
\frac{\tstar}{\epsi} \sup_{s\in[0,\tstar/\epsi]}
\| \T\big(\uhat_{j_1},\uhat_{j_2},\uhat_{j_3} \big)(s) \|_{L^1}
\\
\label{Theorem.error.bound.m3.simple.argument}
&\leq 
\frac{C \tstar}{\epsi} \sup_{s\in[0,\tstar/\epsi]} \prod_{i=1}^3 \| \uhat_{j_i}(s) \|_{L^1}
\leq C \epsi^3.
\end{align}
If $|j|=5$ and $\#J=j$ but $|J|_1>j$ (e.g. if $j=5$ and $J=(3,-1,3)$), we can proceed in the same way. 
The difficult case is that $|j|=5=\#J=|J|_1$.
We consider only $j=5$ and $J=(3,1,1)$, because all other such combinations can be treated analogously. 
Now we cannot use \eqref{Theorem.error.bound.m3.simple.argument},
because Corollary~\ref{Corollary01.m3}  yields only 
$ \prod_{i=1}^3 \| \uhat_{j_i}(s) \|_{L^1} = \| \uhat_3(s) \|_{L^1} 
\| \uhat_1(s) \|_{L^1}^2 \leq C \epsi^2 $, 
which is not enough due to the factor $ \tstar/\epsi$ in 
\eqref{Theorem.error.bound.m3.simple.argument}.

Since $ \uhat_1(t) = \P_\epsi \uhat_1(t) + \P_\epsi^\perp \uhat_1(t) $ and 
since $ \sup_{t\in[0,\tstar/\epsi]}\| \P_\epsi^\perp \uhat_1(t) \|_{L^1} \leq C \epsi $
by Corollary~\ref{Corollary01.m3}, 
the problem boils down to showing the bound
\begin{align} \label{main.difficulty.jmax3}
\Big\| \intl_0^t \TT_{5,\epsi}(s) \T\big(\uhat_{3},\P_\epsi \uhat_{1},\P_\epsi \uhat_{1} \big)(s) \; \dd s \Big\|_{L^1} 
\leq C \epsi^{2}.
\end{align} 
To prove this, we use similar techniques as in step~5 of the proof of Theorem~\ref{Theorem.error.bound.SVEA}.
The strategy is again to identify the oscillatory  ``parts'' of the integrand.

We use the representation
\begin{align}
\label{Theorem.error.bound.m3.representation.uhat3}
\uhat_3(t,k) 
&= \TT_{3,\epsi}^*(t,k) z_3(t,k)
= 
\sum_{\ell=1}^n
\exp\big(-\tfrac{\ii t}{\epsi}\lambda_{3\ell}(\epsi k)\big)z_{3\ell}(t,k)\psi_{3\ell}(\epsi k),
\end{align}
which follows from \eqref{eigendecomposition}, \eqref{Def.z}, and \eqref{Def.TT}.
With \eqref{Def.TT}, \eqref{Shorthand.notation.K}, and \eqref{Def.T.Fourier}, this allows us to reformulate the integral in \eqref{main.difficulty.jmax3} as
\begin{align*}
& \intl_0^t \TT_{5,\epsi}(s) \T\big(\uhat_{3},\P_\epsi \uhat_{1},\P_\epsi \uhat_{1} \big)(s) \; \dd s 
\notag
\\
\notag
&=
\frac{1}{(2\pi)^d} 
\intl_0^t
\exp\big(\tfrac{\ii s}{\epsi}\Lambda_5(\epsi k)\big)  \Psi_5^*(\epsi k)
\intl_{\#K = k} T\left(\uhat_3(s,k^{(1)}),\P_\epsi\uhat_1(s,k^{(2)}),\P_\epsi\uhat_1(s,k^{(3)})\right)
\; \dd K
\; \dd s 
\\
&= 
\sum_{\ell=1}^n \intl_{\#K = k}
\intl_0^t
\exp\big(\tfrac{\ii s}{\epsi}\big[\Lambda_5(\epsi k)-\lambda_{3\ell}(\epsi k^{(1)})I\big]\big)
f_{\epsi,\ell}(s,K)
\; \dd s 
\; \dd K 
\end{align*}
with 
\begin{align*}
f_{\epsi,\ell}(s,K)&=
\frac{1}{(2\pi)^d} 
\Psi_5^*(\epsi k)
 T\left(z_{3\ell}(s,k^{(1)})\psi_{3\ell}(\epsi k^{(1)}),
 \P_\epsi\uhat_1(s,k^{(2)}),\P_\epsi\uhat_1(s,k^{(3)})\right),
 \quad \#K = k.
\end{align*}
The left-hand side of \eqref{main.difficulty.jmax3} can thus be bounded by
\begin{align}
\notag
& \Big\| \intl_0^t 
\TT_{5,\epsi}(s) \T\big(\uhat_3,\P_\epsi \uhat_1,\P_\epsi\uhat_{1}\big)(s)
\; \dd s \Big\|_{L^1}
\\
\label{Theorem.error.bound.m3.Eq.01}
&\leq
\sum_{\ell=1}^n 
\intl_{\IR^d} \intl_{\#K = k}
\Big|
\intl_0^t
\exp\big(\tfrac{\ii s}{\epsi}\big[\Lambda_5(\epsi k)-\lambda_{3\ell}(\epsi k^{(1)})I\big]\big)
f_{\epsi,\ell}(s,K)
\; \dd s 
\Big|_2 \; \dd K \; \dd k.
\end{align}
After setting
\begin{align}
\label{Theorem.error.bound.m3.Def.Delta.ell}
\Delta_\ell(\theta,\theta^{(1)}) &= \Lambda_5(\theta)-\lambda_{3\ell}(\theta^{(1)})I
\qquad \text{for } \theta, \theta^{(1)} \in \IR^d,
\\
\nonumber
g_{\epsi,\ell}(s,K) &= 
\exp\left(\frac{\ii s}{\epsi}\big[
\Delta_\ell(\epsi k,\epsi k^{(1)}) - \Delta_\ell(0,0)\big]
\right)
f_{\epsi,\ell}(s,K),
\end{align}
the inner integral reads
\begin{align*}
\Big| \intl_0^t
\exp\big(\tfrac{\ii s}{\epsi}\big[\Lambda_5(\epsi k)-\lambda_{3\ell}(\epsi k^{(1)})I\big]\big)
f_{\epsi,\ell}(s,K)
\; \dd s \Big|_2
=
\Big| \intl_0^t
\exp\left(\frac{\ii s}{\epsi}\Delta_\ell(0,0)\right)
g_{\epsi,\ell}(s,K)
\; \dd s \Big|_2.
\end{align*}
By Assumption~\ref{Ass:Delta.lambda.j5}
the diagonal matrix $ \Delta_\ell(0,0) = \Lambda_5(0)-\lambda_{3\ell}(0)I $ is regular such that we can integrate by parts and obtain
\begin{align}
\notag
\Big| \intl_0^t
\exp\left(\frac{\ii s}{\epsi}\Delta_\ell(0,0)\right)
g_{\epsi,\ell}(s,K)
\; \dd s \Big|_2
&\leq
C \epsi \Big(| g_{\epsi,\ell}(t,K) |_2 + | g_{\epsi,\ell}(0,K) |_2 \Big)
\\
\label{Theorem.error.bound.m3.Eq.02}
&\quad 
+ 
C \epsi 
\Big| \intl_0^t \exp\left(\frac{\ii s}{\epsi}\Delta_\ell(0,0)\right)
\pt g_{\epsi,\ell}(s,K) \; \dd s \Big|_2.
\end{align}
The term $ C \epsi (| g_{\epsi,\ell}(t,K) |_2 + | g_{\epsi,\ell}(0,K) |_2) $ on the right-hand side leads  
to a contribution of $ \Ord{\epsi^3} $ in \eqref{Theorem.error.bound.m3.Eq.01}, because
\begin{align*}
\sum_{\ell=1}^n | g_{\epsi,\ell}(t,K) |_2 
&= 
\sum_{\ell=1}^n | f_{\epsi,\ell}(s,K) |_2 
\\
&\leq
C 
\sum_{\ell=1}^n \Big| T\left(z_{3\ell}(s,k^{(1)})\psi_{3\ell}(\epsi k^{(1)}),
 \P_\epsi\uhat_1(s,k^{(2)}),\P_\epsi\uhat_1(s,k^{(3)})\right) \Big|_2 
\\
&\leq
C 
|z_{3}(s,k^{(1)})|_2 | \P_\epsi\uhat_1(s,k^{(2)}) |_2 
| \P_\epsi\uhat_1(s,k^{(3)}) |_2
\\
&\leq C\epsi^2
 \end{align*}
 by Corollary~\ref{Corollary01.m3}.
For the other term on the right-hand side of \eqref{Theorem.error.bound.m3.Eq.02}, the product rule gives  
\begin{subequations}
\label{Theorem.error.bound.m3.Eq.03}
\begin{align}
\nonumber
& C \epsi 
\Big| \intl_0^t 
\exp\left(\frac{\ii s}{\epsi}\Delta_\ell(0,0)\right)
\pt g_{\epsi,\ell}(s,K) \; \dd s \Big|_2
\\
\notag
&\leq
C \epsi 
\Big| \intl_0^t 
\exp\left(\frac{\ii s}{\epsi}\Delta_\ell(0,0)\right)
\frac{\ii}{\epsi}\big[\Delta_\ell(\epsi k,\epsi k^{(1)}) - \Delta_\ell(0,0)\big]
g_{\epsi,\ell}(s,K) 
\; \dd s \Big|_2
\\
\notag
&\quad + 
C \epsi 
\Big| 
\intl_0^t 
\exp\left(\frac{\ii s}{\epsi}\Delta_\ell(0,0)\right)
\exp\left(\frac{\ii s}{\epsi}\big[
\Delta_\ell(\epsi k,\epsi k^{(1)}) - \Delta_\ell(0,0)\big]
\right)
\pt f_{\epsi,\ell}(s,K)
\; \dd s \Big|_2
\\
&\leq
\label{Theorem.error.bound.m3.Eq.03.a}
C \epsi 
\Big| \intl_0^t 
\frac{\ii}{\epsi}\big[\Delta_\ell(\epsi k,\epsi k^{(1)}) - \Delta_\ell(0,0)\big]
\exp\left(\frac{\ii s}{\epsi}\Delta_\ell(\epsi k,\epsi k^{(1)})\right)
f_{\epsi,\ell}(s,K) 
\; \dd s \Big|_2
\\
\label{Theorem.error.bound.m3.Eq.03.b}
&\quad + 
C \epsi 
\Big| 
\intl_0^t 
\exp\left(\frac{\ii s}{\epsi}\Delta_\ell(\epsi k,\epsi k^{(1)})\right)
\pt f_{\epsi,\ell}(s,K)
\; \dd s \Big|_2.
\end{align}
\end{subequations}
The Lipschitz continuity \eqref{Ass:lambda.Lipschitz} of the eigenvalues yields
\begin{align*}
\big|\Delta_\ell(\epsi k,\epsi k^{(1)}) - \Delta_\ell(0,0)\big|_2
&\leq
\big| \Lambda_5(\epsi k)-\Lambda_5(0) \big|_2
+ \big| \lambda_{3\ell}(\epsi k^{(1)})-\lambda_{3\ell}(0) \big|
\\
&\leq
C \epsi (|k|_1 + |k^{(1)}|_1),
\end{align*}
and together with Corollary~\ref{Corollary01.m3} it can be shown that \eqref{Theorem.error.bound.m3.Eq.03.a}
causes a contribution of $\Ord{\epsi^2}$ in \eqref{Theorem.error.bound.m3.Eq.01}.

Unfortunately, the term \eqref{Theorem.error.bound.m3.Eq.03.b} requires a bit more efforts. By definition of
$ f_{\epsi,\ell} $, we formally have
\begin{align*}
\pt f_{\epsi,\ell}(s,K)
&=
\frac{1}{(2\pi)^d} 
\Psi_5^*(\epsi k)
\Big[
 T\left(\pt z_{3\ell}(s,k^{(1)})\psi_{3\ell}(\epsi k^{(1)}),
 \P_\epsi\uhat_1(s,k^{(2)}),\P_\epsi\uhat_1(s,k^{(3)})\right)
 \\
 & \hs{25}
 + T\left(z_{3\ell}(s,k^{(1)})\psi_{3\ell}(\epsi k^{(1)}),
 \pt \P_\epsi \uhat_1(s,k^{(2)}),\P_\epsi\uhat_1(s,k^{(3)})\right)
 \\
 & \hs{25}
 + T\left(z_{3\ell}(s,k^{(1)})\psi_{3\ell}(\epsi k^{(1)}),
 \P_\epsi\uhat_1(s,k^{(2)}),\pt \P_\epsi\uhat_1(s,k^{(3)})\right)
 \Big]
 \\
 &=
\frac{1}{(2\pi)^d} 
\Psi_5^*(\epsi k)
 T\left(\pt z_{3\ell}(s,k^{(1)})\psi_{3\ell}(\epsi k^{(1)}),
 \P_\epsi\uhat_1(s,k^{(2)}),\P_\epsi\uhat_1(s,k^{(3)})\right)
 + \Ord{\epsi^2}
\end{align*}
because $\pt \P_\epsi\uhat_1$ is uniformly bounded by \eqref{Lemma.bound.uhat1dot.m3}
and $ z_{3\ell}(s,k^{(1)}) = \Ord{\epsi^2} $ due to Corollary~\ref{Corollary01.m3}.
Proving the desired bound for the $\Ord{\epsi^2}$-part of $ \pt f_{\epsi,\ell}(s,K) $ in 
\eqref{Theorem.error.bound.m3.Eq.03.b} is straightforward because in \eqref{Theorem.error.bound.m3.Eq.03.b} the factor $\epsi$ compensates the integral.
The difficulty is that $ \pt z_{3\ell}(s,k^{(1)}) $ is not $\Ord{\epsi^2}$ in general.
We can only infer from \eqref{zdot.in.terms.of.uhat}, \eqref{Def.F.uhat} 
and Corollary~\ref{Corollary01.m3} that
\begin{align*}
\pt z_3(t) 
&=
\epsi \sum_{\#J=3} \TT_{3,\epsi}(t) \T\big(\uhat_{j_1},\uhat_{j_2},\uhat_{j_3}\big)(t)
\\
&=
\epsi \TT_{3,\epsi}(t) \T\big(\uhat_1,\uhat_1,\uhat_1\big)(t) + \Ord{\epsi^3}
\\
&=
\epsi \TT_{3,\epsi}(t) \T\big(\P_\epsi \uhat_1, \P_\epsi \uhat_1, \P_\epsi \uhat_1\big)(t) + \Ord{\epsi^2}.
\end{align*}
The $\ell$-th entry of the dominating part of $ \pt z_3(t,k^{(1)}) $ is thus
\begin{align*}
&
\epsi \Big[\TT_{3,\epsi}(t,k^{(1)}) \T\big(\P_\epsi \uhat_1,\P_\epsi \uhat_1,\P_\epsi \uhat_1\big)(t,k^{(1)})\Big]_\ell
\\
&=
\epsi \Big[\exp\big(\tfrac{\ii t}{\epsi}\Lambda_3(\epsi k^{(1)})\big) \Psi_{3}^*(\epsi k^{(1)}) 
\T\big(\P_\epsi \uhat_1,\P_\epsi \uhat_1,\P_\epsi \uhat_1\big)(t,k^{(1)})\Big]_\ell
\\
&=
\epsi \exp\big(\tfrac{\ii t}{\epsi}\lambda_{3\ell}(\epsi k^{(1)})\big) \psi_{3\ell}^*(\epsi k^{(1)}) 
\T\big(\P_\epsi \uhat_1,\P_\epsi \uhat_1,\P_\epsi \uhat_1\big)(t,k^{(1)})
\\
&=
\epsi \exp\big(\tfrac{\ii t}{\epsi}\lambda_{3\ell}(\epsi k^{(1)})\big) \phi_\epsi(t,k^{(1)})
\end{align*}
with the abbreviation
\begin{align*}
\phi_\epsi(t,k^{(1)}) = 
\psi_{3\ell}^*(\epsi k^{(1)}) 
\T\big(\P_\epsi \uhat_1,\P_\epsi \uhat_1,\P_\epsi \uhat_1\big)(t,k^{(1)}).
\end{align*}
All in all, it follows that
\begin{align}
\label{Theorem.error.bound.m3.Eq.04}
\pt f_{\epsi,\ell}(s,K)
 &=
  \epsi \exp\big(\tfrac{\ii s}{\epsi}\lambda_{3\ell}(\epsi k^{(1)})\big)\Phi_\epsi(s,K)
  + \Ord{\epsi^2}
\intertext{with}
  \nonumber
  \Phi_\epsi(s,K) &=
\frac{1}{(2\pi)^d} 
\Psi_5^*(\epsi k)
 T\left(
 \phi_\epsi(s,k^{(1)})
 \psi_{3\ell}(\epsi k^{(1)}),
 \P_\epsi\uhat_1(s,k^{(2)}),\P_\epsi\uhat_1(s,k^{(3)})\right).
\end{align}
Substituting the right-hand side of \eqref{Theorem.error.bound.m3.Eq.04} into \eqref{Theorem.error.bound.m3.Eq.03.b} yields
\begin{align*}
& C \epsi 
\Big| 
\intl_0^t 
\exp\left(\frac{\ii s}{\epsi}\Delta_\ell(\epsi k,\epsi k^{(1)}) \right)
\pt f_{\epsi,\ell}(s,K)
\; \dd s \Big|_2
\\
&=
C \epsi^2 
\Big| 
\intl_0^t 
\exp\left(\frac{\ii s}{\epsi}\Delta_\ell(\epsi k,\epsi k^{(1)})\right)
\exp\big(\tfrac{\ii s}{\epsi}\lambda_{3\ell}(\epsi k^{(1)})\big)
\Phi_\epsi(s,K)
\; \dd s \Big|_2 + \Ord{\epsi^2}
\\
&=
C \epsi^2 
\Big| 
\intl_0^t 
\exp\left(\frac{\ii s}{\epsi}
\Lambda_5(\epsi k) 
\right)
\Phi_\epsi(s,K)
\; \dd s \Big|_2 + \Ord{\epsi^2}
\end{align*}
because by definition
$ \Delta_\ell(\epsi k,\epsi k^{(1)}) =  \Lambda_5(\epsi k)-\lambda_{3\ell}(\epsi k^{(1)})I $; 
see~\eqref{Theorem.error.bound.m3.Def.Delta.ell}.
In order to show uniform boundedness of
\begin{align*}
\Big| 
\intl_0^t 
\exp\left(\frac{\ii s}{\epsi}
\Lambda_5(\epsi k) 
\right)
\Phi_\epsi(s,K)
\; \dd s \Big|_2 
=
\Big| 
\intl_0^t 
\exp\left(\frac{\ii s}{\epsi}\Lambda_5(0)\right)
\exp\left(\frac{\ii s}{\epsi}
\big[\Lambda_5(\epsi k)-\Lambda_5(0)\big]
\right)
\Phi_\epsi(s,K)
\; \dd s \Big|_2,
\end{align*}
we can use integration by parts again, because $ \Lambda_5(0) $ is regular by Assumption~\ref{Ass:Delta.lambda.j5} and the time derivative of 
$ \exp\left(\frac{\ii s}{\epsi} [\Lambda_5(\epsi k)-\Lambda_5(0)] \right) \Phi_\epsi(s,K)$ 
is uniformly bounded.
This completes the proof of Theorem~\ref{Theorem.error.bound.m3}.
\qed

% --------------------------------------------------------------------------------
\subsection{Numerical experiment and discussion} \label{Subsec:NumEx.02}
% --------------------------------------------------------------------------------

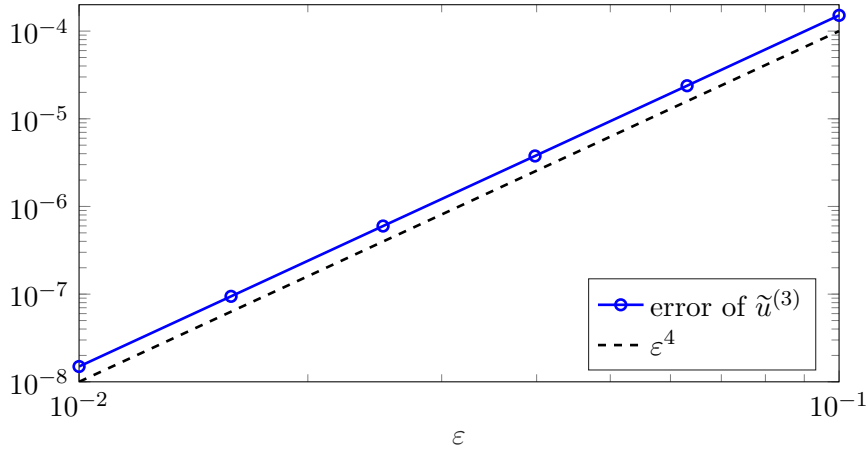
\begin{figure}[htb] 
\begin{center}
	% This file was created by matlab2tikz.
%
%The latest updates can be retrieved from
%  http://www.mathworks.com/matlabcentral/fileexchange/22022-matlab2tikz-matlab2tikz
%where you can also make suggestions and rate matlab2tikz.
%
\begin{tikzpicture}

\begin{axis}[%
width=10cm,
height=5cm,
at={(0cm,0cm)},
scale only axis,
xmode=log,
xmin=0.01,
xmax=0.1,
xminorticks=true,
xlabel style={font=\color{white!15!black}},
xlabel={$\varepsilon$},
ymode=log,
ymin=1e-08,
ymax=2*1e-04,
yminorticks=true,
ylabel style={font=\color{white!15!black}},
ylabel={},
axis background/.style={fill=white},
legend style={legend cell align=left, align=left, draw=white!15!black},
legend pos={south east}
]
\addplot [color=blue, line width=1pt, mark=o, mark options={solid, blue}]
  table[row sep=crcr]{%
0.01	1.49531512727563e-08\\
0.0158489319246111	9.44170007732836e-08\\
0.0251188643150958	5.98148556352207e-07\\
0.0398107170553497	3.76639399946033e-06\\
0.0630957344480193	2.38176453011274e-05\\
0.1	0.000151259697562467\\
};
\addlegendentry{error of $\uutilde^{(3)}$}

\addplot [color=black, dashed, line width=1pt]
  table[row sep=crcr]{%
0.01	1e-08\\
0.0158489319246111	6.30957344480193e-08\\
0.0251188643150958	3.98107170553497e-07\\
0.0398107170553497	2.51188643150958e-06\\
0.0630957344480193	1.58489319246111e-05\\
0.1 0.0001\\
};
\addlegendentry{$\varepsilon^4$}

\end{axis}
\end{tikzpicture}%
\begin{minipage}{12cm}
	\caption{Accuracy of $\uutilde^{(3)}$ for different values of $\epsi$. Parameters, data and discretizations are the same as in the numerical experiment described in Section~\ref{Subsec:NumEx.01}.}
  \label{Fig.error_m3}
\end{minipage}
\end{center}
\end{figure}

We have repeated the numerical experiment described in Section~\ref{Subsec:NumEx.01} with $m=3$ instead of $m=1$,
and with $\tstar=\tend=1$. Figure~\ref{Fig.error_m3} shows that in this example 
the numerical counterpart of the error 
$ \sup_{t\in[0,\tstar/\epsi]} \| \uu(t) - \uutilde^{(3)}(t) \|_{L^\infty} $ scales like $\epsi^4$,
which is better than what the error bound \eqref{Theorem.error.bound.m3assertion02} in Theorem~\ref{Theorem.error.bound.m3} predicts.
We believe, however, that this advantageous error behavior cannot be expected in general, and we briefly sketch the reasons.

If we want to improve \eqref{Theorem.error.bound.m3assertion02} in such a way that $\epsi^3$ is replaced by 
$\epsi^4$, then instead of \eqref{Theorem.error.bound.m3.basic} we have to prove that
\begin{align*}
 \sum_{|j| \in\{5,7,9\}}  \sum_{\#J=j} 
\Big\| \intl_0^t \TT_{j,\epsi}(s) \T\big(\uhat_{j_1},\uhat_{j_2},\uhat_{j_3} \big)(s)
\; \dd s \Big\|_{L^1} \leq C \epsi^3.
\end{align*}
As before, the critical indices are those where $|j|=5=\#J=|J|_1$, and we consider again
$j=5$ and $J=(3,1,1)$ as an example. 
Now instead of \eqref{main.difficulty.jmax3}, we have to show that
\begin{subequations}
\begin{align} 
\label{Theorem.error.bound.m3.basic.better.have2show1}
\Big\| \intl_0^t \TT_{5,\epsi}(s) \T\big(\uhat_{3},\P_\epsi^\perp \uhat_{1},\P_\epsi \uhat_{1} \big)(s) \; \dd s \Big\|_{L^1} 
\leq C \epsi^3
\\
\label{Theorem.error.bound.m3.basic.better.have2show2}
\text{and} \quad 
\Big\| \intl_0^t \TT_{5,\epsi}(s) \T\big(\uhat_{3},\P_\epsi \uhat_{1},\P_\epsi \uhat_{1} \big)(s) \; \dd s \Big\|_{L^1} 
\leq C \epsi^3.
\end{align} 
\end{subequations}
We will now explain why the first inequality \eqref{Theorem.error.bound.m3.basic.better.have2show1}
cannot be true in general.
Since $ \uhat_{3}=\Ord{\epsi^2} $ and 
$ \P_\epsi^\perp \uhat_{1} = \Ord{\epsi} $ by Corollary~\ref{Corollary01.m3}, the integrand is formally $\Ord{\epsi^3}$, but since $t\in[0,\tend/\epsi]$ we need one more factor of $\epsi$ to compensate the long integration interval.
By \eqref{Def.TT} and \eqref{Def.T.Fourier} the integral in 
\eqref{Theorem.error.bound.m3.basic.better.have2show1} reads
\begin{align*} 
& \intl_0^t \TT_{5,\epsi}(s) \T\big(\uhat_{3},\P_\epsi^\perp \uhat_{1},\P_\epsi \uhat_{1} \big)(s) \; \dd s
\\
&=
\frac{1}{(2\pi)^d} 
\intl_0^t
\exp\big(\tfrac{\ii s}{\epsi}\Lambda_5(\epsi k)\big)  \Psi_5^*(\epsi k)
\intl_{\#K = k} T\left(\uhat_3(s,k^{(1)}),\P_\epsi^\perp\uhat_1(s,k^{(2)}),\P_\epsi\uhat_1(s,k^{(3)})\right)
\; \dd K
\; \dd s.
\end{align*} 
After substituting \eqref{representation.Pperpuhat.1}, 
\eqref{Theorem.error.bound.m3.representation.uhat3}, and 
$ \P_\epsi\uhat_1(s,k^{(3)}) = \psi_{11}(\epsi k^{(3)})\ee^{-\ii s \lambda_{11}(\epsi k^{(3)})/\epsi}
z_{11}(s,\epsi k^{(3)}) $
we obtain
\begin{align*} 
& \intl_0^t \TT_{5,\epsi}(s) \T\big(\uhat_{3},\P_\epsi^\perp \uhat_{1},\P_\epsi \uhat_{1} \big)(s) \; \dd s
\\
&= 
\sum_{\ell_1=1}^n \sum_{\ell_2=2}^n \intl_{\#K = k}
\intl_0^t
\exp\big(\tfrac{\ii s}{\epsi}
\big[\Lambda_5(\epsi k)-\big(\lambda_{3\ell_1}(\epsi k^{(1)}) 
+ \lambda_{1\ell_2}(\epsi k^{(2)})
+ \lambda_{11}(\epsi k^{(3)})\big)I
\big]\big)
f_{\epsi,\ell_1,\ell_2}(s,K)
\, \dd s \, \dd K,
\end{align*} 
with a smooth function $ f_{\epsi,\ell_1,\ell_2} $. (Details do not matter at this point.)
% \begin{align*}
% f_{\epsi,\ell_1,\ell_2}(s,K)=
% \frac{1}{(2\pi)^d} 
% \Psi_5^*(\epsi k)
%  T\Big(z_{3\ell_1}(s,k^{(1)})\psi_{3\ell_1}(\epsi k^{(1)}),
%  z_{1\ell_2}(s,k^{(2)})\psi_{1\ell_2}(\epsi k^{(2)}),
%  &
%  \\
%  z_{11}(s,k^{(3)})\psi_{11}(\epsi k^{(3)})
%  \Big), &
% \end{align*}
% for $\#K = k$.
In order to generate an $\epsi$ via integration by parts, we need that the diagonal matrix
\begin{align*}
\Lambda_5(0)-\big(\lambda_{3\ell_1}(0) + \lambda_{1\ell_2}(0) + \lambda_{11}(0)\big)I
\end{align*}
is regular. Since $\lambda_{11}(0) = 0$, this is equivalent to the condition
\begin{align}
\label{Nonresonance.three.terms}
\lambda_{5\ell}(0)-\lambda_{3\ell_1}(0) - \lambda_{1\ell_2}(0) \not= 0
\qquad
\text{for all } \ell, \ell_1, \ell_2 \in \{1, \ldots, n\}, \quad \ell_2\not=1.
\end{align}
This is a non-resonance condition similar to what we have assumed in 
Assumptions~\ref{Ass:Delta.lambda.j3} and \ref{Ass:Delta.lambda.j5}, but now with three terms. 
In contrast to those assumptions, however,
\eqref{Nonresonance.three.terms} is \emph{not} true in case of the Klein--Gordon system with $d>1$, nor for the Maxwell--Lorentz system, as we will show now.
In these applications, the eigenvalues $ \w_\ell(\beta) $ of $ \L(0,\beta)=A(\beta)-\ii E $ 
have the following properties:
\begin{itemize}
\item[(P1)]
The largest eigenvalue $ \w_1(\beta) $ 
is related to the smallest eigenvalue $ \w_n(\beta) $ by $ \w_n(\beta)=-\w_1(\beta) $. 
\item[(P2)]
 $ \L(0,\beta) $ has at least one vanishing eigenvalue, i.e.\ there is an index $\ell_\bullet$ with $ 1\not=\ell_\bullet\not=n$ and
  $ \w_{\ell_\bullet}(\beta) = 0 $ for all $\beta$.
\end{itemize}
Recall that $\w=\w(\kappa)$ is an eigenvalue of $ \L(0,\kappa)=A(\kappa)-\ii E$ (cf.~\eqref{Def.w}), 
and suppose that we have chosen $\w=\w_1(\kappa)$. 
By definition, the eigenvalues of 
\begin{align*}
\L_j(0) = \L(j\w, j\kappa) = -j\w I + \L(0,j\kappa)
\end{align*}
are $\lambda_{j\ell}(0)=-j\w + \w_\ell(j\kappa)$.
If we choose $\ell=\ell_1=\ell_\bullet$ and $\ell_2=n$ in \eqref{Nonresonance.three.terms}, then we obtain
\begin{align*}
\lambda_{5\ell_\bullet}(0)-\lambda_{3\ell_\bullet}(0) - \lambda_{1n}(0) 
&= (-5\w + \w_{\ell_\bullet}(5\kappa)) - (-3\w + \w_{\ell_\bullet}(3\kappa)) - (-\w + \w_n(\kappa))
\\
&= -5\w + 0 + 3\w - 0 + \w - (-\w_1(\kappa))
\\
&= -\w + \w_1(\kappa) = 0,
\end{align*}
which shows that the non-resonance condition \eqref{Nonresonance.three.terms} is not true. 
This is only one counterexample among many others.
The corresponding non-oscillatory terms in the integrand cause contributions of $\Ord{t\epsi^3}$, 
which eventually leads to a contribution of $\Ord{\epsi^2}$ instead of $\Ord{\epsi^3}$
on the left-hand side of 
\eqref{Theorem.error.bound.m3.basic.better.have2show1}. 
Similar resonance problems appear also in the integral in 
\eqref{Theorem.error.bound.m3.basic.better.have2show2}, such that this inequality 
cannot be true for the applications mentioned above.

A noteworthy exception is the Klein--Gordon system in \emph{one} space dimension ($d=1$, $n=2$), which we have used in our numerical experiments. 
Here, the two eigenvalues of the matrix
$ \L(0,\kappa)=A(\kappa)-\ii E \in \IC^{2\times 2}$ are $ w=w_1(\kappa) = \sqrt{\kappa^2 + \gamma^2} $
and $ w_2(\kappa)=-w_1(\kappa)$,
as we have mentioned in Section~\ref{Subsec:NumEx.01}.
These eigenvalues have property (P1), but not property (P2), such that the counterexample does not apply.
We conjecture that in this special case, one could indeed prove that 
\eqref{Theorem.error.bound.m3assertion02} even holds with $\epsi^4$ instead of $\epsi^3$ on the right-hand side, which is the behavior observed in Figure~\ref{Fig.error_m3}.

This discussion raises the question if the convergence behavior predicted by Theorem~\ref{Theorem.error.bound.m3} could be observed in a numerical example with a \emph{two-dimensional} Klein--Gordon equation, because then the eigenvalues have also the property (P2). 
The problem is that in order to test the accuracy of the approximation $\uu \approx \uutilde^{(3)}$,
the PDEs \eqref{PDE.uu} and \eqref{PDE.mfe} have to be solved numerically with such a high precision that the numerical error is negligible compared to the analytical error.
But approximating $\uu$ with sufficiently high precision by applying a standard method to \eqref{PDE.uu} was already hopeless in one space dimension (cf.\ Section~\ref{Subsec:NumEx.01}), and computing
a reference solution via \eqref{Ansatz} and \eqref{PDE.mfe} with $m=5$ was already extremely expensive in the one-dimensional case, because the functions $u_j$ still oscillate in time. For these reasons, we were not able to produce a reliable numerical example in two space dimensions.

The approach to approximate the solution $\uu$ of \eqref{PDE.uu} via \eqref{Ansatz} and \eqref{PDE.mfe} has the advantage that the coefficient functions $ u_j $ do not oscillate in \emph{space}. This gives us the possibility to use a space discretization where the number of grid points depends only on the regularity of $ u_j(t,\cdot) $, but not on $1/\epsi$. 
To realize the full potential of this approach, however, it is important to develop tailor-made time integrators for \eqref{PDE.mfe},
which use non-standard techniques to handle the oscillations in \emph{time}, and which are far more efficient than traditional schemes such as the splitting method used in our numerical examples. 
In a joint work with Johanna M\"odl (KIT), the second author has recently constructed and analyzed such a tailor-made time integrator. This result will be reported elsewhere.

\bibliographystyle{abbrv}
\bibliography{hfwp}

\begin{thebibliography}{10}

\bibitem{alterman-rauch:00}
D.~Alterman and J.~Rauch.
\newblock Diffractive short pulse asymptotics for nonlinear wave equations.
\newblock {\em Phys. Lett. A}, 264(5):390--395, 2000.

\bibitem{alterman-rauch:03}
D.~Alterman and J.~Rauch.
\newblock Diffractive nonlinear geometric optics for short pulses.
\newblock {\em SIAM J. Math. Anal.}, 34(6):1477--1502, 2003.

\bibitem{barrailh-lannes:02}
K.~Barrailh and D.~Lannes.
\newblock A general framework for diffractive optics and its applications to
  lasers with large spectrums and short pulses.
\newblock {\em SIAM J. Math. Anal.}, 34(3):636--674, 2002.

\bibitem{baumstark:22}
J.~Baumstark.
\newblock {\em High-frequency wave-propagation: error analysis for analytical
  and numerical approximations}.
\newblock PhD thesis, Karlsruhe Institute of Technology (KIT), jul 2022.

\bibitem{baumstark-jahnke-2023}
J.~Baumstark and T.~Jahnke.
\newblock Approximation of high-frequency wave propagation in dispersive media.
\newblock {\em SIAM J. Math. Anal.}, 55(2):1214--1245, 2023.

\bibitem{baumstark-jahnke-lubich-2024}
J.~Baumstark, T.~Jahnke, and C.~Lubich.
\newblock Polarized high-frequency wave propagation beyond the nonlinear
  {S}chr\"{o}dinger approximation.
\newblock {\em SIAM J. Math. Anal.}, 56(1):454--473, 2024.

\bibitem{chung-jones-schaefer-wayne:05}
Y.~Chung, C.~K. R.~T. Jones, T.~Sch\"{a}fer, and C.~E. Wayne.
\newblock Ultra-short pulses in linear and nonlinear media.
\newblock {\em Nonlinearity}, 18(3):1351--1374, 2005.

\bibitem{colin-lannes:09}
M.~Colin and D.~Lannes.
\newblock Short pulses approximations in dispersive media.
\newblock {\em SIAM J. Math. Anal.}, 41(2):708--732, 2009.

\bibitem{colin:02}
T.~Colin.
\newblock Rigorous derivation of the nonlinear {S}chr\"{o}dinger equation and
  {D}avey-{S}tewartson systems from quadratic hyperbolic systems.
\newblock {\em Asymptot. Anal.}, 31(1):69--91, 2002.

\bibitem{colin-gallica-laurioux:05}
T.~Colin, G.~Gallice, and K.~Laurioux.
\newblock Intermediate models in nonlinear optics.
\newblock {\em SIAM J. Math. Anal.}, 36(5):1664--1688, 2005.

\bibitem{donnat-joly-metivier-rauch:96}
P.~Donnat, J.-L. Joly, G.~Metivier, and J.~Rauch.
\newblock Diffractive nonlinear geometric optics.
\newblock In {\em S\'{e}minaire sur les \'{E}quations aux {D}\'{e}riv\'{e}es
  {P}artielles, 1995--1996}, S\'{e}min. \'{E}qu. D\'{e}riv. Partielles, pages
  Exp. No. XVII, 25. \'{E}cole Polytech., Palaiseau, 1996.

\bibitem{donnat-rauch:97b}
P.~Donnat and J.~Rauch.
\newblock Dispersive nonlinear geometric optics.
\newblock {\em J. Math. Phys.}, 38(3):1484--1523, 1997.

\bibitem{donnat-rauch:97a}
P.~Donnat and J.~Rauch.
\newblock Modeling the dispersion of light.
\newblock In {\em Singularities and oscillations ({M}inneapolis, {MN},
  1994/1995)}, volume~91 of {\em IMA Vol. Math. Appl.}, pages 17--35. Springer,
  New York, 1997.

\bibitem{joly-metivier-rauch:93}
J.-L. Joly, G.~M\'{e}tivier, and J.~Rauch.
\newblock Generic rigorous asymptotic expansions for weakly nonlinear
  multidimensional oscillatory waves.
\newblock {\em Duke Math. J.}, 70(2):373--404, 1993.

\bibitem{joly-metivier-rauch:96}
J.~L. Joly, G.~Metivier, and J.~Rauch.
\newblock Global solvability of the anharmonic oscillator model from nonlinear
  optics.
\newblock {\em SIAM J. Math. Anal.}, 27(4):905--913, 1996.

\bibitem{joly-metivier-rauch:98}
J.-L. Joly, G.~Metivier, and J.~Rauch.
\newblock Diffractive nonlinear geometric optics with rectification.
\newblock {\em Indiana Univ. Math. J.}, 47(4):1167--1241, 1998.

\bibitem{joly-metivier-rauch:00}
J.-L. Joly, G.~Metivier, and J.~Rauch.
\newblock Transparent nonlinear geometric optics and {M}axwell-{B}loch
  equations.
\newblock {\em J. Differential Equations}, 166(1):175--250, 2000.

\bibitem{kirrmann-schneider-mielke:92}
P.~Kirrmann, G.~Schneider, and A.~Mielke.
\newblock The validity of modulation equations for extended systems with cubic
  nonlinearities.
\newblock {\em Proc. Roy. Soc. Edinburgh Sect. A}, 122(1-2):85--91, 1992.

\bibitem{lannes:98}
D.~Lannes.
\newblock Dispersive effects for nonlinear geometrical optics with
  rectification.
\newblock {\em Asymptot. Anal.}, 18(1-2):111--146, 1998.

\bibitem{lannes:11}
D.~Lannes.
\newblock High-frequency nonlinear optics: from the nonlinear {S}chr\"{o}dinger
  approximation to ultrashort-pulses equations.
\newblock {\em Proc. Roy. Soc. Edinburgh Sect. A}, 141(2):253--286, 2011.

\bibitem{rauch-book:12}
J.~Rauch.
\newblock {\em Hyperbolic partial differential equations and geometric optics},
  volume 133 of {\em Graduate Studies in Mathematics}.
\newblock American Mathematical Society, Providence, RI, 2012.

\bibitem{schneider-uecker:17}
G.~Schneider and H.~Uecker.
\newblock {\em Nonlinear {PDE}s}, volume 182 of {\em Graduate Studies in
  Mathematics}.
\newblock American Mathematical Society, Providence, RI, 2017.

\end{thebibliography}

\end{document}